\documentclass [12pt]{article}      
\usepackage{amsmath} 
\usepackage{amsthm}
\usepackage{latexsym}
\usepackage{pb-diagram}
\usepackage{amssymb}

\usepackage{color}

\usepackage{graphicx}
\usepackage{inputenc}

\usepackage{hyperref}
\hypersetup{
    colorlinks,
    citecolor=blue,
    filecolor=blue,
    linkcolor=blue,
    urlcolor=blue
}
\hypersetup{linktocpage}

\usepackage{todonotes}

\newtheorem{theorem}{Theorem}
\newtheorem*{theorem*}{Theorem}
\newtheorem{corollary}[theorem]{Corollary}
\newtheorem{lemma}[theorem]{Lemma}
\newtheorem{definition}[theorem]{Definition}
\newtheorem{prop}[theorem]{Proposition}

\newtheorem{ex}[theorem]{Example}

\newtheorem{remark}[theorem]{Remark}
\newtheorem{claim}[theorem]{Claim}

\newcommand{\cat}{^\frown}
\newcommand{\rest}{\ensuremath{\upharpoonright}}

\newcommand{\mcd}{{\ensuremath{\mathcal D}}}

\newcommand{\xbmt}{\ensuremath{(X,\mcb,\mu,T)}}
\newcommand{\ycns}{\ensuremath{(Y,\mcc,\nu,S)}}

\renewcommand{\qed}{{\nopagebreak \hfill $\dashv$ \par\bigskip}}

\newcommand{\pf}{{\par\noindent{$\vdash$\ \ \ }}}

\newcommand{\la}{\langle}
\newcommand{\ra}{\rangle}

\newcommand{\mca}{\ensuremath{\mathcal A}}

\newcommand{\mcr}{\ensuremath{\mathcal R}}
\newcommand{\dbar}{\ensuremath{\bar{d}}}

	\newcommand{\poZ}{\mathbb Z}
	
	\newcommand{\mcw}{\mathcal W}
	\newcommand{\nn}{{\mathbb N}}

	\newcommand{\bfni}[1]{\noindent {{\bf{#1}}}}
	
	\newcommand{\inv}{{^{-1}}}

	\newcommand{\rev}[1]{\mathop{\rm rev}({#1})}

	\newcommand{\bk}{{\mathbb K}}
	\newcommand{\bm}{{\mathbb M}}

\newcommand{\qn}{{q_n}}
\newcommand{\qnpo}{{q_{n+1}}}

\newcommand{\kn}{{k_n}}

\newcommand{\mcc}{{\mathcal C}}

\newcommand{\zoo}{{[0,1)}}

\newcommand{\mck}{{\mathcal K}}

\newcommand{\mcb}{\mathcal B}

\newcommand{\wa}{w^\alpha}

\newcommand{\sz}{\Sigma^\poZ}

\newcommand{\bt}{\mathbb T}

\newcommand{\boundary}{\partial}

\newcommand{\bl}{\mathbb L}

\newcommand{\mco}{\mathcal O}
\newcommand{\mcf}{\mathcal F}

\newcommand{\mcu}{\mathcal U}
\newcommand{\mcv}{\mathcal V}

 \title{From Odometers to Circular Systems:\\ A Global Structure Theorem}
 \author{Matthew Foreman, Benjamin Weiss}
 \begin{document}
 \maketitle 
 
 \begin{abstract}{The main result of this paper is that two large collections of ergodic measure preserving systems, the \emph{Odometer Based} and the \emph{Circular Systems} have the same global structure with respect to joinings. The classes are canonically isomorphic by a continuous map that takes factor maps to factor maps, measure-isomorphisms to measure-isomorphisms, weakly mixing extensions to weakly mixing extensions and compact extensions to compact extensions. The first class includes all finite entropy ergodic transformations with an odometer factor. By results in \cite{prequel}, the second class contains all transformations realizable as diffeomorphisms using the untwisted Anosov-Katok method. An application of the main result will appear in a forthcoming paper that shows that the diffeomorphisms of the torus are inherently unclassifiable up to measure-isomorphism. Other consequences include the existence measure distal diffeomorphisms of arbitrary countable distal height.}

\end{abstract}

 \tableofcontents

 \section{Introduction}
 The isomorphism problem in ergodic theory was formulated by von Neumann in 1932 in his pioneering paper \cite{vN}.
  Simply put it asks to determine when two measure preserving actions  are isomorphic, in the sense that there is
  a measure isomorphism between the underlying measure space that intertwines the actions.
  It has been solved completely only for some special classes of transformations.
   Halmos and von Neumann \cite{HvN} used the unitary operators defined by Koopman to completely characterize ergodic measure preserving transformations with pure point spectrum, these transformations can be concretely realized (in a Borel way) as translations on compact groups.
Another notable success was the use of the Kolmogorov entropy to distinguish between measure preserving systems. Ornstein's work showed that entropy completely classifies a large class of highly random systems, such as independent processes, mixing Markov chains and certain smooth systems such as geodesic flows on surfaces of negative curvature.

Closely related to the isomorphism problem is the study of structural properties of measure preserving systems. These including mixing properties and compactness. A famous example is the Furstenberg-Zimmer structure theorem for ergodic measure preserving transformations, which characterizes every ergodic transformation as an inverse limit system of compact extensions followed by a weakly mixing extension. This result is fundamental for studying recurrence properties of measure preserving systems and the related proofs of Szemeredi-type combinatorial theorems (\cite{FuBook}).

In this paper we present a new phenomenon, \emph{Global Structure Theory}. Most structure theorems in ergodic theory consider a single transformation \emph{in vitro}. The approach here is study whole, intact ecosystems of transformations with their inherent relationships. 

Our main result shows that two large collections of measure preserving transformations have exactly the same structure  with respect to factors and isomorphisms (and more generally, joinings). More concretely, define the \emph{odometer based} transformations to be those finite entropy transformations that contain a non-trivial odometer factor. Spectrally, this is equivalent to the associated unitary operator having infinitely many finite period eigenvalues. To each odometer, we can associate a class of symbolic systems, the \emph{circular systems}. 
In \cite{prequel}, it is shown that the circular systems coincide exactly with the ergodic transformations realizable as diffeomorphisms of the torus
using the untwisted method of Approximation-by-Conjugacy, due to  Anosov-Katok (\cite{AK_original}).

We can make two categories by taking the objects  to be these two classes of systems and by taking morphisms to be factor maps (or more generally joinings) that preserve the underlying timing structure. The main theorem of this paper says that these two categories are isomorphic by a map that takes measure-isomorphisms to measure-isomorphisms, weakly mixing extensions to weakly mixing extensions and compact extensions to compact extensions. It follows that it takes distal towers to distal towers.  Moreover the map preserves the simplex of non-atomic invariant measures, takes rank one transformations to rank one transformations and much more. (This will be discussed further in the forthcoming \cite{part4}.)
In other words the global structure of these two categories is identical.

We can get more detail by considering systems based on a fixed odometer map and circular systems based on that odometer map and an arbitrary fast growing coefficient sequence. Doing so gives us collections of pairwise isomorphic categories that can be amalgamated to yield the statement above. The main theorem is framed in this more granular setting.

Our result might be a mere curiosity, were it not for an application which we now describe.
\medskip

Foreshadowed by a remarkable early result by Feldman \cite{feldman}, in the late 1990's a different type of result began to appear: \emph{anti-classification} results that demonstrate in a rigorous way that classification is not possible. This type of theorem requires a precise definition of what a classification is. Informally a classification is a \emph{method} of determining isomorphism between transformations perhaps by computing (in a liberal sense) other invariants for which equivalence is easy to determine.

The key words here are \emph{method} and \emph{computing}. For negative theorems, the more liberal a notion one takes the stronger the theorem. One natural notion is the Borel/non-Borel distinction. Saying a set $X$ or function $f$ is Borel is a loose way of saying that membership in $X$ or the computation of $f$ can be done using a countable (possibly transfinite) protocol whose basic input is membership in open sets. Say that $X$ or $f$ is \emph{not} Borel is saying that determining membership in $X$ or computing $f$ cannot be done with any amount of countable resources.

In the context of classification problems, saying that an equivalence relation $E$ on a space $X$ is \emph{not} Borel is
saying that there is no countable amount of information and no countable transfinite protocol for determining, for arbitrary
$x,y\in X$ whether $xEy$. \emph{Any} such method must inherently use uncountable
resources.\footnote{Many well known classification theorems have as  immediate corollaries that the resulting equivalence
relation is Borel. An example of this is the Spectral Theorem, which has a consequence that the relation of Unitary
Conjugacy for normal operators is a Borel equivalence relation.}

In considering the isomorphism relation as a collection $\mathcal I$ of pairs $(S,T)$ of measure preserving transformations, Hjorth showed that $\mathcal I$ is not a Borel set. However the pairs of transformations he used to demonstrate this were inherently non-ergodic\footnote{The ergodic components of the pairs were rotations of the circle.}, leaving open the essential problem:

\begin{center}
Is isomorphism of ergodic measure preserving transformations Borel?
\end{center}

This question was answered by Foreman, Rudolph and Weiss in \cite{FRW}, where they gave a negative answer. This answer can be interpreted as saying that determining isomorphism between ergodic transformations is inaccessible to countable methods that use countable amounts of information.

In the same foundational paper from 1932 where von Neumann formulated the isomorphism problem he expressed the likelihood that any abstract measure preserving transformation is isomorphic to a continuous measure preserving transformation and
perhaps even to a differentiable one. This brief remark eventually gave rise to one of the  outstanding problems in smooth dynamics, namely:

\begin{center}

 Does every ergodic MPT have a smooth model?
\end{center}

  By a smooth model is meant an isomorphic copy of the transformation which is given by smooth diffeomorphism of a compact manifold preserving a measure equivalent to the volume element.
Soon after entropy was introduced, A. G. Kushnirenko showed that such a diffeomorphism must have finite entropy, and up to now this is the only
restriction that is known. 

This paper is the second in a series of papers whose original purpose was to show that the variety of ergodic transformations that have  smooth models is rich enough so
that the abstract isomorphism relation, when restricted to these smooth systems, is as complicated as it is in general. We show this to be the case even when restricting to
diffeomorphisms of the 2-torus that preserve Lebesgue measure this is the case. In the third paper we will complete  the proof of the following theorem:
\bigskip

\begin{theorem*}[Anti-classification of Diffeomorphisms]

If $M$ is either the torus $\bt^2$, the disk $D$ or the annulus then the measure-isomorphism relation among pairs $(S,T)$ of measure preserving $C^\infty$-diffeomorphisms of $M$ is not a Borel set with respect to the $C^\infty$-topology.
\end{theorem*}
   It was natural for us to try to adapt our earlier work to establish this result. However we were faced at first with the following difficulty.
   The transformations built in \cite{FRW}  were based on odometers (in the sense that the Kronecker factor was an odometer). 
   It is a well known open problem whether it is possible to have
 any smooth transformation on a compact manifold that has a non-trivial odometer factor. Thus proving the anti-classification theorem in the smooth context required constructing a different collection of hard-to-classify transformations and then showing that this collection could be realized smoothly. This is our application of the main result of this paper.

The paper (\cite{prequel})  constructed a new collection of systems, the \emph{Circular Systems}, which are defined as symbolic systems constructed using the \emph{Circular Operator}, a formal operation on words. The main result in \cite{prequel} has as a consequence that uniform circular systems can  be realized as smooth models using the method developed by Anosov and Katok.

The primary theorem of this paper allows us to transfer the general isomorphism structure for odometer based systems to the isomorphism structure for circular systems, at least up to automorphisms of the underlying odometer or rotation. Namely there remains the issue of preserving the timing mechanism. In the forthcoming \cite{part3} it is shown how to construct odometers so that  for the resulting circular systems, up to a small correction factor, all isomorphisms preserve the underlying timing structure. This allows us to conclude  the proof of the anti-classification theorem for diffeomorphisms.

\medskip
Here is a more concrete description of the results in the paper.
 In the present paper we are concerned with  the entire class $\mathcal{OB}$
  of systems based on a fixed odometer and the relations between them. The odometer is
  determined by a sequence of positive integers greater than one, $\la k_n : n \in \mathbb{N}\ra$.  The  the circular operator is determined by an additional
  sequence of integers $\la l_n: n \in \mathbb{N}\ra$. For this paper, the sequence of $l_n$'s can be arbitrary subject to the requirement  that 
  $\sum 1/l_n<\infty$. However for realizing circular systems as diffeomorphisms there is a fixed growth rate, determined by the size of the alphabet of the odometer based system and $\la k_n:n\in\nn\ra$, that the sequence of $l_n$'s must eventually exceed.

  We describe $\mathcal{OB}$ symbolically here, but show in a forthcoming paper that $\mathcal{OB}$ consists of  representations of arbitrary ergodic systems with finite entropy 
  that have the specific odometer as a factor. In the language of ``cutting and stacking" constructions these are those constructions where no
  spacers are introduced. We fix $\la l_n:n\in\nn\ra$, and hence a sequence of circular operators.
  Applying these  to each of the elements of $\mathcal{OB}$ we obtain a second class, 
  $\mathcal{CB}$, of circular systems. This class consists of some of the extensions of
  a fixed irrational rotation which is determined by the circular operator.  As remarked above, for suitably chosen coefficient sequences, this class can be characterized as those transformations realizable as diffeomorphisms using the Anosov-Katok technique. We consider the two classes as categories
  where the morphisms are  graph joinings which are either  the identity of the base or reverse it. These are called \emph{synchronous} and
  \emph{anti-synchronous} joinings respectively. Our main theorem then takes the form:
  
  \begin{theorem}\label{main in intro}
  For a fixed circular coefficient sequence $\la k_n, l_n: n\in\nn\ra$ the categories $\mathcal{OB}$ and  $\mathcal{CB}$ are isomorphic by a functor $\mcf$ that takes synchronous joinings to synchronous joinings, anti-synchronous joinings to anti-synchronous joinings,
isomorphisms to isomorphisms and weakly mixing extensions to weakly mixing extensions.\footnote{E. Glasner showed that the functor takes compact extensions to compact extensions.}   
  \end{theorem}

It is natural to extend the collections of morphisms of $\mathcal{OB}$ and $\mathcal{CB}$ to general synchronous and non-synchronous joinings.   Because the ergodic joinings are not closed under composition, in extending Theorem \ref{main in intro} one is forced to consider  at least \emph{some} non-ergodic joinings. At the end of the paper we discuss how to extend Theorem \ref{main in intro} to expanded categories that have as morphisms arbitrary synchronous and anti-synchronous joinings.  This involve expanding our analysis of generic sequences to non-ergodic joinings. We also describe some detailed analysis of the combinatorics behind the isomorphism $\mcf$.

   We have provided a detailed table of contents which enumerates the contents of the paper. Here is a brief summary. Much of the section following this one is standard, with the exception \S 2.6, which is exposes generic sequences for transformations and extends that notion to joinings. In \S3,  the reader will find an explanation of our two categories and  a proof that circular systems contain a canonical rotation factor. Section 4 is primarily concerned with defining a map $\natural$ that is a symbolic analogue of complex conjugation on the unit circle. In    sections 5 and 6 the mapping $\mcf$ is defined on morphisms, while $\S{7}$ contains the proof of the main result. In $\S{8}$ there is a more detailed analysis of of the dynamical
   properties of our mapping $\mcf$ which may prove useful in the future, and in the final section we collect some problems that are left open.

\subsection{Acknowledgements}

This work was inspired by the pioneering work of our co-author Dan Rudolph, who passed away before this portion of the grand project was undertaken. We owe an inestimable debt to J.P. Thouvenot who suggested using the Anosov-Katok technique to produce our badly behaved transformations rather than directly attacking the ``odometer obstacle." We would like to thank E. Glasner for showing that $\mcf$ preserves compact extensions.  Finally the first author would like to thank Christian  Rosendal for  asking very useful questions about how general our results were. 
 
 \section{Preliminaries}\label{preliminaries}
This section establishes some of the conventions we follow in this paper. There are many sources of background information on this including any standard text or \cite{walters}, \cite{Peterson}. A small portion of the material in this section was presented in \cite{prequel}, but is repeated here in an attempt to be self-contained. The reader is referred to \cite{prequel} for any missing definitions.

\subsection{Measure Spaces}\label{abstract measure spaces} We will call separable non-atomic probability spaces \emph{measure spaces} and denote them $(X,\mcb, \mu)$ where $\mcb$ is the Boolean algebra of measurable subsets of $X$ and $\mu$ is a countably additive, non-atomic measure defined on $\mcb$.\footnote{We will occasionally  make an exception to this by calling discrete probability measures on a finite set \emph{measures}; we hope that context makes the difference clear.} We will often identify two  members of $\mcb$ that differ
by a set of $\mu$-measure $0$ and seldom distinguish between $\mcb$ and the $\sigma$-algebra of classes of measurable sets modulo measure zero unless we are making a pointwise definition and need to claim it is well defined on equivalence classes.

We will frequently use  without explicit mention the 
Maharam-von Neumann result that every standard measure space is isomorphic to $([0,1],\mcb,\lambda)$ where $\lambda$ is Lebesgue measure and $\mcb$ is the algebra of Lebesgue measurable sets.

If $(X, \mcb, \mu)$ and $(Y, \mcc, \nu)$ are measure spaces, an isomorphism between $X$ and $Y$ is a bijection $\phi:X\to Y$ such that $\phi$ is measure preserving and both $\phi$ and $\phi^{-1}$ are measurable. We will ignore sets of measure zero when discussing isomorphisms; i.e.  we allow  the domain and range of $\phi$ to be  subsets of $X$  and $Y$ (resp.) of measure one. 
A measure preserving system is an object $\xbmt$ where $T:X\to X$ is a measure isomorphism. A \emph{factor map} 
between two measure preserving systems $\xbmt$ and $\ycns$ is a measurable, measure preserving  function 
$\phi:X\to Y$ such that $S\circ\phi=\phi\circ T$. A factor map is an \emph{isomorphism} or \emph{conjugacy} between systems iff $\phi$ is a 
measure isomorphism. Following common practice, we will use the word \emph{conjugacy} interchangeably with \emph{isomorphism} in this context.

For a fixed measure space $(X,\mu)$ we can consider the collection of measure preserving transformations $T:X\to X$. These form a group that can be endowed with a Polish topology that has basic open sets described as follows. We fix a finite measurable partition $\mca$ of $X$ and an $\epsilon>0$ and take as a neighborhood of $T$
\[\mathcal N(T,\mca,\epsilon)=_{def}\{S:\sum_{a\in \mca}\mu(Ta\Delta Sa)<\epsilon\}.\]
 Details about this topology can be found in many sources including \cite{halmos}, \cite{walters}.
 \subsection{Joinings}\label{joinings}

We remind the readers of the definitions. Extensive treatments of joinings   can be found in \cite{glasbook} or \cite{DansBook}.  All of the definitions and basic results about joinings necessary for this paper occur in Chapter 6 of the latter reference.

\begin{definition}
 A \emph{joining} between two measure preserving systems $\xbmt$ and $\ycns$ is a measure $\rho$ on $X\times Y$ defined on the product $\sigma$-algebra $\mcb\otimes\mcc$  such that 
 \begin{enumerate}
 \item $\rho$ is $T\times S$ invariant,
 \item for each set $B\in \mcb$, $\rho(B\times Y)=\mu(B)$,
 \item for each set $C\in \mcc$, $\rho(X\times C)=\nu(C)$.
 \end{enumerate}
\end{definition}

The graphs of factor maps provide natural examples of joinings. We characterize these with a definition.

\begin{definition}
A joining $\rho$ is a \emph{graph joining} between $X$ and $Y$ if and only if for all $C\in \mcc$ and all $\epsilon>0$,  there is a $B\in \mcb$ such that 
\[\rho((B\times Y)\Delta (X\times C))<\epsilon.\]

A joining $\rho$ between $\xbmt$ and $\ycns$ is an \emph{invertible graph joining} if and only for all $B\in \mcb$ there is a $C\in \mcc$ such that 
\begin{equation}\label{graph capture}
\rho((B\times Y)\Delta (X\times C))=0
\end{equation}
and vice versa: for all $C\in \mcc$, there is a $B\in\mcb$ such that equation \ref{graph capture} holds.
\end{definition}

Here are some standard facts (see \cite{glasbook}):

\begin{prop}\label{no proof} Let $\mathbb X=\xbmt$ and $\mathbb Y=\ycns$. Then
\begin{enumerate}
\item There is a canonical one-to-one correspondence between the collection of graph joinings of $\mathbb X$ and $\mathbb Y$ and the collection of factor maps from $X$ to $Y$. A graph joining concentrates on the graph of the factor map. We can represent the graph joining corresponding to a measure preserving map $\phi:X\to Y$ by
\[\rho_\phi=\int (\delta_x\times \delta_{\phi(x)}) d\mu(x).\]
\item There is a canonical one-to-one correspondence between the collection of invertible graph joinings of $\mathbb X$ and $\mathbb Y$ and the collection of conjugacies between $\mathbb X$ and $\mathbb Y$.

\item \label{working graph} Suppose that $\mcb'\subseteq \mcb$ and $\mcc'\subseteq \mcc$ are Boolean algebras that generate $\mcb$ and $\mcc$ respectively as $\sigma$-algebras. Let $\rho$ be a joining of $\mathbb X$ with $\mathbb Y$ such that for all $\epsilon>0$ and all $C\in \mcc'$ there  are $B_1, \dots B_n\in \mcb'$ such that we have  $\rho(\bigcup_i(B_i\times Y)\Delta(X\times C))<\epsilon$, then $\rho$ is a graph joining.

\end{enumerate}
\end{prop}
\smallskip

We note that perhaps a more proper term for an invertible graph joining is the earlier usage \emph{diagonal joining}. In view of the results of this section we will often be careless  and say that $\rho$ \emph{is a factor map} or $\rho$ \emph{is a conjugacy}/\emph{isomorphism} to mean that $\rho $ is a graph joining or $\rho$ is an invertible graph joining.

\smallskip

To each joining $\rho$ of $\mathbb X$ and $\mathbb Y$ we can associate its adjoint $\rho^*$, the joining of $\mathbb Y$ with $\mathbb X$ defined for $B\in \mcb$ and $C\in \mcc$ as:
\[\rho^*(C\times B)=\rho(B\times C).\]
If $\rho$ is a graph joining corresponding to a factor map $\pi:X\to Y$, then $\rho^*$  concentrates on $\{(y,x):\pi(x)=y\}$. 

The following is immediate:
\begin{prop}\label{invertible symmetry}
$\rho$ is an invertible graph joining if and only if both $\rho$ and $\rho^*$ are graph joinings. 
\end{prop}
Thus  we can apply Proposition \ref{no proof}, item \ref{working graph} to both $\rho$ and $\rho^*$ to get a criterion for being the joining associated with a conjugacy.

\medskip

\bfni{A potential source of confusion.} Proposition \ref{no proof} allows us to identify graph joinings with factor maps and invertible graph joinings with conjugacies. These joinings are always \emph{ergodic} as joinings.  However, there are non-ergodic conjugacies between ergodic measure preserving transformations. More explicitly: there are ergodic systems $(X,T)$ and $(X,S)$ and non-ergodic isomorphisms $\phi:(X,T)\to (X,S)$.\footnote{The second author has given examples of of isomorphic ergodic transformations where every conjugacy is non-ergodic.}  The associated joining $\rho_\phi$ is, however, ergodic as a $T\times S$-invariant measure. 

\medskip
Let $(X,\mu), (Y,\nu)$ and $(Z,\tilde{\mu})$ be measure spaces and 
$\pi_X:X\to Y$ and $\pi_Z:Z\to Y$ be factor maps. We can define  a canonical 
joining of  $\mathbb X$ and $\mathbb Z$ that reflects the factor structure as 
follows. 
We let $\{\mu_y:y\in Y\}$ and $\{\tilde{\mu}_y:y\in Y\}$ be the disintegrations of $\mathbb X$ and $\mathbb Z$ over $\mathbb Y$ respectively. The \emph{relatively independent joining} of $\mathbb X$ and $\mathbb Z$ over $\mathbb Y$ is the joining $\rho$:
\[\rho=\int (\mu_y\times \tilde{\mu}_y)d\nu(y).
\]
We will sometimes  write this  as $X\times_Y Z$. 

\medskip
We will be concerned about categories of measure preserving systems where the morphisms are joinings. For this we must describe the composition operation. Suppose we are given   joinings $\rho_{XY}$ between $X$ and $Y$ and 
$\rho_{YZ}$ between $Y$ and $Z$. Then $(Y,\nu)$ is a common factor of both  $(X\times Y, \rho_{XY})$ and $(Y\times Z,\rho_{YZ})$ and we can consider the relatively independent joining 
$\rho_{XY}\times_Y\rho_{YZ}$.
\bigskip

We define the \hypertarget{composition of joinings}{composition} of $\rho_{XY}$ and 
$\rho_{YZ}$ to be the projection of the relatively independent joining of $\rho_{XY}$ and $\rho_{YZ}$ to a measure on $X\times Z$. Formally, if $A\subseteq X\times Z$ and $\rho$ is the relatively independent joining, then: 
\[\rho_{XY}\circ\rho_{YZ}(A)=\rho(\{(x,y,z): x, z\in A\}).\]

\begin{ex}\label{for struct}
Suppose that $\pi_0:X\to Y$ and $\pi_1:Y\to Z$ are factor maps. If $\rho_{XY}$ is the joining associated with $\pi_0$ and $\rho_{YZ}$ is the joining associated with $\pi_1$, then 
$(\rho^*_{YZ}\circ\rho^*_{XY})^*$ is the joining associated with the factor map $\pi_1\circ\pi_0:X\to Z$.\footnote{In the following, in the context of factor maps $\pi:X\to Y$ we will be sloppy about whether this is associated with a joining of $X$ with $Y$ or a joining of $Y$ with $X$.}
\end{ex}

The following are standard facts (e.g. in \S 6.2 of \cite{glasbook}):
\begin{prop}\label{laziness}
\begin{enumerate}
\item The operation of composition of joinings  is associative: if $\rho_1, \rho_2$ and $\rho_3$ are joinings, then 
\[(\rho_1\circ\rho_2)\circ\rho_3=\rho_1\circ(\rho_2\circ\rho_3).\]
\item Suppose that $\pi^X:X\to X'$ and $\pi^Z:Z\to Z'$ are factor maps Let 
$\rho_1$ and $\rho_2$ be joinings of $X,Y$ and $Y,Z$ respectively. Let $\rho_1^\pi$ be the projection of $\rho_1$ to a joining of $X'$ and $Y$ via $\pi^X\times id$\ 
and $\rho_2^\pi$ be defined similarly. Finally let $(\rho_1\circ\rho_2)^\pi$ be the projection of the composition of $\rho_1$ and $\rho_2$ to a joining of $X$ with $ Z$. Then:
\[\rho_1^\pi\circ\rho_2^\pi=(\rho_1\circ \rho_2)^\pi.\]
\end{enumerate}
\end{prop}

\subsection{Symbolic Systems}
\label{symbolic shifts}
Let $\Sigma$ be a countable or finite alphabet endowed with the discrete topology. Then $\Sigma^\poZ$ can be given the product topology, which makes it into a separable, totally disconnected space that is compact if $\Sigma$ is finite.  

\medskip

\bfni{Notation:} If $u=\la \sigma_0, \dots \sigma_{n-1}\ra\in \Sigma^{<\infty}$ is a finite sequence of elements of $\Sigma$, then we denote the cylinder set based at $k$ in $\Sigma^\poZ$ by writing $\la u\ra_k$. If $k=0$ we abbreviate this and write $\la u\ra$. 
Explicitly: $\la u\ra_k=\{f\in \Sigma^\poZ: f\rest[k,k+n)=u\}$. The collection of cylinder sets form a base for the product topology on $\Sigma^\poZ$.
\medskip

{ \bfni{Notation:} For a word $w\in \Sigma^{<\nn}$ we will write  $|w|$ for the length of $w$.} We will write $1_{\la w\ra}$ for the characteristic function of the interval $\la w\ra_0$ in $\Sigma^\poZ$.

\medskip

\noindent The shift map:
\[sh:\Sigma^\poZ\to \Sigma^\poZ\]
defined by setting $sh(f)(n)=f(n+1)$ is a homeomorphism. If $\mu$ is a shift invariant Borel measure then the resulting measure preserving system $(\Sigma^\poZ, \mcb,\mu, sh)$ is called a \emph{symbolic system}. The closed support of $\mu$ is a shift invariant closed subset of $\Sigma^\poZ$ called a \emph{symbolic shift} or \emph{sub-shift}.

Symbolic shifts are often described intrinsically by giving a collection of words that constitute a clopen basis for the support of an invariant measure.  
Fix a language $\Sigma$,  and a sequence of collections of words $\la\mcw_n:n\in\nn\ra$ with the properties that:
\begin{enumerate}
\item for each $n$ all of the words in $\mcw_n$ have the same length $q_n$,
\item each $w\in\mcw_{n}$ occurs  at least once as a subword of each $w'\in \mcw_{n+1}$,
\item \label{not too much space} there is a summable sequence $\la \epsilon_n:n\in\nn\ra$ of positive numbers such that for each $n$, every word $w\in \mcw_{n+1}$ can be uniquely parsed into segments 
\begin{equation}u_0w_0u_1w_1\dots  w_lu_{l+1}\label{words and spacers}
\end{equation}
 such that each $w_i\in \mcw_n$, $u_i\in \Sigma^{<\nn}$ and for this parsing
\begin{equation*}
{\sum_i|u_i|\over q_{n+1}}<\epsilon_{n+1}.
\end{equation*}
\end{enumerate}

\noindent The segments $u_i$ in condition \ref{words and spacers} are called the \emph{spacer} or \emph{boundary} portions of $w$.
\begin{definition}A sequence $\la \mcw_n:n\in\nn\ra$ satisfying properties 1.)-3.) will be called a \emph{construction sequence}.
\end{definition}
 Associated with a construction sequence is a symbolic shift defined as follows. Let $\bk$  be the collection of 
 $x\in \Sigma^\poZ$ such that every finite contiguous subword of $x$ occurs inside some $w\in \mcw_n$. Then $\bk$ is a 
 closed shift invariant subset of $\sz$ that is compact if $\Sigma$ is finite.\footnote{  The symbolic shifts built from construction sequences coincide with transformations built by \emph{cut-and-stack} constructions.}

 Formally, we have constructed a symbolic shift.  To get a measure preserving system we find a shift invariant measure $\mu$ concentrating on $\bk$ and write $(\bk,\mu)$. In \cite{prequel} we define the notion of a \emph{uniform} construction sequence and show that the resulting $\bk$ are uniquely ergodic.

 \medskip

We want to be able to unambiguously parse elements of $\bk$. For this we will use construction sequences consisting of uniquely readable words.
\begin{definition}
Let $\Sigma$ be a language and $\mcw$ be a collection of finite words in $\Sigma$. Then $\mcw$ is \emph{uniquely readable} iff whenever $u, v, w\in \mcw$ and $uv=pws$ then either $p$ or $s$ is the empty word.
\end{definition}

 In our constructions we will  restrict our measures to a natural set:

\begin{definition}\label{def of S} Suppose that $\la \mcw_n:n\in\nn\ra$ is a construction sequence  for a symbolic system $\bk$ with each $\mcw_n$ uniquely readable. Let
$S$ be the collection $x\in \bk$ such that there are  sequences of natural numbers 
$\la a_m: m\in\nn\ra$, $\la b_m: m\in\nn\ra$ going to infinity such that for all $m$ there is an 
$n, x\rest [-a_m, b_m)\in \mcw_n$.
\end{definition}

\noindent Note that  $S$ is a dense shift invariant $\mathcal G_\delta$ set.
 The following lemma is routine:

\begin{lemma} 
Fix a construction sequence $\la\mcw_n:n\in\nn\ra$ for a symbolic system $\bk$ in a finite language. Then:
\begin{enumerate}
\item $\bk$ is the smallest shift invariant closed subset of $\Sigma^\poZ$ such that for all $n$,  and $w\in\mcw_n$, $\bk$ has non-empty intersection with the basic open interval $\la w\ra\subset \Sigma^\poZ$.

\item Suppose that there is a unique invariant measure $\nu$ on $S\subseteq \bk$, then $\nu$ is ergodic.

\end{enumerate}
\end{lemma}
\pf Item 1 is clear from the definitions.
If $X$ is a Polish space, $T:X\to X$ is a Borel automorphism and  $D$ is a $T$-invariant Borel set with a unique $T$-invariant measure on $D$, then that measure must be ergodic.
\qed

Let $\la \mcw_n:n\in\nn\ra$ be a uniquely readable construction sequence, and $s\in S$. By the unique readability, for each $n$   either $s(0)$ lies in a well-defined subword  of $s$ belonging to 
$\mcw_n$ or in 
a spacer of a subword of $s$ belonging to some  $\mcw_{n+k}$.

\begin{lemma}\label{principal blocks exist} Suppose that $\bk$ is built from $\la \mcw_n:n\in\nn\ra$ and $\nu$ is a shift invariant measure on $\bk$ concentrating on $S$.  Then for $\nu$-almost every $s$ there is an $N$ for all $n>N$, there are $a_n\le 0< b_n$ such that $s\rest[a_n, b_n)\in \mcw_n$.
\end{lemma}

\pf Let $B_{n}$  be the collection of $s\in S$ such that for some  $a_{n}\le 0<b_{n}$, $s\rest [a_{n}, b_n)\in \mcw_n$ but $s(0)$ is in a boundary portion of $s\rest[a_n,b_n)$. By the Ergodic Theorem and clause 3.) of the definition of a construction sequence $\sum\nu(B_n)<\infty$.  

It follows from the Borel-Cantelli Lemma that for almost all $s$ there is an $N$ such that for all 
$n\ge N$, $s\notin B_n$. Fix an $s\in S$ and such an $N$.  From the definition of $S$ there are arbitrarily large $n^*>N$ and $a_{n^*}\le 0<b_{n^*}$ such that $s\rest[a_{n^*}, b_{n^*})\in \mcw_{n^*}$. 
Using backwards induction from $n^*$ to $N$ and the definition of $B_n$,  this also holds for all $n\in [N, n^*)$. 
\qed

\subsection{Locations}

By Lemma \ref{principal blocks exist} for $\nu$-almost all $x$ and  for all large enough $n$ there is a unique $k$ with $0\le k<q_n$ such that 
$s\rest[-k, q_n-k)\in \mcw_n$.

\begin{definition}\label{def of rn} Let $s\in S$ and suppose that for some $0\le k<q_n, s\rest[-k,q_n-k)\in \mcw_n$.
We define $r_n(s)$ to be the unique $k$ with with this property.   We will call the interval $[-k, q_n-k)$ the \emph{principal $n$-block} of $s$, and $s\rest [-k, q_n-k)$ its \emph{principal $n$-subword}. The sequence of $r_n$'s will be called the \emph{location sequence of $s$}.
\end{definition}
We interpret $r_n(s)=k$ as saying that \emph{$s(0)$ is the $k^{th}$ symbol in the principal $n$-subword  of $s$ containing $0$.} We can view the principal $n$-subword of $s$ as being located on an interval $I$ inside the principal $n+1$-subword. Counting from the beginning of the principal $n+1$-subword, the $r_{n+1}(s)$ position is located at the $r_n(s)$ position in $I$. 

\medskip

\begin{remark} \label{interval coherence}Suppose that $s\in S$ has a principal $n$-block for all $n\ge N$. Let $N\le n<m$. It follows immediately from  the definitions that   $r_n(s)$  and $r_m(s)$ are well defined  and  the $r_m(s)^{th}$ position of the principal $m$-block of $s$  is in the $r_n(s)^{th}$ position inside the principal $n$-block of $s$. 

\end{remark}

The next lemma tells us that an element of $s$ is determined by   knowing  any tail of the  sequence  $\la r_n(s):n\ge N\ra$ together with  a tail of the principal subwords of $s$. 

\begin{lemma}\label{specifying elements} Suppose that $s, s'\in S$ and $\la r_n(s):n\ge N\ra=\la r_n(s'):n\ge N\ra$ and for all $n\ge N$, $s$ and $s'$ have the same principal $n$-subwords. Then $s=s'$.
\end{lemma}

\pf Since $s, s'\in S$ there are sequences $\la a_n, a_n', b_n, b_n':n\ge N\ra$  tending to infinity such that 
$s\rest[-a_n, b_n)\in \mcw_n$ and $s'\rest[a'_n, b'_n)\in \mcw_n$. Since $r_n(s)=r_n(s')$ we know that $a_n=a_n'$ 
and $b_n=b_n'$. Since $s$ and $s'$ have the same principal subwords, $s\rest[a_n, b_n)=s'\rest[a_n', b_n')$. The 
lemma follows.\qed
\begin{remark}We record some consequences of Lemma \ref{specifying elements}: \label{rebuilding}
\begin{enumerate}
\item Suppose that we are given a sequence $\la u_n:M\le n\ra$ with $u_n\in\mcw_n$.   If we specify which occurrence of $u_n$ in $u_{n+1}$ is the principal occurrence, and the distances of the principle occurrence to the beginning of $u_{n+1}$ go to infinity, then $\la u_n:M\le n\ra$ determines an $s\in S\subseteq \bk$ completely up to a shift $k$ with $|k|\le q_M$.

\item  
A sequence $\la r_n:N\le n\ra$ and sequence of words $w_n\in \mcw_n$ comes from an infinite word $s\in S$ if  both $r_n$ and $q_n-r_n$  go to infinity and 
that the $r_{n+1}$ position in $w_{n+1}$ is in the $r_n$ position in  a subword of $w_{n+1}$ identical to $w_n$.

\emph{Caveat}: just because $\la r_n:N\le n\ra$ is the location sequence of some $s\in S$ and $\la w_n:N\le n\ra$ is the sequence of principal subwords of some $s'\in S$, it does not follow that there is an $x\in S$ with location sequence $\la r_n:N\le n\ra$ and sequence of subwords $\la w_n:N\le n\ra$.

\item If
 $x, y\in S$ have the same principal $n$-subwords and $r_n(y)=r_n(x)+1$ for all large enough $n$, then $y=sh(x)$. 
\end{enumerate}
\end{remark}

\subsection{A note on inverses of symbolic shifts}\label{note on inverses}
We define operators we label $\rev{}$, and apply them in several contexts
\begin{definition} If $x$ is in $\bk$, we define the 
reverse  of $x$ by setting $\rev{x}(k)=x(-k)$.  For $A\subseteq \bk$, define:
\hypertarget{reverse}{
\[\rev{A}=\{\rev{x}:x\in A\}.\]}
If $w$ is a word, we define $\rev{w}$ to be the reverse of $w$. If we are viewing $w$ as sitting on an interval, we take $\rev{w}$ to sit on the same interval. Similarly, if $\mcw$ is a collection of words, $\rev{\mcw}$ is the collection of reverses of the words in $\mcw$.
\end{definition}

If $(\bk, sh)$ is an arbitrary symbolic shift then its inverse is $(\bk, sh^{-1})$.   It will be convenient to have all of our shifts go in the same direction, thus:
\begin{prop}\label{spinning}
The map $\phi$ sending $x$ to $\rev{x}$ is a canonical isomorphism between $(\bk, sh\inv)$ and $(\rev{\bk},sh)$. \end{prop}
 We will use the notation $\bl^{-1}$ for the system $(\bl,sh\inv)$ and $\rev{\bl}$ for the system $(\rev{\bl},sh)$.

We can say more. For a fixed symbolic shift $\bk$, the  canonical isomorphism $\phi:\bl^{-1}\to \rev{\bl}$ gives rise to a canonical correspondence 
\[\rho \leftrightarrow \rho'\] between joinings $\rho$ of $(\bk,sh)$ with 
$(\bl,sh\inv)$ and joinings $\rho'$ of $(\bk,sh)$ with $(\rev{\bl},sh)$.

We will also use the following remark.
\begin{remark}\label{nothing remarkable} 
{Assume that there is a unique non-atomic measure on a shift invariant set $S\subseteq\bk$. Then there is also a unique non-atomic shift invariant measure on 
 $\rev{S}$ and  for this measure, which we denote $\nu\inv$, we have 
 $\nu(\la w\ra)=\nu\inv(\la \rev{w}\ra)$.}
\end{remark}

 \subsection{Generic  points and sequences} 
 \label{sequences and points}
 
 Let $T$ be a measure preserving transformation from $(X,\mu)$ to $(X,\mu)$, where $X$ is a compact metric space. Let $C(X)$ be the space of all real valued complex functions. Then a point $x\in X$ is \emph{generic} for $T$ if and only if for all $f\in C(X)$, 
 \begin{equation*}\lim_{N\to \infty}\left({1\over N}\right)\sum_0^{N-1} f(T^n(x))=\int_X f(x)d\mu(x). \end{equation*}
 
 The Ergodic Theorem tells us that for a given $f$ and ergodic $T$ equation  above holds for a set of $\mu$-measure one. Intersecting over a countable dense set of $f$ gives a set of $\mu$-measure one of generic points. For symbolic systems $\bk\subseteq \Sigma^\poZ$ we can describe generic points $x$ as being those $x$ such that the $\mu$-measure of all basic open intervals $\la u\ra_0$ is equal to the density of $k$ such that $u$ occurs in $x$ at $k$.

 The symbolic systems we consider will be built from construction sequences and are characterized by the limiting properties of finite information. We now describe how this works in greater detail. A  more complete discussion of this can be found in \cite{Benjy}.

 Let $\mu$ be a shift invariant measure on a symbolic system $\bk$ defined by a uniquely readable  construction sequence $\la\mcw_n:n\in\nn\ra$ in a finite language $\Sigma$. Assume that $q_n$ is the length of the words in 
 $\mcw_n$. By $\mu_m$ we will denote the discrete measure on the finite 
 set $\Sigma^m$ given by $\mu_m(u)=\mu(\la u\ra)$. By $\hat{\mu}_n(w)$ we will denote the discrete probability measure on $\mcw_n$ defined by
\[\hat{\mu}_n(w)={\mu_{q_n}(\la w\ra)\over \sum_{w'\in\mcw_n} \mu_{q_n}(\la w'\ra)}.
\]

Thus $\hat{\mu}_n(w)$ is the relative measure of $\la w\ra$ among all $\la w'\ra, w'\in \mcw_n$. The denominator is a normalizing constant to account for spacers at stages $m>n$ and for shifts of size less than $q_n$.

Explicitly, if $A_n=\{s\in \bk: s(0)$ is the start of a word in $\mcw_n\}$, then the sets $\{sh^j(A_n)\}_{j=0}^{q_n-1}$ are disjoint and their union has  a measure that tends to one as $n$ grows to infinity. The set $A_n$ is partitioned into $|\mcw_n|$ many sets by the words $w\in \mcw_n$ and  $\hat{\mu}_n$ gives their relative size in $A_n$. Since the measure of an arbitrary finite cylinder set can be calculated along the individual columns represented by a fixed $w$, it is clear that the $\hat{\mu}_n(w)$ determine uniquely the measure $\mu$.
\medskip

Using the unique readability of words in $\mcw_k$ a word $w$ in $\Sigma^{q_{k+l}}$ determines a unique sequence of words $w_j$ in $\mcw_k$ such that , 
\[w=u_0w_0u_1w_1\dots w_Ju_{J+1}.\] 
When $w\in \mcw_{k+l}$,  each $u_j$ is in the region of spacers added  in $\mcw_{k+l'}$, $l'\le l$.
We will denote the \hypertarget{emptiest}{\emph{empirical distribution}} of $\mcw_k$-words in $w$ by EmpDist$_k(w)$. Formally:
\[\mbox{EmpDist}_k(w)(w')={|\{0\le j\le J: w_j=w'\}|\over J+1}, \ w'\in \mcw_k.\]
Then $EmpDist$ extends to a measure on $\mathcal P(\mcw_k)$ in the obvious way.

To finitize the idea of a generic point in $\bk$ we introduce the notion of a generic sequence of words. 

\begin{definition}\label{ED} A sequence $\la v_n\in\mcw_n:n\in\nn\ra$ is a \hypertarget{gen seq}{\emph{generic sequence of words}} if and only if
for all $k$ and $\epsilon>0$ there is an $N$ for all $m,n>N$,
\[\|EmpDist_k(v_m)-EmpDist_k(v_n)\|_{var}<\epsilon.\]
The sequence is generic for a measure $\mu$ if for all $k$:
\[\lim_{n\to \infty}\| \mbox{EmpDist}_k(v_n)-\hat{\mu}_k\|_{var}=0
\]
where $\|\ \|_{var}$ is the variation norm on probability distributions.
\end{definition}
It follows that if $\la v_n:n\in\nn\ra$ is a generic sequence of words then it is generic  for a unique measure $\mu$.
 Even though  Definition \ref{ED} involves only the measures $\hat{\mu}_k$ it is easy to see (using the Ergodic Theorem) that for any $u\in \Sigma^k$, if $\la v_n:n\in\nn\ra$ is generic then the density of the 
occurrences of $u$ in the $v_n$ will converge to $\mu(\la u\ra)$.
\medskip

 We can 
summarize the exact relationship between the empirical distributions and the 
$\mu_{q_k}$ by saying that the empirical distribution is the proportion of occurrences of 
$w'\in \mcw_k$ among the $k$-words that appear in $v_n$, whereas 
$\mu_{q_k}$ is approximately the density  of the locations of the start of $k$-words in $v_n$.  Letting $u\in \mcw_k$, $d$ be the density of the positions where an occurrence of $u$ begins in $v_n$, and $d_s$ be the density of locations of letters in some spacer $u_i$ we see that  these are related by:
\begin{eqnarray*}
d&=&\left({\mbox{EmpDist}(v_n)(u)\over q_k}\right)(1-d_s)
\end{eqnarray*}

\medskip
We record the following consequence of the Ergodic Theorem for future reference:

\begin{prop}\label{generic sequences exist for ergodic}
Let $\bk$ be an ergodic symbolic system with construction sequence $\la \mcw_n:n\in\nn\ra$ and measure $\mu$. Then for any generic $s$ the sequence of principal subwords of s,  $\la w_n:n\in\nn\ra$, is generic for $\mu$. In particular, generic sequences for $\mu$ exist.
\end{prop}

We will need a characterization of when  a generic sequence of words $\la w_n:n\in\nn\ra$ determines an  ergodic measure. 
\begin{definition}\label{def of ergodic sequence}
A sequence $\la v_n:n\in\nn\ra$ with $v_n\in\mcw_n$ is an \hypertarget{ergodic sequence}{\emph{ergodic sequence}} if for any $k$ and $\epsilon>0$ there are $n_0>k$, and $m_0$ such that for all $m\ge m_0$, if
\[v_m=u_0w_0u_1w_1u_2\dots u_Jw_Ju_{J+1}\]
is the parsing of $v_m$ into $\mcw_{n_0}$ words and spacers $u_i$ then there is a subset $I\subseteq \{0,1,2\dots J\}$ with $|I|/J>1-\epsilon$ and for all $j,j'\in I$
\begin{equation}
\|EmpDist_k(w_j)-EmpDist_k(w_{j'})\|_{var}<\epsilon.\label{erg seq equ}
\end{equation}
\end{definition}
Notice that in the definition of an ergodic sequence $\la v_n\ra$ we 
are not assuming that it is a generic sequence for a measure. This follow 
easily (see Lemma \ref{ergodic sequences give ergodic measures}), but we have not made it part of the 
definition to emphasize its finitary nature. In the next lemma we use the fact that the language $\Sigma$ is finite.

\begin{lemma}\label{generic seqs are ergodic}
Any generic sequence $\la v_n:n\in\nn\ra$ for an ergodic measure $\mu$ is an ergodic sequence.
\end{lemma}
\pf Suppose we are given $k$ and $\epsilon>0$. 
For all $\delta>0$ we can
 apply the Ergodic Theorem to find an $N$ much bigger than $q_k$ and a set 
$B$ with $\mu(B)>1-\delta$ such that for all $s\in B$ and all $w\in \mcw_k$:
\[\left| {1\over N}\sum_{0}^{N-1} 1_{\la w\ra}(T^is)-\mu_{q_k}(\la w\ra)\right|<\delta.\]
Fix a generic point $s$ for $\mu$. Let $I=\{i\ge 0: T^is\in B\}$, and define an infinite sequence of disjoint intervals of length $N$ that cover $I$ by inductively letting $i_0=min(I)$,  and $i_{j+1}=min(\{i\in I: i\ge i_j+N\})$. We take the intervals to be the sequence 
\[ [i_0, i_0+N-1], [i_1, i_1+N-1], [i_2, i_2+N-1], \dots \]
Notice that the complement of these intervals in $\poZ^+$ has density less than $\delta$ since their union clearly covers $I$.

Though this is an infinite sequence of intervals, the fact our language is finite implies that only finitely many distinct words of length $N$  occur as subwords of $s$ on these intervals. For each such word $w^*$, the density of those $i$ in the domain of  
$w^*$ such that an occurrence of a $w\in \mcw_k$ starts at $i$ is within $\delta$ of $\mu_{q_k}(\la w\ra)$.\footnote{By taking $N\gg q_k$, we can account for negligible ``end effects" so that
 $\left|{1\over N}\sum_{0}^{N-q_k-1} 1_{\la w\ra}(T^is)-\mu_{q_k}(\la w\ra)\right|<\delta$. We ignore end effects in the rest of the proof.}

Next take $n_0$ large enough that $N/q_{n_0}<\delta$, and  parse $s$ into words from $\mcw_{n_0}$ and the sections of $s$ corresponding to spacers in words in $W_j$ for some $j\ge n_0+1$. By taking $n_0$ large enough we can take the density of locations in $s$ occurring in spacers to be arbitrarily small.  Let $\delta'$ be this density.

The words from $\mcw_{n_0}$ have length much larger than $N$, and we can collect  all those words $w\in \mcw_{n_0}$ that are $(1-\sqrt{\delta})$-covered by the $N$-intervals we chose above into a set $A\subseteq \mcw_{n_0}$. 

The proportion of $s\rest \poZ^+$ not covered by words in  $A$ can be split into the spacer section and the portion inside words $w$ in $B=\mcw_{n_0}\setminus A$. For $w\in B$  the complement of the $N$-intervals has density at least 
$\sqrt{\delta}$.  It follows that the density of sections of $s$ covered by elements of $B$ is less than $\sqrt{\delta}$.

Thus the fraction of $s$ not covered by  words in $A$ is at most $\sqrt{\delta}+\delta'$. It is now clear that if 
$\delta, \delta'$ are chosen to be sufficiently small then 
\begin{equation}\label{whatever}
\sum_{w\in A}\hat{\mu}_{n_0}(w)>1-\epsilon   
\end{equation}
and  all $w\in A$ will have the property that 
\[\|\mbox{EmpDist}_k(w)-\hat{\mu}_k\|_{var}<\epsilon/2\]
which implies inequality \ref{erg seq equ} for pairs of words in $A$. Using inequality \ref{whatever} and the fact that $\la v_n\ra$ is generic 
for $\mu$ gives an $m_0$ so that for all $m\ge m_0$ when $v_m$ is parsed into $n_0$ 
words  a $(1-\epsilon)$-fraction will lie in $A$ 
and this concludes the proof. \qed

We will also need the converse to Lemma \ref{generic seqs are ergodic}, namely that the limiting measure defined by an ergodic sequence is, in fact, ergodic.

\begin{lemma}\label{ergodic sequences give ergodic measures}
An ergodic sequence is generic and the measure $\mu$ defined by an ergodic sequence $\la v_n:n\in\nn\ra$ is ergodic. 
\end{lemma}
\pf  Inequality  \ref{erg seq equ} implies that for each $k$ and $w\in \mcw_k$, the limit of the density of occurrences of $w$ in $v_n$ exists as $n$ goes to infinity. It follows (since $\mcw_k$ is finite) that $\la v_n:n\in\nn\ra$ is a generic sequence and hence it defines a unique measure $\mu$. 

The ergodicity of $\mu$ is equivalent to the fact that the ergodic averages of all $L^2$ functions converge almost everywhere to a constant. Functions of the form 
$1_{\la w\ra}$ where $w\in \bigcup_n\mcw_n$  and their shifts linearly span a dense set in $L^2$ from which it easily follows that if $\mu$ were not ergodic there would be some $k$, and $w\in \mcw_k$ with $({1/N})\sum_0^{N-1}1_{\la w\ra}(T^ix)$ converging $\mu$-a.e. to a non-constant function. This means that there is a $\gamma>0$ and disjoint sets $B_0, B_1$ of positive measure in $\bk$ such that for all large enough $N$ for all $x_0\in B_0,  x_1\in B_1$
\begin{equation}\label{gaining separation}
\left|{1\over N}\sum_0^{N-1}1_{\la w\ra}(T^ix_0)-{1\over N}\sum_0^{N-1}1_{\la w\ra}(T^ix_1)\right|\ge \gamma.
\end{equation}

Take $\epsilon$ small compared to $\gamma$ and $\mu(B_0), \mu(B_1)$.  Find $n_0, m_0$ as in the definition of \hyperlink{ergodic sequence}{ergodic sequence} for this $k$ and $\epsilon$. Choose $N$ large enough that inequality \ref{gaining separation} holds and so that $q_{n_0}/N$ is negligible. Finally take $m\ge n_0$ so that $N/q_m$ is negligible. 

The inequality \ref{gaining separation} depends only on the initial $(N+q_k)$-block of $x_0$ and $x_1$. Thus for large enough $m$ we can compute $\mu(B_0)$ and $\mu(B_1)$ by the empirical distributions of the $(N+q_k)$-blocks in $v_m$.

Since $N$ is large compared to $q_{n_0}$ the frequency of occurrence of $w$ in a block of length $N+q_k$ is determined by its frequencies in the words in $\mcw_{n_0}$ in the $n_0$-parsing of $v_m$. We now get a contradiction to inequality \ref{gaining separation}, since except for an $\epsilon$-fraction, these $w_{n_0}$-words have their $k$-words distributed very close to $\hat{\mu}_k(w)$.\qed

 If $S$ and $T$ are symbolic systems then a joining  $\rho$ of $S$ and $T$ will be a symbolic system, but  may not have well-defined construction sequence, even if $S$ and $T$ do.\footnote{We run into this problem when considering joinings of circular systems and their inverses that project to the $\natural$-map on the canonical factors; these notions are defined in future sections.} Accordingly we must generalize our definition of \emph{empirical distribution} to   take into account the relative locations of words in typical $(s,t)\in \bk\times \bl$.
We express this by shifting one of the basic open sets and considering words $(w, sh^s(v))$, which we view as starting at the locations $(0, s)$.

 Let $\la \mcw_n:n\in\nn\ra$ and $\la\mcv_n:n\in\nn\ra$ be uniquely readable construction 
 sequences for $\bk$ and $\bl$ in the languages $\Sigma, \Lambda$ respectively.  Assume for 
 simplicity that all words in $\mcw_n$ and $\mcv_n$ have the same length. 
 
 Let $n\le n'<n+l$. Then we can uniquely parse a word $w\in \mcw_{n+l}$ as 
 \[w=u_0w_0u_1w_1\dots w_Ju_{J+1}\]
where each $w_j\in \mcw_n$  and each $u_j$ is in the region of spacers for words in $\mcw_{n+l'}$, $l'<l$.  The similar statement holds for $v'_k\in \mcv_{n'}$, and $v\in \mcv_{n+l}$:
\[v=u'_0v'_0u'_1v'_1\dots v'_Ku'_{K+1}.\]
The definition must take into account the relative shifts of $w$ and $v$, the shifts of 
$(w_j, v_k)$ allow spacers to occur in different places and for the possibility  that $J\ne K$. \medskip

Let $n\le n'< n+ l$ be natural numbers, $s, s'\in \poZ$, and $(w', v')\in \mcw_n\times \mcv_{n'}$ and 
$(w,v)\in \mcw_{n+l}\times \mcv_{n+l}$.  Write $w$ and $v$ in terms of $n$ and $n'$-words as above. 
For $s, s'$,  define an \emph{occurrence} of 
$(w', sh^{s'}(v'))$ in $(w, sh^s(v))$ to be a $j\le J$ such that $w_j=w'$ and if $k$ is the location of $w_j$ in $w$, then $v'$ occurs at $k+s'$ in $sh^s(v)$. We note the bijection between occurrences of $(w', sh^{s'}(v'))$ in $(w, sh^s(v))$ and occurrences of $(v', sh^{-s'}(w'))$ in $(v, sh^{-s}(w))$. 

In defining  empirical distributions for joinings we generalize Definition \ref{ED}. The empirical distribution of a shifted pair is defined to be the proportion of times it occurs, relative to the proportion of times arbitrary pairs with the same shift occur.
\begin{definition}\label{empdist for joinings} Fix $w,s,v$
    Let $A$ be the collection  
    \begin{equation*}\{j:\mbox{for some }(w^*, v^*)\in \mcw_n\times \mcv_{n'}, (w^*, sh^{s'}(v^*))
    \mbox{ occurs at j}\\
    \mbox{ in }(w, sh^s(v))\}.
    \end{equation*}
    Assume that $A\ne \emptyset$. For $w'\in \mcw_n$ and $v'\in \mcv_{n'}$, we define: 
    \[EmpDist_{n,n', s'}(w,sh^{s}(v))(w',v')={|\{0\le j\le J:(w',sh^{s'}(v'))\mbox{ occurs at }
    j\}|\over
     |A|}.\]
      \end{definition}
     As before, $EmpDist_{n,n',s'}(w,sh^s(v))$ extends uniquely to a probability measure on $\mathcal P(\mcw_{n}\times \mcv_{n'})$.
    Definition \ref{empdist for joinings} facilitates a notion of a generic sequence for a joining.
     
     \begin{definition}\label{shifted genericity}
     A sequence of  $\la (w_n,v_n,s_n)\in \mcw_n\times \mcv_n\times \poZ:n\in\nn\ra$   is called \emph{generic} iff 
     \begin{enumerate}
        \item $\sum {|s_n|\over q_n}<\infty$ and
        \item  for all $n,n', s'$ and $\epsilon>0$ there is an $N$ for all $m,m'>N$,
     \[\|EmpDist_{n,n',s'}(w_m,sh^{s_{m}}(v_n))-EmpDist_{n,n',s'}(w_{m',}sh^{s_{m'}}(v_{m'}))\|_{var}<\epsilon.\]
     \end{enumerate}
     The definition of an \emph{ergodic sequence of pairs} is done analogously.
     \end{definition}

 It is easy to check that $\la (w_n, v_n, s_n):n\in\nn\ra$ is generic/ergodic   if and only if $\la (v_n,w_n,-s_n):n\in\nn\ra$ is generic/ergodic. For ergodic joinings the analogues of Proposition \ref{generic sequences exist for ergodic}, and Lemmas 
\ref{generic seqs are ergodic} and \ref{ergodic sequences give ergodic measures} hold and are proved in exactly the same way.

We have given these definitions in the case of a product of two symbolic shifts, but they generalize immediately to products of three or more shifts. For example, to consider three shifts with construction sequences $\la \mcu_n\ra_n, \la \mcv_n\ra_n, \la\mcw_n\ra_n $,  we would consider a sequence of the form: 
\[\la (u_n, v_n, w_n,  s_n, t_n):n\in\nn\ra,\]
where the words belong to the respective construction sequences and the $s_n$'s and $t_n$'s give the shifts relative to the first coordinate. 

We will be concerned with compositions of joinings, which involves products of three shifts. To prepare for this we need the notion of a conditional empirical distribution.

\begin{definition}\label{conditional love} Let $n, n'<n+l$.
Given a fixed $w^*\in \mcw_{n'}$ and a pair $(w,v)\in \mcw_{n+l}\times \mcv_{n+l}$ and $(s,s')$ we define the \hypertarget{cond emp dist}{\emph{conditional empirical distribution}} to be:
    \[\mbox{EmpDist}_{n,s'}((w,sh^s(v)|w^*)(v')=\]
    
    \[{|\{0\le j\le J:(w^*,sh^{s'}(v'))\mbox{ occurs at }j\}|\over 
    |\{j\le J: \mbox{for some }v^*\in \mcv_n, (w^*, sh^{s'}(v^*)) \mbox{ occurs at }j\}|}
    \]
for $v'\in \mcw_n.$
\end{definition}
Using the same ideas  we can define the  empirical distribution conditioned on a $v^*\in \mcv_k$ by looking at $(sh^{-s}(w),v)$ and counting occurrences of $(sh^{-s'}(w'),v^*)$ for the $w'\in \mcw_k$.

This definition  generalizes  to products of three or more systems. When working in three or more systems, there will be multiple $s$'s playing the role of $s'$ in Definition \ref{conditional love}. They will refer to the position of the sequences being counted, \emph{relative to the conditioning sequence}. So for example, if $\bk,\bl,\bm$ have construction sequences $\la \mcu_n\ra_n, \la \mcv_n\ra_n, \la\mcw_n\ra_n $ and $\la (u_n, v_n, w_n,  s_n, t_n):n\in\nn\ra$ is a generic sequence for a joining $\rho$ of $\bk, \bl$ and $\bm$, then 
\[EmpDists_{k,k',s,s'}(u_n,sh^{s_n}(v_n),sh^{t_n}(w_n)|v)\]
counts pairs $(sh^s(u),sh^{s'}(w))$, where  $(u,w)\in \mcu_k\times \mcw_{k'}$ have been shifted by  $s$ and $s'$  \emph{relative to }$v$.

\bigskip

 Let $\rho_1$ be a $T_1\times T_2$-invariant measure on $X\times Y$ and $\rho_2$ a $T_2\times T_3$-invariant measure on $Y\times Z$. Recall from Section \ref{joinings} that  the composition of 
 $\rho_1$ and $\rho_2$ is defined to be projection of the relative independent joining of $\rho_1$ and $\rho_2$ over the common factor $Y$ to a measure on $X\times Z$.
We now describe a method for detecting generic sequences for relatively independent joinings.

Suppose that systems $X$ and $Z$ have a common factor $Y$.

\begin{equation*}
\begin{diagram}
\node{(X,\mcb,\mu, T)}\arrow[1]{se}\node{}\node{(Z,\mcd, \tilde{\mu},\tilde{T} )}\arrow[1]{sw}\\
\node{}\node{(Y,\mcc,\nu,S)}
\end{diagram}
\end{equation*}

Let $\rho=X\times_Y Z$ be the relatively independent joining of $X$ and $Y$. Let $\mu_y, \tilde{\mu}_y, \rho_y$ be the distintegrations of $\mu,\tilde{\mu}$ and $\rho$ respectively. Then the relatively independent joining $\rho$ is characterized by the fact that  for $\nu$-a.e $y$, 
\begin{equation}\label{the world is disintegrating}
\rho_y=\mu_y\times \tilde{\mu}_y.
\end{equation}
 Let $\la \mca_n, \tilde{\mca}_n, \mca'_n:n\in\nn\ra$ be sequences of refining partitions that generate $\mcb, \mcd$ and $\mcc$ respectively. Since the sequence of partitions $\mca_n\times \tilde{\mca}'_n$ generates $\mcb\otimes\mcd$, equation \ref{the world is disintegrating} is equivalent to the property that for all $A_k\in \mca_k, \tilde{A}_k\in \tilde{\mca}_k$ and $\nu$-a.e. $y$,
 \begin{equation}
 \mu_y(A_k)\times\tilde{\mu}_y(\tilde{A}_k)=\rho_y(A_k\times\tilde{A}_k)\label{produce}
 \end{equation}
To finitize this we approximate $\mu_y(A_k)$ by $\mu(A_k|A'_m(y))$ for large $m$, where $A'_m(y)$ is the atom of $\mca'_m$ to which $y$ belongs.  We let
$\mu_y(\mca_k)$ be shorthand for the distribution $\la \mu_y(A_k):A_k\in \mca_k\ra$, and 
$\mu(\mca_k|\mca'_m)(y)$ stands for the conditional distribution 
$\mu(A_k|A'_m(y)), A_k\in \mca_k$. (We use similar notation in Lemma \ref{allow us to disintegrate} for the conditional distribution given by $\rho, \mu$ and $\tilde{\mu}$ on various partitions.)

By Martingale convergence,\footnote{See (e.g.) \cite{glasbook}, Theorem 14.26, page 261.} for 
$\epsilon>0$ and fixed $k$ if $m$ sufficiently large, then for $(1-\epsilon)$ proportion of the $y'$ in the same atom as $y$:
\[\|\mu_{y'}(\mca_k)-\mu(\mca_k|\mca'_m)(y)\|_{var}<\epsilon\]
but for a collection of $A'_m$ of whose union has $\nu$-measure less than $\epsilon$.

One can deal similarly with $\tilde{\mu}_n$ and $\rho_y$. We have shown:

\begin{lemma}\label{allow us to disintegrate}
    In the notation above, $\rho$ is the relatively independent joining of $\mu$ and $\tilde{\mu}$ if and only if for all $k, \epsilon>0$, for all large enough $m$, there is a 
    collection of atoms $A_m\in \mca'_m$ of total measure at least $1-\epsilon$ for which:
    \begin{equation}\label{relind}
   \| \rho(\mca_k\times \tilde{\mca_k}|A_m)-\mu(\mca_k|A_m)\times {\tilde{\mu}}(\tilde{\mca}_k|A_m)\|_{var}<\epsilon.
    \end{equation}
\end{lemma}

We now express Lemma \ref{allow us to disintegrate} in terms of sequences of finite words. Suppose that $\la \mcu_n\ra, \la \mcv_n \ra$, and $\la \mcw_n\ra$ are the uniquely readable construction sequences for $X$, $Y$ and $Z$.

\begin{prop}\label{rel ind join}
Let $\la (u_n, v_n, w_n, s_n, t_n)\in \mcu_n\times \mcv_n\times \mcw_n\times\poZ^2:n\in\nn\ra$ be a sequence of words. Suppose that:
\begin{enumerate}
    \item $\la (u_n, v_n, s_n)\ra_n$ is generic for $\rho_1$, 
    \item $\la (v_n, w_n, t_n)\ra_n$ is generic for $\rho_2$. 
    \item \label{hyp 3} for all $\epsilon>0, k$ and $s^*$ 
    for all sufficiently large  $k'$ there is an  $N$ and a set $G_{k'}\subset\mcv_{k'}$ and for each $v\in G_{k'}$ a set of indices $I_v\subseteq [0,q_{k'})$ that satisfies $|I_v|>(1-\epsilon)q_{k'}$ such that for all $n>N$:
    \begin{enumerate}
        \item \label{i} $\sum_{v\in G_{k'}}EmpDist(v_n)(v)>1-\epsilon$ 
        \smallskip 
         
         and 
         \smallskip
         
         \item \label{ii} for all $v\in G_{k'}$ and  $s\in I_{v}$, 
         \begin{eqnarray}
          \|EmpDist_{k,k,s,s+s^*}(u_n, sh^{s_n}(v_n),sh^{t_n}(w_n)|v)\ \ \ \ -\ \ \ \ \ \ \ \ \ \ \ \ \ \ \ \ \ \  &\notag\\
         EmpDist_{k,s}(u_n,sh^{s_n}(v_n)|v)*EmpDist_{k,s+s^*}(v_n,sh^{t_n-s_n}(w_n)|v)&\!\!\!\!\!\|_{var}\notag
         \end{eqnarray}
         is less than $\epsilon$.
    \end{enumerate}

\end{enumerate}
If $\rho$ is the relatively independent joining of $\rho_1, \rho_2$, then $\la (u_n, v_n, w_n, s_n, t_n):n\in\nn\ra$ is a generic sequence for $\rho$.

\end{prop}

\pf Observe that the hypothesis \ref{ii} implies a similar equation for any $k_1<k$ while the other parameters are fixed. Now use hypothesis \ref{i} with a summable sequence of $\epsilon$'s and we can conclude by the Borel-Cantelli lemma that for $\nu$-almost 
every $y\in Y$ for $k'$ sufficiently large,
if $v_{k'}(y)$ is the principal $k'$-block of $y$ with location $r_{k'}$, then the inequality in \ref{ii} will hold for $s=r_{k'}$ and $v=v_{k'}(y)$. 

Now by hypotheses 1 and 2, the single empirical distributions are converging to $(\rho_1)_y$ and $(\rho_2)_y$ respectively (where 
$(\rho_i)_y$ is the disintegration of $\rho_i$ over $y$).

It then follows by integration that the sequence of $(u_n, v_n, w_n, s_n, t_n)$'s is generic for a measure $\rho$ on $X\times Y \times Z$, which is the relatively independent joining.\qed

\begin{remark}
\label{grasping}
It follows immediately from hypothesis \ref{hyp 3} of Proposition \ref{rel ind join} that if we are given a finite set $F$ of natural numbers then for all sufficiently large $k'$ we can find an $N$, $G_{k'}$ and $I_v$ as in hypothesis \ref{hyp 3} so that (a) and (b) hold simultaneously for all $s^*\in F$.
\end{remark}

An immediate corollary of this is:
\begin{corollary}\label{en fin}
Suppose that $\la (u_n, v_n, w_n, s_n, t_n):n\in\nn\ra$ satisfies the hypotheses of Proposition \ref{rel ind join}. Then
$\la (u_n, sh^{t_n}(w_n)):n\in\nn\ra$ is generic for $\rho_1\circ\rho_2$.
\end{corollary}

There is a converse to Proposition 29, namely that a generic sequence
for the relatively independent joining of two odometer based system
satisfies the conditions 1-3 of the Proposition. The first two
are immediate while the third simply expresses the fact that the generic sequence sequence is actually representing the relatively independent joining. For later use we record this as: 
\begin{lemma}\label{ill existe}
Given joinings $\rho_1$ of $X\times Y$ and $\rho_2$ of $Y\times Z$ if $\la (u_n, v_n, w_n, s_n, t_n):n\in\nn\ra$ is generic for the relatively independent joining $\rho$ then it satisfies the hypotheses of Proposition \ref{rel ind join}.
\end{lemma}

\hypertarget{HvN}{\subsection{Unitary Operators}}\label{intro unitary}
We will use spectral tools introduced by Koopman and studied by Halmos and von Neumann. We reprise the basic facts we will use. Readers unfamiliar with this material can find it in \cite{walters} or \cite{glasbook}.
Let $\xbmt$ and $\ycns$ be measure preserving systems.

\medskip
 
If $T:X\to Y$ is a measure preserving transformation then $T$ induces a unitary isometry $U_T:L^2(Y)\to L^2(X)$ by 
setting 
\[U_T(f)=f\circ T.\] 
If $T$ is an isomorphism then $U_T$ is invertible. Moreover if $U:L^2(Y)\to L^2(X)$ is multiplicative on bounded functions then there is a measure preserving transformation $T:X\to Y$ such that $U=U_T$.

If $\pi:X\to Y$ is a factor map, then the map $f\mapsto f\circ\pi$ gives an injection of $L^2(Y)$ into $L^2(X)$, whose range is a closed $U_T$ invariant subspace. Conversely if $M\subseteq L^2(X)$ is a closed $U_T$ invariant  subspace containing 1 that is closed under taking complex conjugates, truncation and multiplication by elements of $M\cap L^\infty(X)$, then there is a factor $Y\subseteq X$ such that $M=L^2(Y)$.

For the rest of this discussion assume that $T$ is ergodic. Then the eigenvalues of $U_T$ all have multiplicity one and form a subgroup $G_T\subseteq \bt$. The group $G_T$ is an isomorphism invariant.

The collection of eigenfunctions generate a closed subspace of $L^2(X)$ corresponding to a factor $K$ of $X$. 
This factor is called the \emph{Kronecker factor}. If $H$ is any subgroup of $G_T$ then there is a further factor
 $K_H$ of $K$ that is canonically determined by the eigenfunctions coming from eigenvalues in $H$. 
 
 Assume  that $\phi$ is an isomorphism from $(X,T)$ to $(Y,S)$. Then $G_T=G_S$ and if $K^X_H, K^Y_H$ are the factors of $X$ and $Y$ determined by $H\subseteq G_T$ then $U_\phi$ determines an unique isomorphism between $K^X_H$ and $K^Y_H$. 
 
 It follows from this that if $\alpha\in \bt$ is an eigenvalue of $U_T$ then there are factors of $X$ and $Y$ isomorphic to rotation $\mcr_\alpha$ of $\bt$ by $\alpha$. Moreover there is a unique isomorphism $U_\phi^\pi:(\bt, \mcb, \lambda, \mcr_\alpha)\to (\bt, \mcb, \lambda, \mcr_\alpha)$ that intertwines $U_\phi$ and the projection maps of $X$ and $Y$ to $(\bt, \mcb, \lambda, \mcr_\alpha)$.
 
 The analogous statement holds for odometers. If $G_T$ consists of finite order eigenvalues and $\mco$ is the corresponding odometer transformation, then there is a unique isomorphism $U_\phi^\pi:\mco\to\mco$ that intertwines $U_\phi$ and the projection maps of $X$ and $Y$ to $\mco$.

  \subsection{Stationary Codes and $\dbar$-Distance}
In this section we briefly describe a standard idea, that of a \emph{stationary code} that we will  use to understand the existence of factor maps and isomorphisms. We review some standard facts here. A reader unfamiliar with this material who wants to see proofs should see \cite{Shields}. 
 \begin{definition}
Suppose that $\Sigma$ is a countable language. A \emph{code} of length $2N+1$ is a function $\Lambda:\Sigma^{[-N, N]}\to \Sigma$, where $[-N, N]$ is the interval of integers starting at $-N$ and ending at $N$.

Given a code $\Lambda$ and an $s\in \Sigma^\poZ$ we define the \emph{stationary code} determined by $\Lambda$ to be $\bar{\Lambda}(s)$ where:
\[\bar{\Lambda}(s)(k)=\Lambda(s\rest[k-N, k+N]).\]
\end{definition} 
Let $(\Sigma^\poZ, \mcb, \nu, sh)$ be a symbolic system.
Suppose we have two codes $\Lambda_0$ and $\Lambda_1$ that are not necessarily of the same length.
Define $D=\{s\in \Sigma^\poZ:\overline{\Lambda}_0(s)(0)\ne \bar{\Lambda}_1(s)(0)\}$ and $d(\Lambda_0, \Lambda_1)=\nu(D)$. Then $d$ is a semi-metric on the collection of codes.  The following is a consequence of the Borel-Cantelli lemma.

\begin{lemma}\label{cauchy codes}Let 
Suppose that $\la \Lambda_i:i\in\nn\ra$ is a sequence of codes such that $\sum_{i}d(\Lambda_i, \Lambda_{i+1})<\infty$. Then there is a shift invariant Borel map $S:\Sigma^\poZ\to \Sigma^\poZ$
such that for $\nu$-almost all $s$, $\lim_{i\to \infty}\overline\Lambda_i(s)=S(s)$

\end{lemma}

A shift invariant Borel map $S:\Sigma^\poZ\to \Sigma^\poZ$, determines a factor  $(\Sigma^\poZ, \mcb, \mu, sh)$ of $(\Sigma^\poZ, \mcb, \nu, sh)$ by setting $\mu=S^*\nu$ (i.e. $\mu(A)=\nu\circ S^{-1}(A)$). Hence a convergent sequence of stationary codes determines a factor of $(\Sigma^\poZ, \mcb, \nu, sh)$.

Let $\Lambda_0$ and $\Lambda_1$ be codes. Define $\dbar(\bar\Lambda_0(s), \bar\Lambda_1(s))$  to be
\[
	\overline{\lim}_{n\to \infty}{|\{k\in [-N, N]:	\bar\Lambda_0(s)(k)\ne 	\bar\Lambda_1(s)(k)\}|\over 2N+1}			
\]
More generally we can define the $\dbar$ metric on $\Sigma^{[a,b]}$ by setting 
\[\dbar_{[a,b]}(x,y)={|\{k\in [a, b):	x(k)\ne 	y(k)\}|\over b-a}.\]  For $x, y\in \Sigma^\poZ$, we set
\[\dbar(x,y)=\overline{\lim}_{N\to \infty}\dbar_{[-N,N]}(x\rest[-N,N], y\rest[-N,N]),\] 
provided this limit exists.

To compute distances between codes we will use the following  application of the Ergodic Theorem.
\begin{lemma}\label{computing code distances} Suppose that $(\Sigma^\poZ, sh,\nu)$ is ergodic and that $\Lambda_0$ and $\Lambda_1$  be codes. Then for almost all $s\in S$:
\[
d(\Lambda_0, \Lambda_1)=\dbar(\bar\Lambda_0(s), \bar\Lambda_1(s))\]
\end{lemma}

We finish with a useful remark:

\begin{remark}\label{cheating on dbar}
If $w_1$ and $w_2$ are words in a language $\Sigma$ defined on an interval $I$ and $J\subset I$ with ${|J|\over |I|}\ge \delta$, then $\dbar_I(w_1, w_2)\ge \delta\dbar_J(w_1, w_2)$.
 
 \end{remark}

\section{Odometer based and Circular Symbolic Systems}\label{odometer based and circular symbolic systems}
Two types of symbolic shifts play central roles for the proofs of our main theorem. We dub them \emph{odometer based} and \emph{circular} systems. In this section we give some general facts about symbolic systems with uniquely readable construction sequences, define odometer and circular systems, and show that every circular system has a canonical rotation factor.
\subsection{Odometer Based Systems}
We recall the definition of an odometer transformation.  Let $\la k_n:n\in\nn\ra$ be a sequence of natural numbers greater than or equal to 2. Let 
\[O=\prod_{n=0}^\infty \poZ/k_n\poZ\]
be the $\la k_n\ra$-adic integers. 
Then $O$ naturally has a compact abelian  group structure and hence carries a Haar measure $\mu$. 
{We make $O$ into a measure preserving system 
$\mco$ by defining $T:O\to O$ to be addition by 1 in the $\la k_n\ra$-adic integers. Concretely, this is the map that ``adds one to $\poZ/k_0\poZ$ and carries right".} Then $T$ is  an invertible transformation that preserves the Haar measure $\mu$ on $\mco$. Let $K_n=k_0*k_1*k_2\dots k_{n-1}$.

The following results are standard:
\begin{lemma} \label{odometer basics}Let $\mco$ be an odometer system. Then:
\begin{enumerate}
\item $\mco$ is ergodic.
\item The map $x\mapsto -x$ is an isomorphism between $(O, \mcb, \mu, T)$ and $(O, \mcb, \mu, T^{-1})$.
\item {Odometer maps are transformations with discrete spectrum and the eigenvalues of the associated linear operator are the $K_n^{th}$ roots of unity ($n>0$).}
\end{enumerate}
\end{lemma}

Any natural number $a$ can be uniquely written as:
\[a=a_0+a_1k_0+a_2(k_0k_1)+ \dots +a_j(k_0k_1k_2\dots k_{j-1})\] 
for some sequence of natural numbers $a_0, a_1, \dots a_j$ with $0\le a_j<k_j$.

\begin{lemma}\label{specifying elements of the odometer} Suppose that 
$\la r_n:n\in\nn\ra$ is a sequence of natural numbers with $0\le r_n<k_0k_1\dots k_{n-1}$ and $r_n\equiv r_{n+1} \mod(K_n)$.  Then there is a unique element $x\in O$ such that $r_n=x(0)+x(1)k_0+\dots x(n)(k_0k_1\dots k_{n-1})$ for each $n$.

\end{lemma}

We now define the collection of symbolic systems that have odometer maps as their timing mechanism. This timing mechanism can be used to parse typical elements of the symbolic system.

\begin{definition}Let $\la \mcw_n:n\in\nn\ra$ be a uniquely readable construction sequence with the properties that {$\mcw_0=\Sigma$} and for all $n, \mcw_{n+1}\subseteq (\mcw_n)^{k_n}$ for some $k_n$. The  associated symbolic system will be called an \emph{odometer based system}.
\end{definition}
Thus odometer
based systems are those built from  construction sequences $\la \mcw_n:n\in\nn\ra$ such that the words in 
$\mcw_{n+1}$ are concatenations of words in $\mcw_n$ of a fixed length $k_n$. The words in $\mcw_{n}$ all have length $K_n$ and the words $u_i$  in  equation \ref{words and spacers} are all the empty words.

Equivalently, an odometer based transformation is one that can be built by a cut-and-stack construction using no spacers. An easy consequence of the definition is that for odometer based systems $\bk$, for all $s\in\bk$ and for all $n\in\nn$, $r_n(s)$ exists.

\begin{prop} \label{kiss your S goodbye}
Let $\bk$ be an odometer based system and suppose that $\nu$ is a shift invariant measure. Then 
$\nu$ concentrates on $S$.
\end{prop}
\pf Let $B=\bk\setminus S$. Then $B$ is shift invariant.  Suppose that $\nu$ gives $B$ positive measure. For $s\in B$ let $a_n(s)\le 0\le b_n(s)$ be the left and right endpoints of the principal $n$-block of $s$. Then for all $s\in B$ there is an $N\in\nn$ such that:
\begin{enumerate}
\item for all $n, -N\le a_n$ or
\item for all $n, b_n\le N$.
\end{enumerate}
We assume that $\nu$ gives the collection $B^*$ of $s$ such that there is an $N\in \nn$ for all $n, -N\le a_n$ positive measure, the other case is  similar.

Define $f:B^*\to \nn$ by setting $f(s)=$ least $N$ satisfying item 1. Then $f$ is a Borel function. Let $B_n=f^{-1}(n)$.  Then the $B_n$'s are disjoint,  $B^*=\bigcup_{n\ge 0}B_n$ and $sh^{-1}(B_n)=B_{n+1}$.  Hence for all $n, m, \nu(B_n)=\nu(B_m)$, a contradiction.
\qed

 The next lemma justifies our terminology.
 \begin{lemma}\label{odometer factor}
 Let $\bk$ be an odometer based system with each $\mcw_{n+1}\subseteq (\mcw_n)^{k_n}$. Then there is a canonical factor map
 \begin{equation*}
 \pi:S\to \mco,
 \end{equation*}
 where $\mco$ is the odometer system determined by $\la k_n:n\in\nn\ra$.
 \end{lemma}
\pf For each $s\in S$, we know that for all $n, r_n(s)$ is defined and both $r_n$ and $k_n-r_n$ go to infinity. By 
Lemma \ref{specifying elements of the odometer}, the sequence $\la r_n(s):n\in\nn\ra$ defines a unique element $\pi(s)$ in 
$\mco$. It is easily checked that $\pi$ intertwines $sh$ and $T$.\qed

In the forthcoming paper \cite{part4} we show a strong converse to this result: if $T$ has finite entropy and an odometer factor then $T$ can be presented by an odometer based system.

Heuristically, the odometer transformation $\mco$ parses the sequences $s$ in $S\subseteq \bk$ by indicating where 
the words constituting $s$ begin and end. Shifting $s$ by one unit shifts this parsing by one. We can understand 
elements of $s$ as being an element of the odometer with words in $\mcw_n$ filled in inductively.

We will use the following remark about the canonical factor of  the inverse of an odometer based system.
\begin{remark}
If $\pi:\bl\to\mco$ is the canonical factor map, then the function $\pi:\bl\to O$ is  also factor map from $(\bl,sh\inv)$ to $\mco\inv$ (i.e. $O$ with the operation ``$-1$").  If $\la\mcw_n:n\in\nn\ra$ is the construction sequence for $\bl$, then $\la \rev{\mcw_n}:n\in\nn\ra$ is a construction sequence for $\rev{\bl}$. If $\phi:\bl\inv \to \rev{\bl}$ is the canonical isomorphism given by Proposition \ref{spinning}, then Lemma  \ref{odometer basics} tells us that the projection of $\phi$ to a map $\phi^\pi:\mco\to \mco$ is given by  $x\mapsto -x$.
\end{remark}
From this remark we immediately see:

\begin{lemma}\label{joining correspondence} Let $\rho\leftrightarrow \rho'$ be the canonical correspondence between joinings of $(\bk,sh)$ and $(\bl, sh^{-1})$ and joinings of $(\bk, sh)$ and $(\rev{\bl},sh)$ given after Proposition \ref{spinning}.
Then the joining $\rho$ concentrates on the set of pairs $(s,t)$ such that $\pi^\bk(t)=-\pi^\bl(s)$ if and only if $\rho'$ concentrates on the collection of $(s,t)$ such that $\pi^\bk(s)=\pi^{\bl^{-1}}(t)$.
\end{lemma}

\subsection{Circular systems} \label{circular systems 1} We now define and  discuss circular systems. The paper \cite{prequel}  showed that the circular systems give symbolic  characterizations of the smooth diffeomorphisms defined by the Anosov-Katok method of conjugacies.
The construction sequences of circular systems have quite specific combinatorial properties that will be important to our understanding of the Anosov-Katok systems and their centralizers in the third paper in this series.

 We call these systems \emph{circular}  because they are closely tied to the behavior of rotations by a convergent sequence of rationals $\alpha_n=p_n/q_n$.  The rational rotation by $p/q$ permutes the $1/q$ intervals of the circle cyclically along a sequence determined by some numbers $j_i=_{def}p^{-1}i$ (mod $q$): the interval $[i/q, (i+1)/q)$ is the $j_i^{th}$ interval in the sequence.\footnote{We assume that $p$ and $q$ are relatively prime and the exponent $-1$ is the multiplicative inverse of $p$ mod $q$.} The operation $\mcc$ which we are about to describe models the relationship between rotations by $p/q$ and $p'/q'$ when $q'$ is very close to $q$.

Let $k, l, p, q$ be positive natural numbers with $p<q$ relatively prime. Set 
\begin{equation}j_i\equiv_q(p)^{-1}i \label{j sub i}
\end{equation}
with $j_i<q$. It is easy to verify that:
\begin{equation}\label{reverse numerology} 
q-j_i=j_{q-i} 
\end{equation}

Let $\Sigma$ be a non-empty set. We define an operation $\mcc$, which depends on $p,q$, an integer $l>1$, and on sequences  $w_0, \dots w_{k-1}$ of words in a language $\Sigma\cup \{b, e\}$ 
by setting:\footnote{We use $\prod$ for repeated concatenation of words.}

\begin{equation}\mcc(w_0,w_1,w_2,\dots w_{k-1})=\prod_{i=0}^{q-1}\prod_{j=0}^{k-1}(b^{q-j_i}w_j^{l-1}e^{j_i}). \label{definition of C}
\end{equation} 
To start our construction we frequently take $p_0=0$ and $q_0=1$. In this case we adopt the convention that $j_0=0$.  Hence 
\begin{eqnarray*}
\mcc(w_0,w_1, \dots w_{k-1})&=&\prod_{j<k}b^qw_j^{l-1}\\
&=&\prod_{j<k}bw^{l-1}.
\end{eqnarray*}

{\begin{remark}\label{word length} We remark:
\begin{itemize}
\item Suppose that each $w_i$ has length $q$, then the length of  $\mcc(w_0, w_1, \dots w_{k-1})$ is $klq^2$. 
\item Every occurrence of an $e$ in $\mcc(w_0, \dots w_{k-1})$ has an occurrence of a $b$ to the left of it. If $p\ne 0$ then every occurrence of a $b$ has an $e$ to the right of it.
\item {Suppose that $n<m$ and $b$ occurs at position $n$ in $\mcc(w_0, w_1,\dots w_{k-1})$  and $e$ occurs at $m$ and neither occurrence is in a $w_i$. Then there must be some $w_i$ occurring between $n$ and $m$.}
\end{itemize}
\end{remark}}

The $\mcc$ operator automatically creates uniquely readable words, as the next lemma shows, however we will need a  stronger unique readability assumption for our definition of circular systems.
\begin{lemma}\label{unique readability}
Suppose that $\Sigma$ is a  language, $b, e\notin \Sigma$, $0<p<q$ and that  $u_0, \dots $ $u_{k-1}$, $v_0, \dots v_{k-1}$ and $w_0 \dots w_{k-1},$ are words in the language $\Sigma\cup \{b, e\}$ of some fixed length $q<l/2$. Let 
\begin{eqnarray*}u&=&\mcc(u_0, u_1, \dots u_{k-1}) \\
v&=&\mcc(v_0, v_1, \dots v_{k-1})\\
w&=&\mcc(w_0, w_1, \dots w_{k-1}).
\end{eqnarray*} Suppose that $uv$ is written as $pws$ where $p$ and $s$ are words in $\Sigma\cup \{b, e\}$. Then either $p$ is the empty word and $u=w, v=s$ or $s$ is the empty word and $u=p, v=w$. 
\end{lemma}
\pf The map $i\mapsto j_i$ is one-to-one. Hence each location in the word of length $klq^2$ is uniquely determined by the lengths of nearby sequences of $b$'s and $e$'s.\qed

In fact something stronger is true: if $\sigma\in \Sigma$ occurs at place $m$ in $w$ then $m$ is uniquely determined by the knowing the $w_0, w_1, \dots w_{k-1}$ and the $kq^l/2 +1$ letters on either side of $\sigma$. 

We now describe how to use the $\mcc$ operation to build a collection of symbolic shifts. Our systems will be defined using a sequence of natural number parameters $ k_n$ and $l_n$ that is fundamental to the version of the Anosov-Katok construction  presented in \cite{katoksbook}.

 Fix an arbitrary sequence of positive natural numbers $\la k_n:n\in\nn\ra$. Let {$\la l_n:n\in\nn\ra$} be an increasing sequence of natural numbers such that $\sum_n 1/l_n<\infty$. From the $k_n$ and $l_n$ we define  sequences of numbers: $\la p_n, q_n, \alpha_n:n\in\nn\ra$. 
 We begin by letting $p_0=0$ and $q_0=1$  and inductively set 
\begin{eqnarray}\label{qn}
\qnpo=\kn l_n\qn^2
\end{eqnarray}
{(thus $q_1=k_0l_0$)} and take 
\begin{equation}\label{pn} p_{n+1}=p_nq_nk_nl_n+1.
\end{equation}
Then clearly $p_{n+1}$ is relatively prime to $q_{n+1}$.\footnote{{$p_n$ and $q_n$ being relatively prime for $n\ge 1$, allows us to define the integer $j_i$  in equation \ref{j sub i}.  For $q_0=1$, $\poZ/q_0\poZ$ has one element, $[0]$, so we set $p_0\inv=p_0=0$.}}

\begin{definition}
A sequence of integers $\la k_n, l_n:n\in\nn\ra\ra$ such that $k_n\ge 2$, $\sum 1/l_n<\infty$  will be called a \hypertarget{circ coef}{\emph{circular coefficient sequence}}.
\end{definition}

\bigskip
 Let $\Sigma$ be a non-empty finite or countable alphabet. We will construct the systems we study by building collections of words $\mcw_n$ in the alphabet $\Sigma\cup \{b, e\}$ by induction as follows:
 \begin{itemize}
 \item Fix a  \hyperlink{circ coef}{circular coefficient sequence} $\la k_n, l_n:n\in\nn\ra\ra$.
 \item Set $\mcw_0=\Sigma$.  
 \item Having built $\mcw_n$ we choose a set $P_{n+1}\subseteq (\mcw_{n})^{k_n}$ and  form $\mcw_{n+1}$ by taking all words of the form $\mcc(w_0,w_1\dots w_{k_n-1})$  with $(w_0, \dots w_{k_n-1})\in P_{n+1}$.\footnote{Passing from $\mcw_n$ to $\mcw_{n+1}$ we use $\mcc$ with parameters $k=k_n, l=l_n, p=p_n$ and $q=q_n$ and take $j_i=(p_n)^{-1}i$ modulo $q_n$.  By Remark \ref{word length}, the length of each of the words in $\mcw_{n+1}$ is $q_{n+1}$.} 
 \end{itemize}
We will call the elements of $P_{n+1}$ \hypertarget{pwords}{\emph{prewords}.}

 \medskip 
 \hypertarget{strong unique readability}{\bfni{Strong Unique Readability Assumption:}}
  Let $n\in \nn$, and view $\mcw_n$ as a collection $\Lambda_n$ of letters. Then each element of $P_{n+1}$ can be viewed as a word  with letters in   $\Lambda_n$. We assume that in the alphabet $\Lambda_n$,  each $P_{n+1}$ is uniquely readable.  
  \begin{definition}A construction sequence $\la \mcw_n:n\in\nn\ra$  will be called \emph{circular} if 
  it is built in this manner using the $\mcc$-operators, a circular coefficient sequence and each 
  $P_{n+1}$ satisfies the strong unique readability assumption.   \end{definition}

It follows from Lemma \ref{unique readability} that each $\mcw_n$  in a circular construction sequence is uniquely readable.
\begin{definition}\label{circular definition}
A symbolic shift $\bk$ built from a  circular construction sequence  will be called a \emph{circular system}. 
\end{definition}
For emphasis we will often write circular construction sequences as $\la \mcw_n^c:n\in\nn\ra$ and the associated circular  shift $\bk^c$. We sometimes write $w^c$ to emphasize that a word is a circular word.

\bigskip

We will need to analyze the words constructed by $\mcc$ in detail. We start by describing the boundary and interior portions of the words.

\begin{definition}
Suppose that $w=\mcc(w_0,w_1,\dots w_{k-1})$.  Then $w$ consists of blocks of $w_i$ repeated $l-1$ times, together with some $b$'s and $e$'s that are not in the $w_i$'s. The \emph{interior} of $w$ is the portion of $w$ in the $w_i$'s.
{The remainder of  $w$ consists of blocks of the form $b^{q-j_i}$ and $e^{j_i}$. We call this portion the  \emph{boundary} of $w$.} 

In a block of the form $w_j^{l-1}$ the first and last occurrences  of $w_j$ will be called the  \emph{boundary} occurrences of the block $w_j^{l-1}$. The other occurrences will be the \emph{interior} occurrences. 
\end{definition}
{While the boundary consists of sections of $w$ made up of $b$'s and $e$'s, not all $b$'s and 
$e$'s occurring in $w$ are in the boundary, as they may be part of a power $w_i^{l-1}$.}

 The boundary of $w$ constitutes a small portion of the word:

\begin{lemma}\label{stabilization of names 1} The proportion of the word $w$ written in
equation \ref{definition of C} that belongs to its boundary is $1/l$. Moreover the proportion of the word that is within $q$ letters of boundary of $w$ is $3/l$. 
\end{lemma}

The next lemma was proved in \cite{prequel} (Lemma 20).

\begin{lemma}\label{dealing with S} Let $\bk^c$ be a circular system and
 $\nu$ be a shift invariant measure on $\bk^c$. Then the following are equivalent: 
	{\begin{enumerate}
	\item $\nu$ has no atoms.   
	\item $\nu$ concentrates on  the collection of $s\in \bk^c$ such that $\{i:s(i)\notin \{b, e\}\}$ is 
unbounded in both $\poZ^-$ and $\poZ^+$.
	\item  $\nu$ concentrates on $S$.
	\end{enumerate}}

\end{lemma}

\begin{remark} Let $\bk^c$ be a circular system.
\begin{enumerate} \item There are only two invariant atomic measures, one concentrates on the constant ``$b$" sequence, the other on the constant ``$e$" sequence.
\item for $\bk^c$,  Lemma \ref{principal blocks exist} can be strengthened to say that for all $s\in S$ for all large enough $n$, the principal $n$-block of $s$ exists.
\item The symbolic shift $\bk^c$ has zero topological entropy.
\end{enumerate}
\end{remark}
\pf
 A direct inspection reveals that the only periodic points in $\bk^c$ are the two fixed points constant ``$b$" and ``$e$".

The second item follows because if $s$ has a principal $n$-block at $[a_n, b_n)$ then it has a principal $n+1$-block at some $[a_{n+1}, a_{n+1}+q_{n+1})$ for an $a_{n+1}$ with $|a_{n+1}|\le |a_{n}|+(q_{n+1}-q_n)$.

 The fact that the topological entropy of $\bk^c$ is  zero follows easily from the fact that the $l_n$ tend to infinity.

\subsection{The structure of the words}\label{word sections}
The words used to form circular transformations have quite specific combinatorial properties.
We begin with an important definition for our understanding of rotations; the three \emph{subscales} at stage $n+1$. Fix a sequence $\la \mcw^c_n:n\in \nn\ra$ defining a circular system. { Using equation \ref{definition of C} we define the \emph{subscales} of a word $w^*\in \mcw_{n+1}$:} 
\begin{enumerate}
\item[]{\bf Subscale 0}  is the scale of the individual powers of $w_j\in \mcw^c_n$ of the form $w_j^{l-1}$; we call each such occurrence of a $w_j^{l-1}$ a \emph{0-subsection}
\item[]{\bf Subscale 1} is the scale of each term in the product $\prod_{j=0}^{k-1}(b^{q-j_i}w_j^{l-1}e^{j_i})$ that has the form $(b^{q-j_i}w_j^{l-1}e^{j_i})$; We call these terms \emph{1-subsections}.
\item[]{\bf Subscale 2} is the scale of each term of $\prod_{i=0}^{q-1}\prod_{j=0}^{k-1}(b^{q-j_i}w_j^{l-1}e^{j_i})$ that has the form 
$\prod_{j=0}^{k-1}(b^{q-j_i}w_j^{l-1}e^{j_i})$; We call these terms \emph{2-subsections}.
\end{enumerate}
\begin{center}
{\bf Summary}
\end{center}

\begin{center}
  \begin{tabular}{| l | c |}
    \hline
  {\bf Whole Word:}  &$\prod_{i=0}^{q-1}\prod_{j=0}^{k-1}(b^{q-j_i}w_j^{l-1}e^{j_i})$\\ \hline
{\bf 2-subsection:} &$ \prod_{j=0}^{k-1}(b^{q-j_i}w_j^{l-1}e^{j_i})$ \\ \hline
{\bf 1-subsection:} &$(b^{q-j_i}w_j^{l-1}e^{j_i})$\\ \hline
{\bf 0-subsection:} & $w_j^{l-1}$\\
    \hline
  \end{tabular}
\end{center}
 By contrast we will discuss \emph{$n$-subwords}  of a word $w$. These will be subwords that lie in 
 $\mcw^c_n$, the $n^{th}$ stage of the construction sequence. We will use \emph{$n$-block} to mean the location of the $n$-subword.

\subsection{The canonical circle factor $\mck$} \label{iso to rotation}
 
 We now define a canonical factor $\mck$ of a circular system and show that this factor is isomorphic to a rotation of the circle by $\alpha$, where $\alpha$ 
 is the limit of $\alpha_n= {p_n\over q_n}$ as $n$ goes to infinity.

\begin{definition}\label{first appearance of circle factor} Let $\la k_n, l_n:n\in\nn\ra\ra$ be a circular coefficient sequence. Let $\Sigma_0=\{*\}$. We define a  circular construction sequence such that each 
$\mcw^c_n$ has a unique element as follows:
\begin{enumerate}
\item $\mcw_0=\{*\}$ and
\item If $\mcw^c_n=\{w_n\}$ then $\mcw^c_{n+1}=\{\mcc(w_n, w_n, \dots w_n)\}$.

\end{enumerate}
Let $\mck$ be the resulting circular system. \end{definition}
It is easy to check that $\mck$ has unique ergodic non-atomic measure, since every $w_n$ occurs exactly $k_n(l_n-1)q_n$ many times in $w_{n+1}$.
\medskip

Let $\bk^c$ be an arbitrary circular system with  coefficients $\la k_n, l_n \ra$. Then $\bk^c$ has a canonical factor isomorphic to $\mck$. This canonical factor plays a role for circular systems analogous to the role  odometer transformations play for odometer based systems.

To see $\mck$ is a factor of $\bk^c$,  we define the following function:

\begin{equation}\label{definition of factor map}
\pi(x)(i) = \left\{ \begin{array}{ll}x(i) &
\mbox{if $x(i)\in \{b, e\}$} \\
* &\mbox{otherwise}
\end{array}\right.
\end{equation}

We record the following easy lemma that justifies the terminology of Definition \ref{first appearance of circle factor}:

\begin{lemma}\label{canonical rotation factor} Let $\pi$ be defined by equation \ref{definition of factor map}. Then:
\begin{enumerate}
\item $\pi:\bk^c\to \mck$ is a Lipshitz  map,
\item $\pi(sh^{\pm 1}(x))=sh^{\pm 1}(\pi(x))$ and thus
\item \label{last item} $\pi$ is a factor map of $\bk^c$ to $\mck$ and $(\bk^c)^{-1}$ to $\mck^{-1}$
\end{enumerate}
\end{lemma}
A variant of item \ref{last item} is also true: $\pi$ can be interpreted as a function from $\rev{\bk^c}$ to $\rev{\mck}$. With this interpretation $\pi$ is also a factor map. We will call $\mck$ the \emph{circle factor} of any circular system with construction coefficients $\la k_n, l_n:n\in \nn\ra$.

 Fix a circular coefficient sequence $\la k_n, l_n:n\in\nn\ra$, and let 
  $\mathcal K$ and $\la \mcw_n^\alpha:n\in\nn\ra$ be given in definition \ref{first appearance of circle factor}.
  Let $\alpha_n=p_n/q_n$ and $\alpha=\lim \alpha_n$.

If $s\in S$, from  $r_n(s)$ we can determine the locations of the beginnings and ends of the words $\wa_n$ that contain $s(0)$.
Since $|\mcw_n^\alpha|=1$ for all $n$,   for all $s\in S$ the sequence $\la r_n(s):n\in\mathbb N\ra$ uniquely determines $s$.

\begin{theorem}\label{rank one description}
Let $\nu$ be the unique non-atomic shift invariant measure on $\mck$. Then 
\[(\mck, \mcb, \nu, sh)\cong (S^1, \mathcal D, \lambda, \mcr_\alpha)\]
 where $\mcr_\alpha$ is the rotation of the unit circle by $\alpha$ and $\mathcal{B, D}$ are the $\sigma$-algebras of measurable sets.\end{theorem}
\pf A more involved geometric proof of this fact is given in \cite{prequel}. Here present a simple algebraic proof. 
As usual we identify the unit circle $S^1$ with $[0,1)$ and use additive notation for the group operations.

By Lemma \ref{principal blocks exist}, the collection $S'$ of $s\in S$ such that for all large enough $n$, the principal 
$n$-block of $s$ exists, has measure one. We define a map ${\phi}_0:S'\to \zoo$ by a limiting process. For $s$ such that $r_n(s)$ exists, we let 
{\begin{equation}\notag
{\rho}_n(s)={p\over q_n}
\end{equation}}
iff 
	\begin{equation}
	p\equiv p_nr_n(s)\mod{q_n}\notag
	\end{equation}
\begin{claim}\label{approximation}
If $r_n$ is defined, then $|{\rho}_{n+1}(s)-{\rho}_n(s)|<2/q_{n}$. 
\end{claim}
\pf
From equation \ref{definition of C}, we see that the position of $s(0)$ in an $n+1$-block is determined by the parameters $i\in [0,q_n-1), j\in [0, k_n-1), l^*\in[0,l-1]$ and $r_n$, which determine its location among  the 2-subsections, 1-subsections, 0-subsections and inside the $n$-words $w_n$ respectively. Explicitly:
\[r_{n+1}(s)=i(k_nl_nq_n)+j(l_nq_n)+(q_n-j_i)+l^*q_n+r_n(s),\]
where $r_n(s)$ is the position of $s(0)$ in its principal $w_n$-word. 

From the definition of $\rho_{n+1}$, and working \emph{mod 1}:
\begin{eqnarray*}
\rho_{n+1}&=&r_{n+1}(s)\left({p_{n+1}\over q_{n+1}}\right)\\
&=&r_{n+1}(s)\left({p_n\over q_n}+{1\over q_{n+1}}\right)
\end{eqnarray*}
Expanding this, using our formula for $r_{n+1}(s)$ and the fact that all but two terms of $r_{n+1}(s)$ are divisible by $q_n$,  we get:
\begin{eqnarray}
\rho_{n+1}&=&\left(-j_i\left({p_n\over q_n}\right)+ r_n(s)\left({p_n\over q_n}\right)\right) +\left({i\over q_n} + \delta\right)\label{alg exer}
\end{eqnarray}
where 
\[\delta={j\over k_nq_n}+{1\over k_nl_nq_n}+{l^*\over k_nl_nq_n}+{r_n(s)-j_i\over k_nl_nq_n^2}.\]
The first and third terms of equation \ref{alg exer} cancel, thus:
\[\rho_{n+1}=\rho_n+\delta.\]
Since $\delta<2/q_n$, the claim follows.\qed

\bigskip
Since the sequence $1/q_n$ is summable, for almost all $s, \la {\rho}_n(s):n\in \omega\ra$ is Cauchy. We define 
\begin{equation}
\phi_0(s)=\lim_n{\rho}_n(s).\notag
\end{equation} 
It is easy to check that $\phi_0$ is one-to-one.
By the unique ergodicity of the rotation $\mcr_\alpha$,  Theorem \ref{rank one description} will be proved when we establish:
\begin{claim}
The map $\phi_0:S\to [0,1)$ satisfies:
\begin{equation}
\phi_0(sh(s))=\mcr_\alpha(\phi_0(s)).\notag
\end{equation} In particular, if $\nu$ is the unique invariant measure on $S$
\begin{equation}
\label{rotation iso}
(\mck, \mcc, \nu, sh)\cong ([0,1), \mcb, \lambda, \mcr_\alpha).\notag
\end{equation}
\end{claim}

\pf Suppose that $r_n(s)$ and $r_n(sh(s))$ both exist. Then $r_n(sh(s))=r_n(s)+1$. If follows that $\rho_n(sh(s))=\rho_n(s)+p_n/q_n$. Taking limits we see that $\phi_0(sh(s))=\phi_0(s)+\lim_n \alpha_n=\phi_0(s)+\alpha$. \qed
\noindent This finishes the proof of Theorem \ref{rank one description}.\qed

{\subsection{Kronecker Factors}	
Both odometer transformations and irrational rotations of the circle are ergodic discrete spectrum transformations. Because  the odometer transformation based on $\la k_n:n\in\nn\ra$ is a factor of any odometer based system $T$ and  the rotation $\mcr_\alpha$ is a factor of any circular system $S$, both are factors of the respective Kronecker factors of $T$ or $S$.  
In general it is not the whole Kronecker factor in either case. }

We make the following lemma explicit in the case of odometer based transformations. In the case of systems with a circle factor the exactly analogous results hold.
\begin{lemma}\label{isos induce isos of odometer factor} Let $(\bk,\mcb,\mu, T)$ and $(\bl, \mcc, \nu, S)$ be measure preserving systems.
Suppose that $\bk$ has an  odometer factor $\mco$ and that $\phi:\bk\to \bl$ is an isomorphism. Then there is a unique odometer factor $\mco^*$ of $\bl$ with an isomorphism $\phi^\pi:\mco\to\mco^*$  such that the following diagram commutes:

\[
\begin{diagram}
	\node{\bk}\arrow{e,r}{\phi}\arrow{s,r}{\pi^\bk}\node{\bl}\arrow{s,r}{\pi^\bl}\\
	\node{\mco} \arrow{e,t}{\phi^\pi}\node{\mco^*}
\end{diagram}
\]

If each finite order eigenvalue of $\bl$ has multiplicity 1 (e.g. if $\bl$ is ergodic), then  $\mco^*$  is the unique odometer factor of $\bl$ isomorphic to $\mco$.
\end{lemma}
\pf 
Since the unitary operator $U_\phi:L^2(\bk)\to L^2(\bk)$ takes eigenfunctions to eigenfunctions, we know that 
$U_\phi$ takes the subspaces of $L^2(\bk)$ corresponding to $\mco$ to a subspace of $L^2(\bl)$ corresponding to an isomorphic copy of $\mco$. The lemma follows.\qed
An immediate corollary of Lemma \ref{isos induce isos of odometer factor} is that if $\bk$ and $\bl$ are ergodic odometer based systems over the same odometer $\mco$, with projections $\pi_K$ and $\pi_L$, then $\phi^\pi$ is an isomorphism between  the canonical odometer factors.
\medskip

We record the following consequences for later use; 
\begin{prop}\label{preservation jazz}
Suppose that $\bk$ and $\bl$ are  both ergodic odometer based systems with coefficients $\la k_n:n\in\nn\ra$. Then any isomorphism $\phi:\bk\to\bl$ takes the canonical odometer factor $\mco^\bk$ of 
$\bk$ to the canonical odometer factor $\mco^\bl$ of $\bl$. 

Similarly if $\bk^c$ and $\bl^c$ are both ergodic circular systems with the same coefficient sequences $\la k_n, l_n:n\in\nn\ra\ra$, then any isomorphism between $\bk^c$ and $\bl^c$ takes the canonical rotation $\mck^\bk$ to the canonical rotation factor $\mck^\bl$ \end{prop}
\pf In the first case there is a unique factor of $\bk$ and $\bl$ corresponding to the eigenvalues of 
$\mco^\bk$ and $\mco^\bl$. Any isomorphism must preserve the factor corresponding to these eigenvalues. The same argument works for $\mck$, as it is isomorphic to the rotation by $\alpha=\lim_n p_n/q_n$. \qed

\subsection{Uniform Systems}
In \cite{prequel} it is established that  the \emph{strongly uniform} circular systems with sufficiently fast growing $\la l_n:n\in\nn\ra$, are realizable as measure preserving diffeomorphisms of the torus. Strongly uniform systems are those for which each word in $\mcw_n$ occurs  the same number of times in each word in $\mcw_{n+1}$. These systems carry unique non-atomic invariant measures, simplifying much of what we do later in this paper. For example the correspondence between the measures $\nu$ on uniform odometer systems $\bk$ and $\nu^c$ on their  uniform circular system counterparts $\bk^c$ given in equation \ref{nu vs nuc}, is automatic. 

In the forthcoming \cite{part4} we show that arbitrary (i.e. non-uniform) circular systems are realizable as measure preserving diffeomorphisms of the torus, provided that the measures of the words in $\mcw_n$ go to zero.

 \section{Details of Circular Systems}\label{CS 2}
This section examines the circular systems defined in section \ref{circular systems 1} in more detail.  Initially we are given a circular coefficient sequence $\la k_n, l_n:n\in\nn\ra\ra$  and $\la q_n:n\in\nn\ra$ where  $q_n$ satisfies the inductive definition in  equation \ref{qn}. When $n$ is fixed, we again let $j_i=(p_n)^{-1}i$ modulo $q_n$ and $0\le j_i<q_n$. {Without significant loss of generality it is convenient to assume that $\sum 1/q_n<1/10$.}

 To understand joinings of circular systems we will be comparing generic elements $(s,t)$ of  circular 
 $\bk^c$ and $\bl^c$, and their parsings into subwords. We will use the following terminology:
 
 \begin{definition}
  Let $u, v$ be finite sequences of elements of $\Sigma\cup \{b, e\}$ having length $q$.  Given intervals $I$ 
  and $J$ in $\poZ$ of length $q$ we can view $u$ and $v$ as functions having domain $I$ and $J$ 
  respectively. We will say that $u$ is \emph{shifted by $k$} relative to $v$ iff $I$ is the shift of the interval 
  $J$ by $k$.  We say that $u$ is the \emph{$k$-shift} of $v$ iff $u$ and $v$ are the same words and $I$ is 
  the shift of the interval $j$ by $k$.
 \end{definition}

\subsection{Understanding the words}
\label{understanding the words}

We elaborate on the descriptions given in Section \ref{word sections}.
Our first combinatorial lemma is the following:

\begin{lemma}\label{gap calculation}Let $w=\mcc(w_0, \dots w_{k_n-1})$ 
for some $n$ and $q=q_n, k=k_n, l=l_n$. View $w$ as a word in the alphabet  $\Sigma\cup\{b, e\}$ lying on the interval of integers $[0, klq^2)$. 
\begin{enumerate}
\item If $m_0$ and $m_1$ are the locations of the beginnings of $0$-subsections in the same 2-subsection, then $m_0\equiv_qm_1$.
\item If $m_0$ and $m_1$ are such that $m_0$ is  the location of the beginning of a $0$-subsection  occurring in a $2$-subsection 
$\prod_{j=0}^{k-1}(b^{q-j_i}w_j^{l-1}e^{j_i})$ and $m_1$ at the 
i beginning of a $0$-subsection  occurring in the next 2-subsection $\prod_{j=0}^{k-1}(b^{q-j_{i+1}}w_j^{l-1}e^{j_{i+1}})$ then $m_1-m_0\equiv_q -j_1$.
\end{enumerate}
\end{lemma}
\pf To see the first point, the indices of the beginnings of $0$-subsections in the same $2$-subsection  differ by multiples of $q$ coming from powers of a $w_j$ and intervals of $w$ of the form $b^{q-j_i}e^{j_i}$.

To see the second point, let $u$ and $v$ be consecutive $2$-subsections. In view of the first point it suffices to consider the last $0$-subsection of $u$ and the first $0$-subsection of $v$.  But these sit on either side of an interval of the form $e^{j_i}b^{q-j_{i+1}}$. Since
$j_i+q-j_{i+1}\equiv_q (p)^{-1}i-p^{-1}(i+1)\equiv_q-p^{-1}\equiv_q-j_1$, we see that $m_0-m_1\equiv_q q+j_i+q-j_{i+1}\equiv_q-j_1$.
\qed

Assume that $u\in \mcw_{n+1}$ and $v\in \mcw_{n+1}\cup\rev{\mcw_{n+1}}$ and $v$ is shifted with respect to $u$.
 On the overlap of $u$ and $v$, the 2-subsections of $u$ split each 2-subsection of $v$ into either one or two pieces. Since all of the 2-subsections in both words have the same length, the number of pieces in the splitting  and the size of each piece is constant across the overlap except perhaps at the two ends of the overlap. If $u$ splits a 2-subsection of $v$ into two pieces, then we call the left  piece of the  pair  the even piece and the right piece the odd piece.

If $v$ is shifted only slightly, it can happen that either the even piece or the odd piece does not contain a $1$-subsection. In this case we will say that split is \emph{trivial on the left} or \emph{trivial on the right}

\begin{lemma}\label{numerology lemma}
Suppose that the $2$-subsections of $u$ divide the $2$-subsections of $v$ into two non-trivial pieces. Then
\begin{enumerate}
\item the boundary portion of $u$ occurring between each consecutive pair of 2-subsections of 
$u$ completely overlaps at most one $0$-subsection of $v$ 
 \item  there are two numbers $s$ and $t$ such that the positions of the $0$-subsections of $v$ in even pieces are  shifted relative to the  $0$-subsections of $u$ by $s$ and the positions of the  $0$-subsections of $v$ in  odd pieces are shifted relative to the $0$ subwords of $u$ by $t$.  Moreover $s\equiv_q t -j_1$.
\end{enumerate}
\end{lemma}

\pf This follows easily from Lemma \ref{gap calculation}\qed
In the case where  the split is trivial we get Lemma \ref{numerology lemma} with just one coefficient, $s$ or $t$.

A special case Lemma \ref{numerology lemma} that we will use is:

\begin{lemma}\label{weak numer}
Suppose that the $2$-subsections of $u$ divide the $2$-subsections of $v$ into two pieces and 
 that for some  
 occurrence of an $n$-subword  of $v$
 in an even (resp. odd) piece is lined up with an occurrence of some  $n$-word  
 in $u$.  Then 
 every occurrence of an $n$-word 
 in an even (resp. odd) piece of $v$ is either:
 \begin{enumerate}
 \item[a.)] lined up with some $n$-subword of $u$ or
 \item[b.)] lined up with a portion of a $2$-subsection that has the form $e^{j_i}b^{q-j_i}$.  
 \end{enumerate}
 Moreover, no $n$-subword in an odd (resp. even) piece of $v$ is lined up with a $n$-subword in $u$.
 \end{lemma}

\subsection{Full measure sets for circular systems}\label{full measure for ccs}

Fix a summable sequence $\la \varepsilon_n:n\in\nn\ra$ of numbers in $\zoo$ and a circular coefficient sequence $\la k_n, l_n:n\in\nn\ra$.
As we argued in the proof of Lemma \ref{stabilization of names 1},  the proportion of boundaries that occur in words of $\mcw^c_n$ is always summable, 
independently of the way we build $\mcw^c_n$. Recall the set $S\subseteq \bk^c$ given in Definition \ref{def of S}, where $\bk^c$ is the symbolic shift defined from a construction sequence.

\begin{definition}We define some  sets that a typical generic point for a  circular system eventually avoids. Let:
\begin{enumerate}
\item  $E_n$ be the collection of $s\in S$ such that $s$ does not have a principal $n$-block or $s(0)$ is in the boundary of that $n$-block, 
\item  $E^0_n=\{s:s(0)$ is in the first or last $\varepsilon_nl_n$ copies of $w$ in a power of the form $w^{l_n-1}$ where $w\in\mcw_n\}$, 
\item $E^1_n=\{s:s(0)$ is in the first or last $\varepsilon_nk_n$  1-subsections of the 2-subsection in which $s(0)$ is located$\}$, 
\item $E^2_n=\{s:s(0)$ is in the first or last $\varepsilon_nq_n$ 2-subsections of the principal $n+1$-block of$s\}$.
\end{enumerate}
\end{definition}

\begin{lemma}\label{bc1}  Assume that $\sum 1/l_n<\infty$. Let $\nu$ be a shift invariant measure on $S\subseteq \bk^c$, where $\bk^c$ is a circular system. Then:
\begin{enumerate} 
\item \[\sum_n\nu(E_n)<\infty.\] 
\begin{center}
Assume that $\la \varepsilon_n\ra$ is a summable sequence, then for $i= 0, 1, 2$:
\end{center}
\item
\[\sum_n\nu(E^i_n)<\infty.\]
\end{enumerate}
\end{lemma}
\pf This is an application of the Ergodic Theorem.\qed
In particular we see:
\begin{corollary}\label{bc2}
For $\nu$-almost all $s$ there is an $N=N(s)$ such that for all $n>N$, 
\begin{enumerate}
\item $s(0)$ is in the interior of its principal $n$-block,
\item $s\notin E^i_n$. 

In particular, for almost all $s$ and all large enough $n$:
\item if $s\rest [-r_n(s),-r_n(s)+q_n)=w$, then 
\begin{equation}\notag
s\rest[-r_n(s)-q_n, -r_n(s))=s\rest [-r_n(s)+q_n, -r_n+2q_n)=w.
\end{equation}

\item $s(0)$ is not in a string of the form $w_0^{l_n-1}$ or $w_{k_n-1}^{l_n-1}$.
\end{enumerate}
\end{corollary}
\pf This follows from the Borel-Cantelli Lemma.\qed
{The elements 
$s$ of $S$ such that   some shift $sh^k(s)$ fails one of the conclusions 1.)-4.) of Corollary \ref{bc2}  form a measure zero set.} Consequently we work on those elements of $S$ whose whole orbit satisfies the conclusions of Corollary \ref{bc2}. Note, however that the $N(sh^k(s))$ depends on the shift $k$.
\begin{definition}\label{mature}
We will call $n$ \emph{mature} for $s$ (or say that \emph{$s$ is mature at stage $n$}) iff $n$ is so large  that 
 $s\notin E_m \cup \bigcup_{0\le i\le 2}E^i_m$
 for all $m\ge n$.
\end{definition}
Thus if $s$ is mature at stage $n$ then for all $m>n$ the principal $m$-block of $s$ exists and conclusions 1-4 of Corollary \ref{bc2} hold.

Recall that in Section \ref{circular systems 1}, we defined a canonical factor of a circular system which we called the circle factor. Since the notion of maturity only involves the punctuation of the words involved, it is an easy remark that for all $s\in S$, $n$ is mature for $s$ just in case $n$ is mature for $\pi(s)$, where $\pi$ is the canonical factor map.
\bigskip

 {{For the following definition and lemma, we view $s\in S$ as a function with domain $\poZ$, and $s\in \mcw_n$ as a function with domain $[0,q_n)$ or, sometimes, an interval $[k, k+q_n)$. In each of these cases we use \emph{dom($s$)} to mean the domain of $s$.}}

{\begin{definition}\label{def bound} We will use the symbol $\boundary_n$ in multiple equivalent ways. If $s\in S$ or $s\in \mcw^c_m$ we define $\boundary_n=\boundary_n(s)$ to be the collection of $i$ such that $sh^i(s)(0)$ is in the boundary portion of an $n$-subword of $s$. This is well-defined by our unique readability lemma.  In the spatial context we will say that $s\in \boundary_{n}$ if $s(0)$ is the boundary of an $n$-subword of $s$. 
\end{definition}}
\noindent For $s\in S$
\[\boundary_n(s)\subseteq\bigcup\{[l,l+q_n):{l\in \mbox{ dom}(s)} \mbox{ and }s\rest[l,l+q_n)\in \mcw_n\}.\] 
An integer, $i\in \boundary_n(s)\subseteq \poZ$ iff $sh^i(s)$, viewed as an element of $\bk^c$, belongs to the $n$-boundary, $\boundary_n$.

 In what follows we will be considering a generic point $s$ and all of its shifts. We will use the fact if $s$ is mature at stage $n$, then we can detect  locally  those $i$ for which the $i$-shifts of $s$ are mature.  
\begin{lemma}\label{getting old}
Suppose that $s\in S$, $n$ is mature for $s$ and $n<m$.

\begin{enumerate}
 
\item Suppose that $i\in [-r_m(s), q_m-r_m(s))$. Then $n$ is mature for $sh^i(s)$ iff 
\begin{enumerate}
\item $i\notin \bigcup_{n\le k\le m} \boundary_k$ and
\item $sh^i(s)\notin \bigcup_{n\le k< m}(E^0_k\cup E^1_k\cup E^2_k)$.

\end{enumerate}
\item For all but at most $(\sum_{n< k\le m}1/l_k) +   (\sum_{n\le k< m}6\varepsilon_kq_{k+1})/q_m$ portion of the $i\in [r_m(s), q_m-r_m(s))$, the point $sh^i(s)$ is mature for $n$.

\end{enumerate}

In particular, if  $\varepsilon_{n-1}>sup_{m}({1/q_m})\sum_{k=n}^{m-1} 6\varepsilon_kq_{k+1}$, $1/l_{n-1}>\sum_{k=n}^\infty 1/l_k$  and $n$ is mature for $s$,  the upper density of those $i\in \poZ$ for which the $i$-shift of $s$ is not  mature for $n$ is less than $1/l_{n-1}+\varepsilon_{n-1}$.
\end{lemma}
Similarly:
\begin{lemma}\label{whole $n$-blocks 1} 
Suppose that $s\in S$ and $s$ has a principal $n$-block. Then $n$ is mature provided that $s\notin \bigcup_{n\le m}E^0_m\cup E^1_m\cup E^2_m$. In particular, if $n$ is mature for $s$ and $s$ is not in a boundary portion of its principal $n-1$-block or in $E^0_{n-1}\cup E^1_{n-1}\cup E^2_{n-1}$, then $n-1$ is mature for $s$.
\end{lemma}

\subsection{The $\natural$ map} 
\label{definition of natural}

 Proposition \ref{preservation jazz} implies that any  isomorphism $\phi$ between an ergodic 
 $(\bk^c,sh)$ and 
$(\bk^c,sh^{- 1})$ induces an isomorphism $\phi^\pi$ between $(\mck, sh)$ and 
$(\mck, sh^{- 1})$, where $\mck$ is the canonical circle factor. 
Because $(\mck, sh\inv)$ is canonically isomorphic with $(\rev{\mck},sh)$ 
(Proposition \ref{spinning})
 and $(\mck, sh)$ is isomorphic to the rotation $\mcr_\alpha$ of the circle, we see that 
 $(\rev{\mck},sh)$ is isomorphic to the rotation $\mcr_{-\alpha}$.
 
 We use a specific isomorphism $\natural:(\mck, sh)\to (\rev{\mck}, sh)$  as a benchmark for 
understanding of potential  maps $\phi:\bk^c\to\rev{\bk^c}$. If we view $\mck$ as a rotation $\mcr_\alpha$ of the unit circle by $\alpha$ radians one can view the transformation $\natural$  as a symbolic analogue of complex 
conjugation $z\mapsto \bar{z}$ on the unit circle, which is an isomorphism between $\mcr_\alpha$ and $\mcr_{-\alpha}$.  Copying 
$\natural$ over to a map on the unit circle gives an isomorphism $\phi$ between $\mcr_\alpha$ and $\mcr_{{-\alpha}}$. Such an isomorphism must be of the form 
\[\phi(z)=\bar{z}e^{2\pi i \beta}\]
for some $\beta$. It follows immediately from this characterization that $\natural$ is an involution, however for completeness we prove this directly (and symbolically) in  Proposition \ref{alas necessary}. 

As usual we find it more convenient to work on the unit interval $I=[0,1)$ rather than the unit circle. The complex conjugacy map $z\mapsto \bar{z}$
corresponds to the map $x\mapsto -x$ on $[0,1)$.

We begin by recalling from equation \ref{definition of C} the formula for a $w\in \mcw^c_{n+1}$ that is of the form $\mcc(w_0, \dots w_{k_n-1})$:
\begin{equation}\label{C again}
  w=\prod_{i=0}^{q-1}\prod_{j=0}^{k-1}(b^{q-j_i}w_j^{l-1}e^{j_i})
\end{equation}
where $q=q_n,k=k_n,  l=l_n$ and  $j_i\equiv_{q_n}(p_n)^{-1}i$ with $0\le j_i<q_n$. By examining this formula we see that 
\begin{equation*}
\rev{w}=\prod_{i=1}^{q}\prod_{j=1}^ke^{j_{q-i}}\rev{w_{k-j}}^{l-1}b^{q-j_{q-i}}.
\end{equation*}
Applying the identity in formula \ref{reverse numerology}, we see that this can be rewritten as\footnote{We take $j_q=0$.}
\begin{equation}\label{reverse mcc}
\rev{w}=\prod_{i=1}^{q}\prod_{j=1}^k(e^{q-j_i}\rev{w_{k-j}}^{l-1}b^{j_i}).
\end{equation}
We can reindex again and get another form of equation \ref{reverse mcc}:
\begin{equation}\label{foreshadow}
\rev{w}=\prod_{i=0}^{q-1}\prod_{j=0}^{k-1}(e^{q-j_{i+1}}\rev{w_{k-j-1}}^{l-1}b^{j_{i+1}}).
\end{equation}

We can now state the basic lemma about the way $w$ lines up with a shift of $\rev{w}$. 

\begin{lemma}\label{two steps left} Let $w\in \mcw^c_{n+1}$ and view $w$ as sitting at location $[0, q_{n+1})\subseteq \poZ$. 
Let $q=q_n$ and $k=k_n$. Consider $sh^{-j_1}(\rev{w})$ as being the word $\rev{w}$ in location $[j_1, q_{n+1}+j_1))\subseteq \poZ$. 
For all but at most $2kq$ of the occurrences of an $n$-subword $w_j$ of $w$ starting in a location $r\in [0, q_{n+1})$, the reversed 
word $\rev{w_{k-j-1}}$ occurs in $sh^{-j_1}(\rev{w})$ starting at $r$.
\end{lemma}

\pf The word $w$ starts with a block of $q$ $b$'s and then a block of $l-1$ copies of $w_0$, whereas $\rev{w}$ starts with a block of $q-j_1$ $e$'s followed by $l-1$ copies of $\rev{w_{k-1}}$. Hence if we shift $\rev{w}$ to the right by $j_1$ (to get $sh^{-j_1}(\rev{w})$) the first copy of $\rev{w_{k-1}}$ is aligned with the first copy of $w_0$ in $w$. Hence all of the copies of $\rev{w_{k-1}}$ in the first 1-subsection are aligned with the copies of $w_0$ in the first 1-subsection of $w$. 
Because the consecutive blocks of $b$'s and $e$'s (or $e$'s and $b$'s)  in the 2-subsections  add up to $q$ we see that every copy of $\rev{w_{k-j-1}}$ in the first 2-subsection of   $sh^{-j_1}(\rev{w})$ is aligned with with a copy of $w_j$.

We now argue as in Section \ref{understanding the words}.   At the end of each 2-subsection, $w$ has a block of $e$'s of length $j_i$, followed at the beginning of the next 2-subsection, by a block of $b$'s of length $q-j_{i+1}$. Together the $e$'s and $b$'s form a block of length $j_i+q-j_{i+1}$, which is equivalent mod($q$) to $-j_1$. Similarly the combined length of a block of $b$'s and $e$'s finishing and starting consecutive 2-subsections of $\rev{w}$ is equal to $-j_1$ mod($q$).

Both the beginning of the block of $e$'s ending the $k^{th}$ 2-subsection and the end of the block of $b$'s starting the $k+1^{st}$ 2-subsection are of distance less than $q$  from the location of the end of the $k^{th}$ 2-subsection. 
It follows from this  and the comments in the previous paragraph, that if $S_1$ and $S_2$ are consecutive 2-subsections of $w$ and $S_1'$ and $S_2'$ are the corresponding 2-subsections of $\rev{w}$ then  the beginning of the first occurrence of $\rev{w_{k-1}}$ in $S_2'$ is within $2q$ of the first occurrence of $w_0$ is $S_2$ and their locations are equivalent mod($q$). Hence inside the first 1-subsection, the 0-subsections are lined up except for at most $2$ copies of $w_0$.  This pattern is continued through $S_2$, giving at most $2k$ locations of $n$-blocks that are not aligned in $S_2$. 

Since there are less than $q$ 2-subsections with potential misalignments, the Lemma is proved.
\qed
 The next proposition gives a somewhat more detailed view into situation of Lemma \ref{two steps left}. 
 \begin{prop}\label{first slip} Let $w,w'\in \mcw^c_{n+1}$ and suppose that 
\begin{eqnarray*}
w=\mcc(v_0,v_1,\dots v_{k_n-1}) & \mbox{ and }&w'=\mcc(v_0', v_1',\dots v'_{k_n-1}).
\end{eqnarray*}
We look at the relative positions of $n$-words in $w$ and $sh^{-j_1}(\rev{w'})$.
\begin{enumerate}
\item Each occurrence of $v_i$ in $w$ is either lined up with an occurrence of $\rev{v'_{k_n-i-1}}$ or entirely lined up with a section of $\boundary_n$ inside $sh^{-j_1}(\rev{w'})$.
\item There is a number $C$ such that for all $i$ the number of occurrences of $v_i$ lined up with an occurrence of $\rev{v'_{k_n-i-1}}$ is $C$.
\end{enumerate}
\end{prop}
\pf
The first part is clear from the proof of Lemma \ref{two steps left}. The second part follows because all of the 1-subsections in a given 2-subsection of $w$ have the same alignment relative to $sh^{-j_1}(\rev{w'})$.\qed

Since the total number of occurrences of $n$-subwords in $klq$, the proportion of $n$-subwords lined up with $\boundary_n$ in $sh^{-j_i}(\rev{w'})$ is at most $2/l$.

Suppose that $\mck$ is given by the canonical construction sequence $\la \mcw_n^\alpha:n\in\nn\ra$. 
We define a sequence of functions $\la\Lambda_n:n\in\nn\ra$ and argue that they converge to an isomorphism from $\mck$ to $\rev{\mck}$. 

 We begin by defining an increasing sequence of natural numbers. Recall the definition of the Anosov-Katok coefficients $p_n$ and $q_n$ given in equations \ref{pn} and \ref{qn}.  Since $p_n$ and $q_n$ are relatively prime we can define $(p_n)^{-1}$ in $\poZ/q_n\poZ$. For the following definition we will view $(p_n)^{-1}$ as a \emph{natural number} with $0\le (p_n)^{-1}<q_n$.\footnote{In the notation used to define $\mcc$, $(p_n)\inv=j_1$. However  the notation $j_1$ is ambiguous (it depends on $n$), so  we use $(p_n)\inv$ in this context.}

 We let $A_0=0$ and 
\begin{equation}
A_{n+1}=A_n-(p_n)^{-1}. \label{code coefficients}
\end{equation}

\begin{lemma}\label{An is small}
If $A_n$ is defined as above, then $|A_{n+1}|<2q_{n}$.
\end{lemma}
\pf This is proved inductively using the fact that $q_{n+1}>2q_{n}$.\qed

Let $\mck$ be the  circular system in the language $\Sigma=\{*\}$, as given in Definition \ref{first appearance of circle factor}. We now define a stationary code $\overline{\Lambda}_n$ with domain $S$ that approximates elements of $\rev{\mck}$
 by  defining
\begin{equation}\label{definition of Lambdan}
\Lambda_n(s)=\left\{\begin{array}{ll} sh^{A_n+2r_n(s)-(q_n-1)}(\rev{s})(0) & \mbox{if $r_n(s)$ is defined}\\
				b&\mbox{otherwise}
				\end{array}\right.
\end{equation}

Since for all $s\in S$ and all large enough $n$, $r_n(s)$ is defined, the  default value is only obtained for finitely many $n$.
\begin{lemma}
${\Lambda}_n$ is given by a finite  code.
\end{lemma}
\pf To check whether $r_n(s)$ is defined one need only examine $s$ on the interval $[-q_n, q_n]\subseteq \poZ$. The relevant portion of $\rev{s}$ necessary to compute $\Lambda_n(s)$ is contained in $s\rest[-q_n-A_n, q_n+A_n]$. Hence $\Lambda_n$ is determined by a finite code.\qed

 The formula in equation \ref{definition of Lambdan} can be understood as follows.
Suppose that $s\in S$ and $s$ has a principal $n$-block. Then the element $s^*$ defined as
 $sh^{2r_n(s)-(q_n-1)}(\rev{s})$ belongs to $\rev{{\mck}}$, has a principal $n$-block that is the reverse of the principal $n$-block of $s$ and moreover, the principal $n$-block of $s^*$ is exactly lined up with the principal $n$-block of $s$. 

The reverse of the principal $n$-block of $s$ begins with a block of $q_{n-1}-(p_{n-1})^{-1}$ many $e$'s, 
and hence if $s'=sh^{(-(p_{n-1})^{-1})+2r_n(s)-(q_n-1)}(\rev{s})$ then the  first $n-1$-subword of the 
principal $n$-block of $s'$ is lined up with the first $n-1$-subword of the principal $n$-block of $s$. The rest of the terms used to define $A_n$ (coming from $A_{n-1}$) are used for lower order adjustments inside this principal $n$-block.

\medskip
\noindent Thus, a qualitative description of $\bar\Lambda_{n}(s)$ can be given as follows: 
\begin{enumerate}
	\item It first reverses the principal $n$-block of $s$ leaving it exactly lined up. 
	 \item It then adjusts the result by shifting so that the first occurrence of a 
	 reverse $n-1$-block lines up with the  first $n-1$-subword of the principal $n$-block of $s$. (So far we have 
	 described $sh^{(-(p_{n-1})^{-1})+2r_n(s)-(q_n-1)}(\rev{s})$.)
	 By Lemma \ref{two steps left},  we get a sequence where the principal $n$-block of $\Lambda_n(s)$ has the vast majority of its $n-1$-blocks lined up with the $n-1$-blocks of $s$: all of them except those that span a section of boundary at the juncture of two 2-subsections of the principal $n$-word of $s$.
	 
	 \item Finally it shifts by $A_{n-1}$ which is the cumulative adjustment at earlier stages.
 \end{enumerate}

The next lemma follows from this description: 

\begin{lemma} 
\label{88bis}
Let $n<m$ and suppose that $s\in \mck$ has a principal $m$-block. Let $s'=sh^{2r_m-q+A_m-A_n}(rev(s))$. Then at least
\[\prod_{n}^{m-1}(1-{2\over (l_i-1)})\]
proportion of the $n$-blocks in the principal $m$-block of $s$ are lined up with $n$-blocks in $s'$. \end{lemma}
\pf We first consider $m=n+1$. By Lemma \ref{two steps left}, all but $2k_nq_n$ of the $n$-blocks in $w$ are aligned with the $n$-blocks in $sh^{-j}(\rev{w})$.  This is proportion
\begin{equation*}
1-{2k_nq_n\over k_nq_n(l_n-1)}=1-{2\over l_n-1}.
\end{equation*}
The general result follows by induction.\qed

\begin{theorem}\label{mr natural}
Suppose that $\la k_n,l_n:n\in\nn\ra$ is a circular coefficient sequence. Then
 the sequence of stationary codes $\la \overline\Lambda_n:n\in\nn\ra$  converges to a 
 {shift invariant function} $\overline{\natural}:\mck\to (\{*\}\cup \{b, e\})^{{\poZ}}$  that induces an 
 isomorphism  $\natural$ from $\mck$ to $\rev{\mck}$.
\end{theorem}
\pf 
We first show that the sequence $\la \overline\Lambda_n:n\in\nn\ra$ converges, which will follow if we show  that the code distances between the $\Lambda_n$ and $\Lambda_{n+1}$ are summable. For notational simplicity, let $q=q_n, k=k_n, l=l_n$ and $j\equiv_{q}(p_n)^{-1}$ with $0\le j<q$.
\bigskip

\bfni{Claim:} There is a summable sequence of positive numbers $\delta_n$ such that for almost all $s$, the $\dbar$-distance between $\bar{\Lambda}_n(s)$ and $\bar\Lambda_{n+1}(s)$ is bounded by $\delta_n$, and $\bar\Lambda_n(s)$ and 
$\bar\Lambda_{n+1}(s)$ agree on all but at most $\delta_n$ proportion of the $n$-blocks of $s$.
\medskip

We use Lemma \ref{computing code distances}, which tells us that for a typical $s\in S$, the code distance between $\Lambda_n$ and $\Lambda_{n+1}$ is 
$\dbar(\overline\Lambda_n(s),\overline\Lambda_{n+1}(s))$, which is defined to be   the density of
\begin{equation}\label{real distance}
D=_{def}\{k:\Lambda_n(sh^k(s))(0)\ne \Lambda_{n+1}(sh^k(s))(0)\}.
\end{equation}

{Because  $|\mcw_n^\alpha=1$
for each $n$, there is only one possible $n$-subword at any location of any element of 
$\rev{\mck}$. Thus to compute $\dbar$-distance,  it suffices count positions where the   $\overline\Lambda_m$'s disagree on the  \emph{locations} of the $n$-subwords.}

 By Lemma \ref{getting old} for a typical $s\in S\subseteq\mck$ and all $n$, $I_{n}=_{def}\{i:n $ is not mature for $sh^i(s)\}$ has density at most $1/l_{n-1}+\varepsilon_{n-1}$, hence we can neglect these $i$ when computing the density of $D$. 

This allows us to assume that $r_{n+1}(s)$ is defined. 
We compute the density of the difference between $\bar\Lambda_n$ and $\bar\Lambda_{n+1}$ as they pass across an $n+1$-block in $s$.  If this number is $d$ then 
the distance between $\Lambda_n$ and $\Lambda_{n+1}$ is bounded by the sum of $d$ and the density of $I_{n}$.

As  $\Lambda_{n+1}$ 
crosses an $n+1$-block it produces the reverse $n+1$-block shifted by $A_{n+1}$. Explicitly, if $w$ is the 
$n+1$-block of $s$, as $\Lambda_{n+1}$ crosses $w$ it produces $sh^{A_{n+1}}(\rev{w})$. As
 $\Lambda_n$ passes across this same section, each time it crosses an $n$-block $w'$ it produces
  $sh^{A_n}(\rev{w'})$. If $w'$ starts at $r$ then the beginning of this copy of 
$sh^{A_n}(\rev{w'})$ is $r-A_n$.

We begin by rewriting $sh^{A_{n+1}}(\rev{w})$ as $sh^{A_n}(sh^{-j}(\rev{w}))$ where $j=(p_{n})^{-1}$.
 By Lemma \ref{two steps left}, all but $2kq$ of the $n$-blocks in $w$ are aligned with the $n$-blocks in $sh^{-j}(\rev{w})$. Hence, relative to the complement of $I_n$,  the portion of the principal $n+1$-block $w$ of $s$ that lies in an  $n$-block aligned with an $n$-block of $sh^{-j}(\rev{w})$ is 
\begin{equation}\label{crude estimate}{k(l-1)q^2-2kq\over k(l-1)q^2}=1-{2\over (l-1)q}
\end{equation}
{Because there is only one possible $n$-word, } whenever $sh^{A_n}(\rev{w'})$ is aligned with $sh^{A_n}(sh^{-j}(\rev{w}))$ they are equal.

Putting this altogether, we see that $\Lambda_n$ and $\Lambda_{n+1}$ agree on all of the $n$-subwords of the principal $n+1$-block of $s$ that are aligned with $sh^{-j}(\rev{w})$. The disagreements are limited to the $n$-subwords that are not aligned and the boundary. The total length of the disagreements is therefore bounded by
\[(2kq)*q+kq^2=3kq^2.\]
This has proportion $3kq^2/klq^2=3/l$.

Thus the distance between $\Lambda_n$ and $\Lambda_{n+1}$ is bounded by 
  $1/l_{n-1}+\varepsilon_{n-1}+3/l_n$. In particular the distances are summable and the sequence $\la \bar\Lambda_n:n\in\nn\ra$ converges almost everywhere to a function
  {$\natural:\mck \to (\Sigma\cup\{ b, e\})^\poZ$}.  
\medskip

We now show that $\natural$ is an isomorphism between $\mck$ and $\rev{\mck}$. Since $\bar\Lambda_n$  takes an $n$-block to a shift of the reverse $n$-block, it makes sense to discuss the \emph{principal $n$-block} of $\bar\Lambda(s)$.  Since the $r_n$'s cohere as in Remark \ref{interval coherence}, for $n<m$, $r_m(\bar\Lambda_m(s))$ is in the $r_n(\bar\Lambda_m(s))^{th}$ position of the principal $n$-block of $\bar\Lambda_m(s)$ (provided both $r_n$ and $r_m$ are defined). An application of the Ergodic Theorem shows that if $D_n$ is defined to be the collection of $s$ such that:
\[r_n(\bar\Lambda_n(s))\mbox{ exists and the principal $n$-words of $\bar{\Lambda}_n(s)$ and $\bar\Lambda_{n+1}(s)$ disagree}\]
then $\sum\nu(D_n)<\infty$.   From the Borel-Cantelli Lemma, it follows that for almost every $s$ for all large enough $n$ the principal $n$-blocks of $\bar\Lambda_n(s)$ and $\bar\Lambda_{n+1}(s)$ are the same, and thus that for {$s\in S, \natural(s)\in \rev{\mck}$.}

We now argue that if $s$ is typical and $s^*=\natural(s)$, then $s^*\in \rev{S}$. 
It suffices to show that $\lim_{n\to \infty}-r_n(s^*)=-\infty$ and $\lim_{n\to \infty}q_n-r_n(s^*)=\infty$.\footnote{We are adopting the convention that in defining $r_n(s^*)$ for $s^*\in \rev{S}$ we count $r_n$ from the left end of an $n$-block. Thus the position $r$ in a word $w\in \mcw^\alpha_n$ corresponds to the position $q-1-r$ in $\rev{w}$.}

If $n$ is mature for $s$ and large enough that for $m>n, \bar{\Lambda}_m(s)$ and $\bar\Lambda_n(s)$ have the same principal $n$-blocks,  then $r_n(s^*)=r_n(s)+A_n$ unless $r_n(s)\in [0,|A_n|)$.  Assuming that $r_n(s)\ge |A_n|$, we know from Lemma \ref{An is small} that 
\[r_n(s)-2q_{n-1}<r_n^*(s)<r_n(s).\] 
Hence, $-r_n^*(s)\le 2q_{n-1}-r_n(s)$ and $q_n-r_n^*(s)\ge q_n-r_n(s)$.
Applying Lemma \ref{getting old} (using the fact that $\sum nq_{n-1}/q_n<\infty$, and hence $\sum |A_n|/q_n<\infty$) we see that for large $n$, $r_n(s)>|A_n|$ and that $r_n(s)-2q_n\to \infty$. Since $q_n-r_n(s)\to \infty$ we have shown that $s^*\in \rev{S}$.

 As noted before Theorem \ref{rank one description}, if $s\in S$ then $s$ is determined by any tail of the sequence $\la r_n(s):n\in \nn\ra$. In particular, if we know a tail of $\la r_n(s^*):n\in \nn\ra$ we can determine $s^{*}$. Since for large $n$, 
 $r_n(s^*)=r_n(s)+A_n$, $\natural$ is one-to-one on a set of measure one.

We can now conclude that $\natural$ is an isomorphism. It is shift invariant since it is a limit of stationary codes, it maps from $S$ to $\rev{S}$, and is one-to-one on a set of 
$\nu$-measure one. If we define a measure 
$\mu$ on the Borel sets of $\rev{\mck}$ by setting $\mu(A)=\nu(\natural^{-1}(A))$, then $\mu$ is a shift invariant, non-atomic measure on $\rev{S}$. Since $S$ is uniquely ergodic,  $\rev{S}$ is as well and thus $\mu$ must be equal to the unique invariant measure $\nu$.
We have shown that $\natural$ is an isomorphism between $\mck$ and $\rev{\mck}$.
\qed

\begin{definition}\label{def of natural}
We denote the limit of $\la \bar{\Lambda}_n:n\in\nn\ra$ by $\natural:\mck\to \rev{\mck}$.
\end{definition}

We describe the  qualitative behavior of $\natural$ in a remark that we will use later:

\begin{remark}\label{natural coding} There is a summable sequence $\la \delta_n\ra$ such that for all but $1-\delta_n$ measure of $s\in S\subseteq \mck$, there is an interval $I$ containing 0 in $\overline{\Lambda}_n(s)$ such that $s\rest I\in \mcw^\alpha_n$, and moreover $\overline{\Lambda}_{n+1}(s)$ and $\overline{\Lambda}_n(s)$ agree on this interval. It follows from the Borel-Cantelli Lemma that 
 for almost all $s$ and large enough $n$, 
$\natural(s)$ agrees with $\bar\Lambda_{n}(s)$ on the principal $n$-block of $s$. Thus for a typical $s$ and large enough $n$, the map $\natural$ reverses the principal $n$-block while keeping its location and then shifts it by $A_n$.
\end{remark}

As noted at the beginning of this section, the next proposition follows immediately from Theorem \ref{rank one description}, however we include a symbolic proof for completeness.
\begin{prop}\label{alas necessary}
The map $\natural$ is an involution.
\end{prop}
\pf It is immediate from the qualitative description of $\bar{\Lambda}_n$ given before Lemma \ref{88bis}, that each $\bar{\Lambda}_n$ is an involution. To see that $\natural^2$ is the identity, let $\epsilon>0$. We can choose an $m_0$ large enough that for all $m\ge m_0$, $\bar{\Lambda}_m$ and $\natural$ agree with $\bar{\Lambda}_{m_0}$ on all but $\epsilon$ proportion of the $m_0$-blocks and $\bigcup_{m_0+1}^\infty \boundary_k$ has measure $\epsilon*10^-6$. Then $\natural\circ\bar{\Lambda}_{m_0}$ is equal to the identity on a set of density at least $1-\epsilon$. Letting $\epsilon\to 0$ and  $m_0\to \infty$ completes the argument.
\qed

%

 \subsection{Synchronous and Anti-synchronous joinings}
 Every odometer based system has a built in metronome: its odometer factor defined in Lemma \ref{odometer factor}. Correspondingly circular systems can be timed by their canonical rotation factor defined in Lemma \ref{canonical rotation factor}. 
 
 Joinings between odometer based and circular systems may induce non-trivial automorphisms of the underlying timing structure.  To avoid this complication we restrict ourselves to synchronous and anti-synchronous joinings: those  which preserve or exactly reverse the underlying timing. We now make this idea precise.
 
  Both the odometer transformations and rotations of a circle have easily understood inverse transformations and the isomorphisms between transformations and their inverses  are given by the maps $x\mapsto -x$ and $\rev{}\circ\natural$ respectively.
If $\bk$ and $\bl$ are either odometer based or circular systems let $\bk^\pi$ and 
$\bl^\pi$ be the corresponding odometer or rotation systems on which they are based.
\begin{definition}
\begin{itemize}
\item Let $\bk$ and $\bl$ be odometer based systems with the same coefficient sequence, and $\rho$ a joining between $\bk$ and $\bl^{\pm1}$. Then $\rho$ is 
\emph{synchronous} if $\rho$ joins $\bk$ and $\bl$ and the projection of $\rho$ to a joining on $\bk^\pi\times \bl^\pi$ is the graph joining determined by the identity map (the diagonal joining of the odometer factors); $\rho$ is \emph{anti-synchronous} if $\rho$ is a joining of $\bk$ with $\bl^{-1}$ and its projection to 
$\bk^\pi\times (\bl^{-1})^\pi$ is the graph joining determined by the  map $x\mapsto -x$.
\item Let $\bk^c$ and $\bl^c$ be circular systems with the same coefficient sequence and $\rho$ a joining between $\bk^c$ and $(\bl^c)^{\pm 1}$. Then $\rho$ is 
\emph{synchronous} if $\rho$ joins $\bk^c$ and $\bl^c$ and the projection to a joining of $(\bk^c)^\pi$ with $(\bl^c)^\pi$ is the graph joining determined by
the identity map of $\mathcal K$ with $\mathcal L$, the underlying rotations; $\rho$ is \emph{anti-synchronous} if it is a 
joining of $\bk^c$ with $(\bl^c)^{-1}$ and projects to the graph joining determined by $\rev{}\circ\natural$ on $\mck\times \mathcal L^{-1}$.
\end{itemize}
\end{definition}
\medskip

\noindent There is always a synchronous joining of odometer systems with the same underlying timing factor $\mco$:  \begin{definition}
Suppose that $\bk$ and $\bl$ are based on $\mco$. Then the relatively independent joining of $\bk$ and $\bl$ over $\mco$ is a synchronous joining, which we will call the \emph{synchronous product joining}. The relatively independent joining of 
$\bk$ and $\bl^{-1}$ over the map $x\mapsto -x$ we will call the \emph{anti-synchronous product joining}. 
 We will use the same terminology for the independent joinings of circular systems over the identity and $\rev{}\circ\natural$.  
\end{definition}

\section{Building the Functor $\mcf$}\label{building functor}

The main result of this paper concerns two categories whose objects are odometer based systems and circular systems respectively. The morphisms in these categories will be graph joinings. We will show that there is a  functor taking odometer systems to circular systems that preserves the factor and conjugacy structure. In this section we focus on defining the function from odometer based systems to circular systems that underlies the functorial isomorphism between these categories.

We begin by defining a function from the odometer based symbolic shifts $\bk$ to the circular symbolic shifts $\bk^c$. After having done so we define $\mcf$ on the pairs $(\bk, \mu)$ where $\mu$ is an invariant measure on $\bk$. Finally we define $\mcf$ on synchronous and anti-synchronous graph joinings.

We will use the notation that $K_n=\prod_{i<n}k_i$. 
 Then the $K_n$'s are the lengths of the odometer based words in $\mcw_n$ and the $q_n$'s are the lengths of the circular words in $\mcw_n^c$.
\medskip

Except where otherwise stated we will assume that we are working with a fixed \hyperlink{circ coef}{circular coefficient sequence} $\la k_n, l_n:n\in\nn\ra$. 
\medskip

Let $\Sigma$ be a language and $\la \mcw_n:n\in\nn\ra$ be a   construction sequence  for an 
odometer based system with coefficients $\la k_n:n\in\nn\ra$.    Then for each $n$ the operation $\mcc_n$ is well-defined. We define a  construction sequence 
$\la \mcw_n^c:n\in\nn\ra$ and bijections 
$c_n:\mcw_n\to \mcw_n^c$ by induction as follows:

\begin{enumerate}
\item Let $\mcw^c_0=\Sigma$ and $c_0$ be the identity map.
\item Suppose that  $\mcw_n, \mcw_n^c$ and $c_n$ have already been defined. 
\[\mcw_{n+1}^c=\{\mcc_n(c_n(w_0),c_n(w_1), \dots c_n(w_{k_n-1})):w_0w_1\dots w_{k_n-1}\in \mcw_{n+1}\}.\]

Define the map $c_{n+1}$ by setting
 \[c_{n+1}(w_0w_1\dots w_{k_n-1})=\mcc_n(c_n(w_0),c_n(w_1), \dots c_n(w_{k_n-1})).\] 
 
 \end{enumerate}
We note in case 2 the \hyperlink{pwords}{{prewords}} are:
\[P_{n+1}=\{c_n(w_0)c_n(w_1)\dots c_n(w_{k_n-1}): w_0w_1\dots w_{k_n-1}\in \mcw_{n+1}\}.\]

\begin{definition}\label{def of functor}
Define a map $\mcf$ from the set  of odometer based systems (viewed as subshifts) to  circular systems (viewed as subshifts) as follows. Suppose that 
$\bk$ is built from a construction sequence $\la \mcw_n:n\in\nn\ra$. Define 
\[\mcf(\bk)=\bk^c\]
where $\bk^c$ has construction sequence $\la \mcw_n^c:n\in\nn\ra$.
\end{definition}

Suppose that $\bk^c$ is a circular system with coefficients 
$\la k_n, l_n:n\in\nn\ra$.  We can recursively recursively build 
functions $c_n\inv$ from words in $\Sigma\cup \{b,e\}$ to words in $\Sigma$. The 
result is a  
odometer based system $\la \mcw_n:n\in\nn\ra$ with coefficients $\la k_n:n\in\nn\ra$.\footnote{We are using the strong unique readability assumption on the $P_n$'s to see the unique readability of the words in the sequence $\la \mcw_n:n\in\nn\ra$.}

 If $\bk$ is the resulting odometer based system then 
$\mcf(\bk)=\bk^c$. Thus we see:
\begin{prop}\label{bijection from readability}
The map $\mcf$ is a bijection between odometer based symbolic systems with coefficients $\la k_n:n\in\nn\ra$ and circular symbolic systems with coefficients $\la k_n, l_n:n\in\nn\ra$.
\end{prop}
\pf That $\mcf$ is one-to-one follows from the unique readability of words occurring in the construction sequence $\la W_n:n\in\nn\ra$.\qed

\begin{remark}
It is clear from  Definition \ref{def of functor} that $\mcf$ preserves uniformity and strong uniformity (see \cite{prequel} for these notions). In fact it preserves much more: the simplex of non-atomic invariant measures, rank one transformations and so on. We verify much of this in this paper and more in the forthcoming \cite{part4}.
\end{remark}

To understand the correspondence between measures on $\bk$ and $\bk^c$ we will have to understand the structure of basic open intervals. Recall that we write $\la u\ra_L$ to mean the basic open interval of $\bk$ determined by $u$ sitting on the interval $[L, L+|u|)\subseteq \poZ$. Without the subscript $L$, $\la u \ra$ is shorthand for $\la u\ra_0$. We adopt the same conventions for $\bk^c$, that the subscripts correspond to the beginning of the sequence and without a subscript the sequence begins at zero.

\subsection{Genetic Markers}\label{GMs and coding}
To see that $\mathcal F$ can be extended to a map from invariant measures on odometer based systems to invariant measures on circular systems, we begin by recalling how to identify elements of a symbolic system. 
Suppose that $\la \mcw_n:n\in \nn\ra$ is a construction sequence for an odometer 
based transformation $\bk$. Let $\la\mcw_n^c:n\in\nn\ra$ be the corresponding 
circular construction sequence for $\bk^c$. By Lemma \ref{specifying elements}  to 
specify a typical $s\in \bk$ or  $s^c\in \bk^c$, it suffices to give a tail of the sequence 
of principal $n$-blocks $\la w_n(s):N\le n\in\nn\ra$ or $\la w_n^c(s^c):N\le n\in\nn\ra$ 
along with the locations $\la r_n(s):N\le n\ra$ or $\la r_n(s^c):N\le n\ra$.

\begin{definition}\label{def of gms}
Suppose that $u, v$ are words in $\mcw_n$ and $\mcw_{n+1}$ respectively and $u$ occurs as an $n$-subword of $v$ in a particular location. Viewing $v$ as a concatenation $w_0w_1\dots w_{n_k-1}$ of $n$-subwords, there is a $j$ such that $u=w_j$. Let $j_n^*=j$ and call $j_n^*$ the \emph{genetic marker} of $u$ in $v$. 

 Suppose that $u\in \mcw_n$ and $v\in \mcw_{n+k}$ and $u$ is an $n$-subword of $v$ occurring at a particular location. Then there is a sequence of words $u_n=u, u_{n+1}, \dots  u_{n+k-1}, u_{n+k}=v $ such that $u_i$ is a $n+i$-subword of $v$ at a definite location  and the location of $u$ in $v$ is inside $u_i$. Let $j_{n+i}^*$ be the genetic marker of $u_{n+i}$ inside $u_{n+i+1}$. We call the sequence
  $\vec{j}^*=\la j_n^*, j_{n+1}^*, \dots j_{n+k-1}^*\ra$ the \emph{genetic marker} of $u$ in $v$. If $\vec{j}^*$ is the genetic marker of some $n$-word inside and $m$-word, we will call it an $(n,m)$-genetic marker.
  \end{definition}
If $u$ occurs as a subword of $v$ then  the genetic marker $\la j^*_n, j^*_{n+1} \dots j^*_{n+k-1}\ra$ of that occurrence codes its location inside $v$.

Suppose that $s\in \bk$ has principal $n$-blocks $\la w_n:n\in\nn\ra$. Each $w_{n+1}$ is a concatenation of words $v_0v_1\dots v_{k_{n}-1}$. Let 
	\begin{equation}j'_n=_{def}{r_{n+1}(s)-r_n(s)\over K_n}
	\label{jns}
	\end{equation}
or equivalently
\begin{equation}
r_{n+1}(s)=r_n(s)+j_n'K_n \label{inductive rns}.
\end{equation}
Each $w_{n+1}$ is a concatenation of words $v_0v_1\dots v_{k_n-1}$, and we see that $s(0)$ belongs to $v_{j_n'}$. In particular, the genetic marker of $w_n$ inside $w_{n+k}$ is the sequence $\la j_n', j_{n+1}', \dots j_{n+k-1}'\ra$.

\bigskip
\noindent{\bfni{Genetic markers for regions of words in $\mcw_{n+k}^c$:}} In circular words, genetic markers code regions rather than subwords.
Given $u$ and $v$ as above, we can consider the construction of $c_{n+k}(v)$ starting with the collection 
$\{c_n(u):u$ is an $n$-subword of $v\}$.
 Each of the genetic markers  $\la j^*_{n}, j^*_{n+1}, \dots j^*_{n+k-1}\ra$ of a subword $u$ of $v$ 
 determines a \emph{region} of $n$-subwords of $c_{n+k}(v)$. More explicitly, in the first step of the construction we 
 put $u$ into the $(j^*_{n})^{th}$ argument of $\mathcal C_n$. At the next step we put the result into the 
 $j^*_{n+1}$ argument of $\mcc_{n+1}$ and so on. Thus we see that there are bijections between
 \begin{enumerate} 
 \item sequences  $\la j^*_{n}, j^*_{n+1}, \dots j^*_{n+k-1}\ra$ with $0\le j^*_m< k_m$,
 \item $n$-subwords $u$ of $v$,
 \item the regions of $v^c$ occupied by the  occurrences of powers $(u^c)^{l_n-1}$ where $u^c$ is the element of $\mcw_n^c$ determined by $\la j^*_{n}, j^*_{n+1}, \dots j^*_{n+k-1}\ra$.
\end{enumerate}
Thus genetic markers give the correspondence between the regions of $c_{n+k}(v)$ that are not in $\bigcup_{n< m\le n+k}\boundary_m$ and particular occurrences of an $n$-word $u$ in $v$.

The next lemma computes the number of occurrences of a $c_n(u)$ with a given genetic marker $\la j^*_{n}, j^*_{n+1}, \dots j^*_{n+k-1}\ra$ in $c_{n+k}(v)$.

\begin{lemma}\label{products of imagination}
Suppose that $u^c$ occurs in $v^c$ with genetic marker $\la j^*_{n}, j^*_{n+1}, \dots j^*_{n+k-1}\ra$. Then the 
number of occurrences of $u^c$ in $v^c$ with the same genetic marker $\la j^*_{n}, j^*_{n+1}, \dots j^*_{n
+k-1}\ra$ is
\begin{equation}\label{helical product}
\prod_n^{n+k-1}q_i(l_i-1).
\end{equation}
\end{lemma}
\pf Fix $m$ and $v^c\in \mcw_m^c$. We prove equation \ref{helical product} for $n=m-k$ by induction on $k\ge 1$. If $k=1$ then we have a single genetic marker $j^*_{m-1}$. By formula \ref{definition of C} for $\mcc_{m-1}$  we see that the $j_{m-1}^*$ argument occurs in $v^c$ exactly $q_n(l_n-1)$ times.

Suppose now that we know that formula \ref{helical product} holds for $k-1$. We show it for $k$. Let $n=m-k$ and $u^c$ be the $n$-subword of $v^c$  with genetic marker$\la j^*_{n}, j^*_{n+1}, \dots j^*_{n
+k-1}\ra$. Let $w^c$ be the subword of $v^c$ with genetic marker $\la  j^*_{n+1}, \dots j^*_{n
+k-1}\ra$. Then:
	\[
	|\{\mbox{occurrences of $u^c$ in $v^c$ with marker }\la j^*_n, j^*_{n+1}, \dots j^*_{n
	+k-1}\ra\}|\]
is equal to 
	\[|\{\mbox{occurrences of $u^c$ in $w^c$ with marker $j^*_n$}\}|\times \] 
	\[|\{\mbox{occurrences of $w^c$ in $v^c	$ with marker }\la j^*_{n+1}, \dots j^*_{n
	+k-1}\ra\}|
	\]
The lemma follows.\qed

Since particular $(n,m)$-genetic markers   $\la j^*_{n}, j^*_{n+1}, \dots j^*_{n+k-1}\ra$ correspond to powers of $u^c$'s that occur with the same multiplicity in $v^c$, independently of the marker, we see that for a given $u$ and $v$:
\begin{equation}\label{up or down}
{|\{\mbox{occurrences of $u$ in } v\}|\over |\{n\mbox{-subwords of }v\}|}=
{|\{\mbox{occurrences of $c_n(u)$ in } c_{n+k}(v)\}|\over |\{\mbox{circular $n$-subwords of }c_{n+k}(v)\}|}
\end{equation}
We can restate equation \ref{up or down} in the language of section \ref{sequences and points}. It says that 
\begin{equation}\label{all is equal}
\mbox{EmpDist}(v)(u)=\mbox{EmpDist}(c_{n+k}(v))(c_n(u)).
\end{equation}

\medskip

In particular, if we fix a set $S^*$  of genetic markers we can compare the number of occurrences of a word with genetic marker in $S^*$ in $v\in\mcw_{n+k}$ with the number of occurrences in the corresponding $v^c\in \mcw_{n+k}^c$. Specifically,  the number of occurrences of a word $u^c$ in $v^c$ at some genetic marker in $S^*$ is 
$|S^*|*\prod_n^{n+k-1}q_i(l_i-1)$. The proportion of $n$-words occurring with a genetic marker in $S^*$ relative to all $n$-words occurring in $v^c$ is the same as the proportion of $n$-words with genetic markers in $S^*$ occurring in 
$v$  relative to the total number of genetic markers.  The number of $(n,m)$-genetic markers is 
$\prod_n^{n+k-1}k_i$  so this proportion is equal to
\begin{equation}\label{keeping things in proportion}
{|S^*|\over \prod_{n}^{n+k-1}k_i}.
\end{equation}
This is simply a restatement of our discussion involving empirical distributions in Section \ref{sequences and points}.
\bigskip

We introduce some notation that allows us to compare densities of various sets between odometer based and circular words. For  sets $A\subseteq [0,K_m)$ and $A^c\subseteq[0,q_m)$ we denote their densities by:
\begin{eqnarray*}
d_m(A)&=&|A|/K_m\\
d_m^c(A^c)&=&|A^c|/q_m
\end{eqnarray*}
Then $d_m$ and $d_m^c$ can be viewed as discrete probability measures on the sets $[0,K_m)$ and $[0,q_m)$ respectively.

\hypertarget{71}{\begin{lemma} \label{dracula}
Let $n\le m$, $w\in \mcw_m$ and $w^c=_{def}c_m(w)\in \mcw^c_m$. {We view $w$ as sitting on the interval $[0,K_m)$ and $w^c$ as sitting on $[0,q_m)$} Let $S^*$ be a collection of $(n,m)$-genetic markers, $g$ the total number of $(n,m)$-genetic markers
 and $d=|S^*|/g$.} If:
\begin{itemize}
\item $A=\{k\in [0,K_m):$  some $u\in \mcw_n$ with genetic marker in $S^*$ begins at $k$ in $w\}$
\item $A^c=\{k\in [0, q_m):$ some $u^c\in \mcw_n^c$ with genetic marker in $S^*$ begins at $k$ in $w^c\}$,
\end{itemize}
then the following equations hold:
\begin{eqnarray}
d_m(A)&=&{d\over K_n}\label{second mess}\\
d_m^c(A^c)&=&{d\over q_n}\prod_{p=n}^{m-1}(1-1/l_p)\label{first mess}\\
d_m(A)&=&\left({d_m^c(A^c)\over \prod_{p=n}^{m-1}(1-1/l_p)}\right)\left({q_n\over K_n}\right)\label{third mess}\\
d_m^c(A^c)&=&d_m(A)\left(\prod_{p=n}^{m-1}(1-1/l_p)\right)\left({K_n\over q_n}\right).\label{fourth mess}
\end{eqnarray}

\end{lemma}

\pf
We prove  equation \ref{first mess}. Equation \ref{second mess} is similar but easier. The other two equations follow algebraically. 

{ The union of the  boundary regions $\boundary_p$ for $p=n$ to $m-1$ consist exactly of the elements of $[0,q_m)$ that are not part of any $n$-word. We denote the complement of  $\bigcup_{p=n}^{m-1}\boundary_p$ by $(\bigcup_{p=n}^{m-1}\boundary_p)\tilde{}$. The various $\boundary_p$ are pairwise disjoint and for each $n^*$, 
 $(\bigcup_{p=n^*}^{m-1}\boundary_p)\tilde{}$ 
 consists of the locations of entire $n^*$-words. Starting with $p=m-1$, iteratively deleting boundary sections as $p$ decreases to $n$, and using 
   Lemma \ref{stabilization of names 1} we see that the $d^c_m$-measure of $(\bigcup_{p=n}^{m-1}\boundary_p)\tilde{}$ is $\prod_{p=n}^{m-1}(1-1/l_p)$.}

{ Let $B=\{k\in [0,q_m):k$ is at the beginning of an $n$-word$\}$. Then $B$ consists of a $1/q_n$ portion of the regions made up of $n$-words; i.e. $(\bigcup_{p=n}^{m-1}\boundary_p)\tilde{}$. We note that $A^c\subseteq B$ and $B$ is disjoint from $\bigcup_{p=n}^{m-1}\boundary_p$.}

{
By Lemma \ref{products of imagination}, the number $C_1$ of $n$-words occurring in $w^c$ with a given genetic marker does not depend on the marker.  Let $C_2$ be the total number of $n$-words occurring in $w^c$. Then:
\begin{eqnarray*}  
{|A^c|\over |B|}&=&{|\{n\mbox{-words with genetic marker in }S^*\}\over C_2}\\
&=&{|S^*|*C_1\over g*C_1}\\
&=&d.
\end{eqnarray*}
We compute conditional expectations  to get equation \ref{first mess}:}

\begin{eqnarray*}
d_m^c(A^c)&=&d_m^c(A^c\ |\ (\bigcup_{p=n}^{m-1}\boundary_p)\tilde{}\ )\ d_m((\bigcup_{p=n}^{m-1}\boundary_p)\tilde{}\ )\\
&=&d_m^c(A^c\ |\ B,(\bigcup_{p=n}^{m-1}\boundary_p)\tilde{}\ {})\ d_m(B|\ (\bigcup_{p=n}^{m-1}\boundary_p)\tilde{}\ {})\ d_m((\bigcup_{p=n}^{m-1}\boundary_p)\tilde{}\ {})\\
&=&d\ \left({1\over q_n}\right) \prod_{p=n}^{m-1}(1-1/l_p)
\end{eqnarray*}
Equation \ref{second mess} is similar and \ref{third mess}, \ref{fourth mess} follow from the first two equations by substitution.
\qed
The following relationship between  pairs of measures $\nu$ on $\bk$ and $\nu^c$ on $\bk^c$
	\begin{equation}\notag
	\nu^c(\la c_n(u)\ra)=\left({K_n\over q_n}\right)\nu(\la u\ra)\left(1-\sum_n^\infty \nu^{c}(\boundary_m)\right)
	\end{equation}
is the limit of equation \ref{fourth mess} as $m$ goes to infinity. This relationship will hold for a correspondence between measures that we build in forthcoming sections.

We note that since $\boundary_m$ has a density that depends only on the circular coefficient sequence, the measures of $\boundary_m$ is the same for all invariant measures. If we set
$d^{\boundary_n}$ be this density, then we can rewrite the previous equation as:

	\begin{equation}	
\nu^c(\la c_n(u)\ra)=\left({K_n\over q_n}\right)\nu(\la u\ra)\left(1-\sum_n^\infty d^{\boundary_n}\right)\label{nu vs nuc}
	\end{equation}

A consequence of equation \ref{nu vs nuc} is that  for all basic open sets $u$, $\nu(\la u\ra)$ determines $\nu^c(\la c_n(u)\ra)$ and vice versa.

For counting arguments the following inequalities will be helpful.
\begin{lemma}\label{deep in the heart of things} Let $n$ be a number greater than $0$. Then there are constants $K_n^U, K_n^L$ between 0 and 1 such that for all $k>0$ and  $w^c\in \mcw^c_{n+k}$ and all collections $S^*$ of $(n,n+k)$-genetic markers,
	
if 
\[A^c=\{i:i \mbox{ is the location of a start of an $n$-subword of $w^c$ indexed in } S^*\}\]
then
\begin{equation}\label{computing measure}
K^L_n|S^*|\le \left({|A^c|\over q_{n+k}}\right)\left(\prod_{m=0}^{n+k-1}k_{m}\right)\le K^U_n|S^*|
\end{equation}	

\end{lemma}
\pf 
By equation \ref{helical product}
there are
\[|A^c|=|S^*|*\prod_{m=0}^{k-1}q_{n+m}(l_{n+m}-1)\]
many $i$
that occur at the beginning of occurrences of $n$-subwords with genetic markers in $S^*$. Since 
\[q_{n+k}=k_nl_nq_n^2\left(\prod_{m=1}^{k-1}k_{n+m}l_{n+m}q_{n+m}\right)\]
we have:
\begin{equation}
{|A^c|\over q_{n+k}}=|S^*|*\left({1\over q_n}\right)\left(\prod_{m=1}^{k-1}(1-{1\over l_{n+m}})\right)\left({1\over \prod_{m=0}^{k-1}k_{n+m}}\right).\notag
\end{equation}
Since the $\la 1/l_n\ra$ is a summable sequence, $\prod_{m=1}^{k-1}(1-{1\over l_{n+m}})$ converges as $k$ goes to $\infty$.  The inequality \ref{computing measure} follows.\qed

Since $K_{n+k}=\prod_{m=0}^{n+k-1}k_{m}$, inequality \ref{computing measure} can be rewritten as:
\begin{equation}\label{explicit}
K_n^L{|S^*|\over K_{n+k}}\le {|A^c| \over q_{n+k}}\le K_n^U{|S^*|\over K_{n+k}}
\end{equation}

\bigskip

\hypertarget{itms}{\bfni{Infinite genetic markers:}} Suppose that we are given a construction 
sequence $\la \mcw_n:n\in\nn\ra$ for an 
odometer based or circular system $\bk$,  $s\in S$ and an occurrence of an $n$-word $u$ in $s$.  Then 
we can inductively define an  infinite sequence of words $\la u_m:n\le m\in\nn\ra$,  letting 
$u_n=u$, and $u_{m+1}$ to be the $m+1$-subword of $s$ that contains $u_m$.  
For each $n<m$ we get a genetic marker $\la j_n^*, j_{n+1}^*, \dots j_{m-1}^*\ra$, and these 
cohere as $m$ goes to infinity. We define the \emph{infinite genetic marker} to be 
$\vec{j}^*=\la j^*_m:n\le m\in\nn\ra$.

 If an $n$-word $u$ occurs inside an occurrence of an $m$-word $v$ in $s$, then $v=u_m$. Thus  their infinite genetic markers agree on the tail $\la j_i^*:m\le i\in\nn\ra$.

As in Remark \ref{rebuilding}, if we are given a sequence of words $\la u_m:n\le m\ra$, with $u_m\in \mcw_m$,  and an infinite sequence $\la j_m:n\le m\ra$ such that the genetic marker $j_{m}$ denotes an instance of $u_m$ in $u_{m+1}$ then we can find an 
$s\in \bk$ with $\la u_m:m\ge n\ra$ as a tail of its principal subwords. If $\bk$ is odometer then $s$ is unique up to a shift of size less than or equal to $K_m$. A similar statement holds for circular systems.

\subsection{$TU$ and $UT$.}
To understand the relationships between $\bk$ and $\bk^c$, we define maps $TU:S\to S^c$ and $UT:S^c\to S$ where $S\subseteq \bk$ and $S^c\subseteq \bk^c$ are as in definition \ref{def of S}. The map $TU$ will be one-to-one but $UT$ will not,  in general it is continuum-to-one. Nevertheless  $UT\circ TU$ will be the identity map.

We begin by considering a element $s\in S$. Let $u_n$ be the principal $n$-subword of $s$.
 The sequence $\la u_n:n\in\nn\ra$ determines a sequence of circular words $\la u_n^c:n\in\nn\ra$ which we assemble to define $TU(s)$. Let $\vec{j}=\la j_n:n\in\nn\ra$ be the infinite genetic marker of $s(0)$. To describe $TU(s)$ completely we need to define $\la r^c_n:n\in\nn\ra$. Set $r^c_0=0$, and inductively define $r^c_{n+1}$ to be the $(r^c_n)^{th}$ position in the first occurrence of an $n$-word with genetic marker $j_n$ in $u^c_{n+1}$. \hypertarget{def of TU}{Set $TU(s)$} to be the element of $\bk^c$ with principal subwords $\la u_n^c:n\in\nn\ra$ and location sequence $\la r_n^c:n\in\nn\ra$.

We \hypertarget{def of UT}{define a map $UT$} that associates an element of 
$\bk$ to each element of $S^c$. 
 Given such an $s^c\in S^c$, let $\la u_n^c:n\ge N\ra$ be its sequence of principal $n$-subwords. For 
each $n\ge N, u^c_n$ occurs as $u_{j^*_n}$ in the preword corresponding to $u_{n+1}^c$.  Let $u_n=c_n^{-1}(u^c_n)$. Then the sequence of words $\la u_n:n\in\nn\ra$ and genetic markers $\la j_n^*:n\ge N\ra$ determine an element of $s\in\bk$ except for the location of 0 in the double ended sequence. (The sequence is double ended because $s\in S^c$.)

We determine this location arbitrarily in a manner that makes the sequence of $u_n$'s the principal $n$-blocks of $s$ ($n\ge N$) and the $j_n^*$ the sequence of genetic markers of these $n$-blocks.
Let $\bar{0}$ be a 
sequence of zeros of length $N$. Then $\bar{0}\cat \la j_n^*:n\ge N\ra$ is a well-defined member of the 
odometer $\mco$ associated with $\bk$. 
From equation  \ref{inductive rns}, $\bar{0}\cat \la j_n^*:n\ge N\ra$ determines a sequence $\la r_{n}:n\in \nn\ra$. Thus by Lemma \ref{specifying elements},  the pair $\la u_n:n\ge N\ra$ and $\bar{0}\cat \la j_n^*:n\ge N\ra$ determines a unique element $s$ of $ \bk$ which we will denote $UT(s^c)$. 
It is easy to check that $UT\circ TU=id$ and that for each $s\in S$, there is a perfect set of $s^c$ with $UT(s^c)=s$.

We can get more precise information about correspondences between $\bk$ and $\bk^c$ by noting that if we are given a sequence $\la u_n:n\in\nn\ra$ of principal subwords of an $s\in S$, 
the genetic markers $\la j_n:n\in\nn\ra$ define an element $s^c$ of $\bk^c$ up to a choices 
 $(s^c)^\pi\in \mck$. Specifically, suppose that $s^*\in \mck$ is such that the infinite genetic marker   of $s^*(0)$ is 
 $\la j_n:n\in\nn\ra$. Then there is an $s^c\in \mck^c$ that has a sequence of principal $n$-blocks $\la u_n^c:n\in\nn\ra$.

The following lemma will be useful for understanding joinings.

\begin{lemma} \label{wrapping translations} Let $s\in S$. Then $\{TU(sh^k(s)):k\in \poZ\}\subseteq \{sh^k(TU(s)):k\in \poZ\}$.
If
 $s\in S$, $s^c=TU(s)$ and $u\in \mcw_n$, then there is a canonical correspondence between occurrences of $u$ in $s$ and finite regions of $s^c$ where $u^c$ occurs. The occurrences of $u^c$ in these finite regions have the same infinite genetic marker 
$\la j_m:m>n\ra$ in $s^c$ as $u$ does in $s$.  

\end{lemma}
\pf Given an  $s\in S$ and a $k$, the shift $sh^k(s)$ and $s$ have a tail of the principal $n$-blocks $\la u_n:N\le n\ra$ in common. Moreover the genetic markers associated with this tail are the same for both $s$ and $sh^k(s)$. It follows that $TU(sh^k(s))$ is a shift of $TU(s)$.

We can describe the correspondence as follows. If $u$ occurs in $s$ at $k$, then $u$ is the principal $n$-word of $sh^k(s)$. Choose an $N$ so large that  some $N$-word $u^*$ is the principal $N$-word of both $s$ and $sh^k(s)$. 
Then $(u^*)^c$ is the principal $N$-block of $s^c$. Let $\vec{j}$ be the genetic marker of the occurrence of $u$ (at $k$) in $u^*$. The region of $s^c$ corresponding to this occurrence of $u$ is the collection of occurrences of $u^c$ with the genetic marker $\vec{j}$ in the principal $N$-block of $s^c$.\qed

\subsection{Transferring measures up and down, I}
In this section we develop the tool we need for lifting measures on $\bk$ to measures on $\bk^c$. This will also allow us  to 
establish a one-to-one correspondence between synchronous joinings on odometer systems  and synchronous joinings on the corresponding circular systems.
Throughout this section we will use $\pi$ to denote either the projection of an odometer based system to its canonical odometer factor or a circular system to its canonical circular factor.

 We begin with a proposition relating sequences of words in a construction sequence for an odometer based system to sequences of  words in a construction sequence for a circular system.

\begin{prop}\label{words tell the story}
Let $\la v_n:n\in\nn\ra$ be a sequence  with $v_n\in\mcw_n$. Let $v_n^c=c_n(v_n)$. Then:
	\begin{enumerate}
	\item $\la v_n:n\in\nn\ra$ is an ergodic sequence iff $\la v_n^c:n\in\nn\ra$ is an 
	ergodic sequence.
	\item \label{gen iff gen} $\la v_n:n\in\nn\ra$ is a generic sequence for a measure $\nu$ iff 
	$\la v_n^c:n\in\nn\ra$ is a generic sequence for a measure  $\nu^c$. In case either sequence is generic, the measures  $\nu$ and $\nu^c$ satisfy equation \ref{nu vs nuc}.	\end{enumerate}
\end{prop}

\pf Both parts follow immediately from the definitions using equations \ref{all is equal} and \ref{keeping things in proportion} to relate the frequencies of $k$-words $w\in \mcw_k$ in $n$-words $u\in \mcw_n$, for $k<n$ to the frequencies of $c_k(w)$ in the corresponding $c_n(u)$. Equation \ref{nu vs nuc} follows from the Ergodic Theorem and Lemma \ref{dracula}. \qed

We endow that collection of invariant measures on a symbolic system
 $(\bk, sh)$ with the weak* topology.
 
 \begin{theorem}\label{k and kc}
Let $\la \mcw_n:n\in\nn\ra$ be a uniquely readable construction sequence for an odometer based system $\bk$ and $\la \mcw_n^c:n\in\nn\ra$ be the associated circular construction sequence for 
$\bk^c$. Then there is a canonical affine homeomorphism $\nu\mapsto \nu^c$  between shift invariant measures $\nu$ concentrating on $\bk$ and non-atomic, shift invariant measures $\nu^c$  such that equation \ref{nu vs nuc} holds between $\nu$ and $\nu^c$. 
\end{theorem}

\pf By Proposition \ref{kiss your S goodbye} and Lemma \ref{dealing with S} we can assume that $\nu$  and $\nu^c$ concentrate on $S$ and $S^c$ respectively.

We begin by defining the correspondence for ergodic measures. 
Suppose that we are given an ergodic measure $\nu$ and we want to associate a measure $\nu^c$. Let $s\in S$ be a generic point for $(\bk, \nu)$. Let $\la v_n:n\in\nn\ra$ be the sequence of principal $n$-blocks of $s$. By Proposition \ref{generic sequences exist for ergodic} this sequence is generic for $\nu$.  By Proposition 
\ref{words tell the story}, if we let $v_n^c=c_n(v_n)$, then $\la v_n^c:n\in\nn\ra$ is an 
ergodic sequence. Let $\nu^c$ be the measure associated with 
$\la v_n^c:n\in \nn \ra$. Then $\nu^c$ is ergodic and equation \ref{nu vs nuc} holds by Proposition \ref{words tell the story}.

The other direction is similar, let $s^c\in S^c$ be generic for $\nu^c$. Propositions \ref{generic sequences exist for ergodic} and  \ref{words tell the story} imply that if $\la v_n^c:n\in\nn\ra$ is the sequence of principal $n$-blocks of $s^c$ and $v_n=c_n^{-1}(v^c_n)$, then $\la v_n:n\in\nn\ra$ is ergodic and generic for a measure  $\nu$. Again equation \ref{nu vs nuc} holds by Proposition \ref{words tell the story}.

Suppose now that $\nu$ is an arbitrary measure on $\bk$. Write the ergodic decomposition of $\nu$ as:
\[\nu=\int \nu_i d\mu(i).\]
We define $\nu^c$ by
\[\nu^c=\int \nu_i^c d\mu(i)\] 
which gives a corresponding measure on $\bk^c$. Since equation \ref{nu vs nuc} holds between corresponding ergodic components $\nu_i$ and $\nu_i^c$, it holds between $\nu$ and $\nu^c$. 

By the ergodic decomposition theorem the map $\nu\mapsto \nu^c$ is a surjection. Since the map is invertible, it is a bijection. The map is affine by construction.

It remains to show that it is a homeomorphism. To see that $\nu\mapsto \nu^c$ is weak* continuous 
it suffices to show that for all $\epsilon>0$ and $n\in \nn$ there is a $\delta$ and an $m$ such that  for all invariant $\mu, \nu$, if for all $u\in \mcw_m$ 
\[|\mu(\la u\ra)-\nu(\la u\ra)|<\delta\] we know that
for all $v\in \mcw_n$ we have
\[|\mu^c(\la v^c\ra)-\nu^c(\la v^c\ra)|<\epsilon.\]
But the equation \ref{nu vs nuc} easily implies this taking $m=n$ and 
\[\delta<\left({K_n\over q_n}\right)\left(1-\sum_n^\infty d^{\boundary_n}\right)*\epsilon/4.\]
The argument that the inverse is continuous is the same.
\qed
\begin{definition}We will call a pair $(\nu,\nu^c)$ constructed as in Theorem \ref{k and kc}
\hypertarget{corresponding measures}{\emph{corresponding measures}}.
\end{definition}

\begin{remark}\label{benjy's remark} It follows from Proposition \ref{words tell the story} that if $\nu$ and $\nu^c$ are corresponding measures on $\bk$ and $\bk^c$ and $s\in \bk$ is arbitrary then $s$ is generic for $\nu$ iff $TU(s)$ is generic for $\nu^c$. The point $s$ is generic just in case its sequence of principal subwords is generic for $\nu$. By item \ref{gen iff gen} of  Proposition \ref{words tell the story}, this holds just in case the sequence of principal subwords of $TU(s)$ is generic; i.e. $TU(s)$ is generic. 
\end{remark}

We can use Theorem \ref{k and kc} to characterize the possible simplexes of invariant measures for circular systems. By a theorem of Downarowicz (\cite{downar}, Theorem 5), every non-empty compact metrizable Choquet simplex  is affinely homeomorphic to the simplex of  invariant probability measures for a dyadic Toeplitz flow. Note that the space of invariant probability measures is always a compact Choquet simplex, hence this theorem is optimal.

Since Toeplitz flows are special cases of odometer based systems it follows immediately that every  non-empty compact metrizable Choquet simplex is affinely homeomorphic to the simplex of invariant measures of a 2-symbol odometer based system.

Let $K$ be a compact Choquet simplex and $\bk$ an odometer based system having its simplex of invariant probability measures affinely homeomorphic to $K$. Let $\bk^c$ be a circular system corresponding to an odometer based system $\bk$. Then the non-atomic measures on $\bk^c$  are a Choquet simplex isomorphic to $K$. There are two additional ergodic measures, the atomic measures concentrating on the constant ``$b$" sequence and on the constant ``$e$" sequence. These  two atomic measures are isolated among the ergodic measures. 
\medskip

In the forthcoming \cite{part4} we discuss the question of invariant measures further and show that $\mcf$ preserves several other properties, such as being \emph{rank one}.

\section{$\mathbf{P^-}, \mathbf{P^\natural}$, genetic markers and the $\natural$-map}

Our goal is to understand the structure of synchronous and anti-synchronous joinings between pairs of ergodic systems $(\bk, \bl^{\pm 1})$. 
We will use Theorem \ref{k and kc} to define a bijection between synchronous joinings of odometer based systems and synchronous joinings of circular systems.  This  is 
relatively easy: to a joining of $\bk$ with $\bl$ that projects to the identity we can directly associate 
an odometer system $(\bk,\bl)^\times$ with a measure $\nu$ such that the corresponding measure 
$\nu^c$ on $((\bk,\bl)^\times)^c$ can be identified with a measure on $\bk^c\times \bl^c$ that projects to the identity. We carry this construction out in detail in section \ref{the categories} and show that the map $\nu\mapsto \nu^c$ given by Theorem \ref{k and kc} gives a bijection between synchronous joinings of the two kinds of systems.

The situation for anti-synchronous joinings of $\bk$ and $\bl\inv$ is more complicated. In Lemma \ref{joining correspondence}, we remarked that the anti-synchronous joinings of $\bk$ and $\bl\inv$ can be identified with joinings of $\bk$ and $\rev{\bl}$ that concentrate on $\{(s,t):\pi s=\pi t\}$. Similarly we can identify the anti-synchronous joinings of $\bk^c$ and $(\bl^c)\inv$ with joinings of $\bk^c$ with $\rev{\bl^c}$ that concentrate on 
 $\{(s^c,t^c):\pi t^c=\natural(\pi s^c)\}$. We give notation for these sets:
 \begin{enumerate}
 \item Let $\mathbf{P^-}$ be the collection of anti-synchronous joinings $\rho$ of $\bk$ and $\bl\inv$. 
 \item Let $\mathbf{P^\natural}$ be the collection of anti-synchronous joinings $\rho^c$ of $\bk^c$ and $(\bl^c)\inv$. 
 \end{enumerate}
 
 To understand the relationship between $\mathbf{P^-}$ and $\mathbf{P^\natural}$ we need an analogue of Lemma \ref{dracula}, and the corresponding analogue of equation \ref{nu vs nuc}. We now describe the tools we use to do this. 
 
Fix construction sequences for $\la\mathcal U_n:n\in\nn\ra$ and $\la\mathcal V_n:n\in\nn\ra$ for $\bk$ and $\bl$ respectively based on $\la k_n: n\in\nn\ra$ and $\bk^c, \bl^c$ the corresponding circular systems based on $\la k_n, l_n:n\in\nn\ra$.

 Let $(s,t)$ be 
an arbitrary point in $\bk\times \bl$ with $\pi t=-\pi s$ and $s\in S^\bk, t\in S^\bl$.
Let $\la u_n:n\in\nn\ra$ and $\la v_n:n\in\nn\ra$ be 
the sequence of principal subwords of $s$ and $t$ 
respectively. If $s^c=TU(s)$ and $t^c=TU(t)$, then 
$\la u_n^c:n\in\nn\ra$ and $\la v_n^c:n\in\nn\ra$ are 
the sequences of principal subwords of $s^c$ and 
$t^c$.

 Let 
$x=\natural(\pi s^c)$. Then $x\in \rev{\mck}$ and set 
$r_n=r_n(x)$. 

\begin{definition}\label{t-hat}
Define 
$\hat{t}\in \rev{\bl^c}$ by taking $\la \rev{v_n^c}:n\in \nn\ra$  as its principal $n$-subword
sequence and $\la r_n:n\in \nn\ra$ as its 
location sequence.
\end{definition} 
 We will study the relationship between $\mathbf{P^-}$ and $\mathbf{P^\natural}$ via the function taking 
$(s,t)$ to $(s^c, \hat{t})$.

 \subsection{Genetic Markers revisited} 
 To understand the relationship between joinings $\rho$ in $\mathbf{P^-}$ and $\rho^c$ in $\mathbf{P^\natural}$ we need to take into account the manner that $\natural$ shifts the reverse of the second coordinate of a the image of a generic pair $(s,t)$ for $\bk\times \bl\inv$
and the interplay between the map 
$\natural$ and genetic markers. Let $n<m$.
 Suppose that $(u',\rev{v'})$ is a pair of $n$-words  coming from $\mcu_n\times\rev{\mcv_n}$ that occur  aligned inside $m$-words $(u, \rev{v})\in \mcu_m\times \rev{\mcv_m}$. If $u'$ and $\rev{v'}$ occur at the same location in $(u,\rev{v})$, then $\vec{j}_{u'}$ determines $\vec{j}_{v'}$ in the following way:
 
 \begin{quotation}
 \noindent for $n\le r<m$ we must have 
  	\begin{equation}(j_{v'})_r=k_r-(j_{u'})_r-1 \label{conjugal gm's}
	\end{equation}
 (where $\vec{j}_{u'}=(j_n, j_{n+1}, \dots j_{m-1})$). 
\end{quotation}
 \begin{definition}\label{conjugate pairs def} Let $(u',v')\in \mcw_n$ and $(u,v)\in \mcw_m$. Define the $(n,m)$-\emph{genetic marker of an occurrence of the pair} $(u',\rev{v'})$ in $(u,\rev{v})$ to be $(\vec{j}_{u'},\vec{j}_{v'})$ where $\vec{j}_{u'}$ is the genetic marker of $u'$ in $u$ and $\vec{j}_{v'}$ is the genetic marker of ${v'}$ in ${v}$.\footnote{Note that the genetic marker $\vec{j}_{u'}$ denotes a different position inside $\rev{v}$ then it does in 
 $u$.} We call $\vec{j}_{u'}$ and $\vec{j}_{v'}$ a \emph{conjugate pair}. 
 \end{definition}
Being a conjugate pair is equivalent to satisfying the numerical relationship given in 
equation \ref{conjugal gm's} and thus  either element of a conjugate pair determines the other. Hence for purposes of counting conjugate pairs we need only use the first coordinates, $j_{u'}$.

Let $(u,\rev{v})\in \mcu_m\times \rev{\mcv_m}$ be words that occur in a pair $(s,\rev{t})\in \bk\times\rev{\bl}$.
Then the relative alignment of $u^c$ and 
$\rev{v^c}$ in  $(s^c,\hat{t}\ )$  is determined by the 
$\natural$-map. This is approximated with a high degree of accuracy by where the code 
$\Lambda_m$ sends  intervals. Accordingly:

\begin{definition}\label{uvc} Define the pair $(u,\rev{v})^c$  
to be $(u^c, sh^{A_m}(\rev{v^c})$. 
\end{definition}
Thus  $\la(u,\rev{v})^c\ra_l$ determines a basic open interval in $\bk^c\times \rev{\bl^c}$ which we might also write as $(\la u^c\ra_l\times \rev{\bl^c})\cap (\bk^c\times \la\rev{v^c}\ra_{l+A_m})$. Alternatively we could write this as:
\begin{eqnarray*}\{(f,g)\in \bk^c\times \rev{\bl^c}:f\rest[l,l+q_m)=u^c \mbox{ and }\\ g\rest[l+A_m, l+A_m+q_m)=\rev{v^c}\}.
\end{eqnarray*}

We now have a lemma extending Lemma \ref{first slip} which  says that if $u$ and $v$ belong to $\mcu_{n+1}$ and $\mcv_{n+1}$ then, relative to $sh^{-j_1}(\rev{v})$, all occurrences of $(u')^c\in \mcu_n^c$ in $u^c$ are either  lined up with an occurrence of a $\rev{(v')^c}$ for some $(v')^c\in \mcv_n^c$  or a boundary section of $sh^{-j_1}(\rev{v^c}))$. The lemma also says that  if $(u')^c, \rev{(v')^c}$ are lined up then $\vec{j}_{u'}$ and $\vec{j}_{v'}$ form a conjugate pair.\footnote{In this case both $\vec{j}_{u'}$ and $\vec{j}_{v'}$ are of length one.}

\begin{prop}\label{fair and equal}
Let $n<m$ and $u\in\mcu_m, v\in \mcv_m$. Then for $u'\in \mcu_n, v'\in\mcv_n$ we consider occurrences of $(u',\rev{v'})^c$ in $(u,\rev{v})^c$.\footnote{Since $A_m\ne A_n$ we are considering different shifts in $(u,\rev{v})^c$ and $(u',\rev{v'})^c$.} 
\begin{enumerate}
\item If $(u',\rev{v'})^c$ occurs in $(u,\rev{v})^c$, then $\vec{j}_{u'}$ and $\vec{j}_{v'}$ form a conjugate pair.
\item There is a constant $C=C(n,m)$ such that all conjugate pairs occur $C$ times.
\item Fix a conjugate pair $(j_{u'},j_{v'})$ of genetic markers  of $(u',\rev{v'})$. If $k$ is a location of an occurrence of $(u')^c$ in $u^c$ with genetic marker $\vec{j}_{u'}$, but not a location of 
$(u',\rev{v'})^c$, then the section of $sh^{A_m}(\rev{v^c})$ in the interval $[k+A_n, k+A_n+q_n)$ is contained in $\bigcup_{i=n+1}^m\boundary_i$.

\end{enumerate} 
\end{prop}
\pf Item 1 is immediate from the definitions.

 The latter items are asking about pairs of the form $((u')^c, sh^{A_n}(\rev{(v')^c}))$ occurring in 
$(u^c,sh^{A_m}(\rev{v^c})$. Such a pair occurs at $k$ if and only if the pair $((u')^c, \rev{(v')^c})$ occurs 
aligned in $(u,sh^{A_m-A_n}(\rev{v^c}))$ at $k$.   Item 3 is equivalent to 
saying that $(u')^c$ is lined up with a portion of $sh^{A_m-A_n}(\rev{v^c})$ contained in  
$\bigcup_{i=n+1}^m\boundary_i$.

We fix $m$ and prove 2 and 3 by induction on $m-n$. The case that $m=n+1$ is the content of Lemma \ref{first slip}. Suppose that the proposition is true for $m$ and $n+1$, we prove it for $m$ and $n$. 

A pair of $n+1$-circular words $(w_0,w_1)^c$ lined up in the shifted pair
 $(u^c, sh^{A_m-A_{n+1}}(\rev{v^c}))$  must have conjugate genetic markers. Moreover any there is a number $C_0$ such that any pair with conjugate genetic markers occurs lined up $C_0$ many times.

Fix  an occurrence $k$ of an $n+1$-word $w_0$ so that no word in $sh^{A_m}(\rev{v^c})$ occurs at $[k+A_{n+1}, k+A_{n+1}+q_{n+1})$, i.e  $w_0$ is not  lined up with the reverse of an $n+1$-word in $sh^{A_m-A_{n+1}}(\rev{v^c})$. Then $w_0$ is lined up with   a segment of $sh^{A_m-A_{n+1}}(\rev{v^c})$ that is a subset of in $\bigcup_{n+2}^m \boundary_i$. 
To pass from $A_m-A_{n+1}$ to $A_m-A_n$ we shift by $-j_1$, where $j_1=p_{n}^{-1} \mod{q_{n}}$. Noting that each reversed $n+1$-word ends with a string of $b$'s of length 
$q_n$, we see that after the additional shift there can be no $n$-subwords inside $w_0$ lined up with anything besides a portion of $sh^{A_m-A_n}(\rev{v^c})$  contained in
$\bigcup_{n+1}^m\boundary_i$.

Suppose that $u'$ and $v'$ are $n$-words and we have an occurrence of $(u')^c$ and $\rev{(v')^c}$  lined up in the pair
$(u^c, sh^{A_m-A_n}(\rev{v^c}))$. If $\vec{j}_{u'}=k_0\cat \vec{j^*_{u'}}$ and $\vec{j}_{v'}=k_1\cat\vec{j^*_{v'}}$, we let 
$(w_0,w_1)$ be the occurrence of  $n+1$-subwords of $(u,v)$ with genetic markers 
$\vec{j^*_{u'}}$ and $\vec{j^*_{v'}}$ that contain $u'$ and $v'$. It follows from the previous paragraph that the genetic markers of $w_0$ and $w_1$ are conjugate and $w_0^c, \rev{w_1^c}$ are aligned in $(u^c,sh^{A_m-A_{n+1}}(\rev{v^c}))$.  By Lemma \ref{first slip}, $k_0$ and $k_1$ are conjugate and thus $\vec{j}_{u'}$ and $\vec{j}_{v'}$ are conjugate. 

Further each conjugate pair occurs aligned the same number $C_1$ of times in the pair 
$(w^c_0, sh^{-j_1}(\rev{w^c_1}))$.  The number $C_1$ is independent of $w_0, w_1$ and 
$k_0$ and $k_1$.  
It follows now that given a conjugate pair of genetic markers $(\vec{j}_{u'}, \vec{j}_{v'})$, the 
number of occurrences of a pair of circular $n$-words with genetic marker 
$\vec{j}_{u'}$ in $u^c$ aligned with an occurrence of a circular word with genetic marker 
$\vec{j}_{v'}$ is in $v^c$ is $C_0*C_1$.

To finish we note that the unaligned $n$-words are in two categories, those that are not aligned because the $n+1$-words that contain them are not aligned, or those that are not aligned by the final shift $-j_1$. In each case, the unaligned $n$-words in $u$ occur across from boundary sections in the word $sh^{A_m-A_n}(\rev{v^c})$. \qed

Thus, using the backwards $\mcc$-operation to wrap words around the circle in opposite directions introduces some \emph{slippage}, but the slippage is uniform and predictable. 

\begin{definition}\label{slippery definition}
Suppose that $\vec{j}$ and $\vec{j'}$ are a conjugate pair of $(n,m)$-genetic markers and 
$u^c\in \mcu^c_m,v^c\in\mcv^c_m$. Let $(u')^c$ and $(v')^c$ have genetic markers $\vec{j}$ and $\vec{j'}$ in $u^c, v^c$ respectively.   Then the set of locations $k$ such that  
$(u')^c$  occurs in $u^c$ starting at $k$ with genetic marker $\vec{j}$ but $\rev{(v')^c}$ does not occur starting at $k+A_n$ in $sh^{A_m}(\rev{v^c})$ is called the $(n,m)$-slippage of $\vec{j}$. 
\end{definition}

A location $k$ can belong to the slippage of $\vec{j}$ for two mutually exclusive  reasons. Either, for some proper tail segment $\vec{j^*}$ of $\vec{j}$, $k$ is part of the slippage of the subword of $u^c$ with genetic marker $\vec{j^*}$ or $k$ is part of the slippage of the $j_{n}$ inside the $n+1$ word containing $u$ caused by $sh^{-j_1}$.

Let $SL_{n,m}$ stand for the $(n,m)$-slippage of $n$-subwords of $u^c$; i.e. the locations  $k$ in $u^c$ of some $n$-word $(u')^c$ such that there there is no 
$n$-word $\rev{(v')^c}$ at position $k+A_m$.  Inside an $m$-word $u^c$ we find multiple copies of $SL_{n,n+1}$ corresponding the location of each $n+1$ word in $u^c$. Denote the union of these copies as $SL^m_{n,n+1}$. Then it follows that:
\begin{equation}
SL_{n,m}=\bigcup_{k=n}^{m-1}SL^m_{k,k+1}\cap\{\mbox{locations of $n$-words}\} \label{slippage account}
\end{equation}
and moreover the union is disjoint.

The slippage is the portion of of the words that we have no control over when counting, so we want to be able to estimate the proportion of words in the slippage. 
Let
\begin{equation}
\varpi^m_n={|SL_{n,m}|\over |\{n-\mbox{subwords of }u^c\}}\label{all of slnm}
\end{equation}
The next proposition allows us to control the $(n,m)$-slippage by controlling the successive $(n,n+1)$-slippages.
\begin{prop}
 \begin{equation}\label{piling up}
1-\varpi^m_n=\prod_n^{m-1}(1-\varpi^{i+1}_i).
\end{equation}
  \end{prop}
 \pf We begin by noting that for $n^*$ between $n$ and $m$, all pairs $(u^*,\rev{v^*})$ of  $n^*$-words have the same proportion of slippage of $n$-words in $(u^*, \rev{v^*})^c$.  Thus
 $\varpi_{n}^{n^*}$ is equal to the  proportion of slippage of all of the $n$-words occuring in pairs $(u^*,\rev{v^*})^c$ of  $n^*$-subwords of  $(u,\rev{v})^c$.

 The argument is similar to Lemma \ref{dracula}. Starting with $n^*=m-2$ and decreasing until $n^*=n+1$, using that fact that the union in equation \ref{slippage account} is disjoint, one inductively demonstrates that:
 \[(1-\varpi^m_n)=(1-\varpi_{n}^{n^*})\prod_{n^*}^{m-1}(1-\varpi_{i}^{i+1}).\]
 \qed

We can combine item 3 of  Lemma \ref{fair and equal} with equation \ref{piling up} to see that if $k$ is in $SL_{n,m}$, then 
$[k+A_n, k+A_n+q_n)$ is a subset of $\bigcup_{i={n+1}}^m \boundary^{v^c}_i$.
It thus follows from Lemma \ref{88bis} that: 
\begin{equation}
1-\varpi^m_n\ge\prod_n^{m-1}(1-{2\over(l_i-1)}).\label{all that work for this?}
\end{equation}

Because the definition of $\varpi^m_n$ was made entirely in terms of genetic markers, the whole discussion could have been carried out simply by considering 
$\mck^c\times \rev{\mck^c}$. The numerics depend only on the circular coefficient sequence, not on particular construction sequences $\la \mcu_n, \mcv_n:n\in\nn\ra$. 

Viewing the operator $\natural$ as the limit of the codes $\Lambda_m$, we can pass to infinity and define
\hypertarget{slninf}{$SL_{n}^\infty$} similarly and let $\varpi^\infty_n$ be the proportion of locations $k$ of $n$-subwords of a typical $s\in \bk^c$ such that  no $n$-subword of $\natural(\rev{\pi(s)})$
occurs at $k+A_n$. 

Then:
\begin{eqnarray}
(1-\varpi^\infty_n)&=&\prod_n^\infty(1-\varpi^{i+i}_i)\notag \\
&\ge&\prod_n^\infty{(1-2/l_i)}\label{squeeze past}\\
&>&0.\notag
\end{eqnarray}
It follows that $\sum_1^\infty \varpi_i^{i+1}<\infty$.

\bigskip
We now formulate and prove the version of Lemma \ref{dracula} involving the $\natural$ map. One might expect that   would require considering arbitrary pairs of genetic markers $\vec{j}$ and $\vec{j'}$. However, by Proposition \ref{fair and equal}, if $u'$ occurs in $u$ with $(n,m)$-genetic marker $\vec{j}$, then the only genetic marker it can occur lined up with in $\rev{v}$ is its conjugate pair. Similarly  either of the genetic markers of aligned words $(u')^c$ occurring in $u^c$ and $sh^{A_n}(\rev{(v')^c})$ occurring in $sh^{A_m}(\rev{v})$ determine the other member of the conjugate pair.  

It follows that we need only consider pairs $(u',\rev{v'})$ whose genetic markers are conjugate in $(u,\rev{v})$. Since the map $\vec{j}$ to $\vec{j'}$ is a bijection we will refer to either of $\vec{j}$ or $\vec{j'}$ as the genetic marker of a pair $(u',\rev{v'})$ or equivalently $(u',\rev{v'})^c$. 

We are reduced to  considering sets $S^*\subseteq \{(n,m)\mbox{-genetic markers}\}$ rather than sets of pairs of genetic markers. 
Let $n<m$ and let $S^*$ be a set of $(n,m)$-genetic markers of pairs of $n$-words in $(u,\rev{v})$. Let 

	\begin{eqnarray}
	A&=&\{k\in [0,K_m): \mbox{some $u'$ with with genetic marker}\notag\\
	&& \mbox{in $S^*$ begins at $k$ in } u \} \notag
	\end{eqnarray}
and
	 \begin{eqnarray}	
	A^c&=&\{k\in[0,q_m):\mbox{for some $u'$ with genetic marker in $S^*$,}\notag\\
	&&\mbox{ there is a $v'$ such that $(u')^c$ occurs beginning at k in }u^c\notag\\
	&&\mbox{and $\rev{(v')^c}$ occurs beginning at }k+A_n \mbox{ in } sh^{A_m}(\rev{v^c})\} \notag
	\end{eqnarray}
and define
	\begin{eqnarray}
	d_m(A)&=&|A|/K_m\notag\\
	d^c_m(A^c)&=&|A^c|/q_m\notag
	\end{eqnarray}
If $(u')^c$ occurs at $k$ in $u$ and $\rev{(v')^c}$ occurs at $k+A_n$ in $sh^{A_m}(v)$ then $(u',\rev{v'})^c$ occurs at $k$ in $(u^c,sh^{A_m}(\rev{v^c})$

\begin{lemma}\label{prince of darkness}
Let $n< m$ and $(u,v)\in \mcu_m\times \mcv_m$. Let $S^*$ be a collection of $(n,m)$-genetic markers, $g$ the total number of $(n,m)$-genetic markers\footnote{As before it is easy to check that $g=\prod_n^{m-1}k_i$.} and $d=|S^*|/g$. Then (in the notation above):
\begin{eqnarray}
d_m(A)&=&{d\over K_n}\label{deja second mess}\\
d_m^c(A^c)&=&{d\over q_n}\prod_{p=n}^{m-1}(1-1/l_p)(\prod_{i=n}^{m-1}(1-\varpi^{i+1}_i)\label{deja first mess}\\
d_m(A)&=&\left({d_m^c(A^c)\over \prod_{p=n}^{m-1}(1-1/l_p)(\prod_{i=n}^{m-1}(1-\varpi^{i+1}_i)}\right)\left({q_n\over K_n}\right)\label{deja third mess}\\
d_m^c(A^c)&=&d_m(A)\left(\prod_{p=n}^{m-1}(1-1/l_p)\right)\left(\prod_{i=n}^{m-1}(1-\varpi^{i+1}_i)\right)\left({K_n\over q_n}\right).\label{deja fourth mess}
\end{eqnarray}

\end{lemma}
\pf The proof is essentially the same as the proof of Lemma \ref{dracula}, indeed the proof of equation \ref{deja second mess} \emph{is} the same.  Because all genetic markers occur with the same frequency, after allowing for the portions $u^c$ in boundary sections and in slippage (which are disjoint), $d/q_n$ is the density of locations $k$ of occurrences of words with genetic markers in $S^*$. Once again equations \ref{deja third mess} and \ref{deja fourth mess} follow from \ref{deja second mess} and \ref{deja first mess} by substitution.\qed

\noindent The equation relating $\rho\in \mathbf{P^-}$ and $\rho^c\in \mathbf{P^\natural}$ that  corresponds to  equation \ref{nu vs nuc} is:
\begin{equation*}
\rho^c(\la (u, \rev{v})^c\ra)=\left({K_n\over q_n}\right)\rho(\la (u,v)\ra)(1-\sum_n^\infty \rho^c(\boundary_m))(1-\varpi^\infty_n).
\end{equation*}
Once again $\rho^c(\boundary_m)$ is independent of the choice of $\rho^c$. Setting $d^{\boundary_m}_\rho=\rho^c(\boundary_m)$, we can write the previous equation as: 

\begin{equation}
\rho^c(\la (u, \rev{v})^c\ra)=\left({K_n\over q_n}\right)\rho(\la (u,v)\ra)(1-\sum_n^\infty d^{\boundary_m}_\rho)(1-\varpi^\infty_n).\label{rho vs rhoc}
\end{equation}
\bigskip

Understanding empirical distributions of joinings along the natural map involves studying how the slippage affects each pair of $n$-words. 
 Fix $u'\in \mcu_n, v'\in \mcv_n$ and $u\in \mcu_m, v\in \mcv_m$ where $n<m$. 
Let the conjugate pair $(\vec{j}, \vec{j'})$ be the genetic marker of $(u',\rev{v'})$ in $(u,\rev{v})$. Then, as remarked earlier $\vec{j'}$ is determined by $\vec{j}$, since they are a conjugate pair.
Define $SL_{n,m}(u', \rev{v'})$ to be the collection of locations $k\in SL_{n,m}$  of $n$-subwords of $u^c$ that have genetic marker $\vec{j}$. 
Item 2 of Proposition \ref{fair and equal} implies that $|SL_{n,m}(u',\rev{v'})|$ is the same for all choices of $(u',\rev{v'})$. Since $SL_{n,m}$ is the union over all possible pairs of $SL_{n,m}(u', \rev{v'})$, we see that 
\begin{eqnarray}
\varpi_n^m&=_{def}&{|SL_{n,m}|\over |\{n-\mbox{subwords of }u^c\}|}\notag\\
&=_{\ \ \ }&{|SL_{n,m}(u',\rev{v'})|\over |\{\mbox{subwords of }u^c\mbox{ with genetic marker } \vec{j}\}|}
\label{local slnm}
\end{eqnarray}
From the definition:
\[EmpDist_{n,n,A_n}((u,\rev{v})^c)((u')^c,(\rev{v'})^c)\]
is equal to
\[	{
	|\{\mbox{occurrences of }(u',\rev{v'})^c
	\mbox{ in }(u,\rev{v})^c\}|\over
	|\mbox{for some }(u^*,v^*)\in \mcw_n\times \mcv_n, (u^*,\rev{v^*})^c \mbox{ occurs in }(u,\rev{v})^c\}|}\]
 This in turn is equal to:
 \begin{equation}
 {(1-\varpi_n^m)|\{\mbox{subwords of }u^c\mbox{ with genetic marker} \vec{j}\}|\over
 (1-\varpi_n^m)|\{n\mbox{-subwords of }u^c\}|}
 \notag
 \end{equation}
 which in turn is equal to 
 \[EmpDist_{n,n,0}(u,\rev{v})(u',\rev{v'}).\]
 For notational convenience we write:
 \[EmpDist(u,\rev{v})(u',\rev{v'})=_{def}EmpDist_{n,n,0}(u,\rev{v})(u',\rev{v'})\]
 and 
 \[EmpDist((u,\rev{v})^c)((u',\rev{v'})^c)=_{def}\]
 \[EmpDist_{n,n,A_n}((u,\rev{v})^c)((u')^c,(\rev{v'})^c).\]
 
 \bfni{Summarizing:}
 \begin{equation}\label{just another miracle}
 EmpDist(u,\rev{v})(u',\rev{v'})=EmpDist((u,\rev{v})^c)(u',\rev{v'})^c
 \end{equation}

 \subsection{Transferring measures up and down, II}
 In this section we describe the correspondence between  joinings in $\mathbf{P^-}$ and $\mathbf{P^\natural}$. We do this by considering generic points for the joinings and transferring them up or down.
 
For the reader's convenience we repeat a definition.  Let $(s,t)$ be 
an arbitrary point in $\bk\times \bl$ with $\pi t=-\pi s$ and $s\in S^\bk, t\in S^\bl$.
Let $\la u_n:n\in\nn\ra$ and $\la v_n:n\in\nn\ra$ be 
the sequence of principal subwords of $s$ and $t$ 
respectively.  Then  
$\la u_n^c:n\in\nn\ra$ and $\la v_n^c:n\in\nn\ra$ are 
the sequences of principal subwords of $s^c=TU(s)$ and 
$t^c=TU(t)$.
If 
$x=\natural(\pi s^c)$, then $x\in \rev{\mck}$ and we can  set 
$r_n=r_n(x)$. Recall that we defined  
$\hat{t}\in \rev{\bl^c}$ by taking $\la \rev{v_n^c}:n\in \nn\ra$  as its principal $n$-subword
sequence and $\la r_n:n\in \nn\ra$ as its 
location sequence.

The following follows immediately from equation \ref{all is equal}:
\begin{lemma}\label{t and t-hat} The sequence $t$ is generic for an invariant measure $\mu$ on $\bl$ if and only if $\hat{t}$ is generic for an invariant measure $\mu^*$ on $\rev{\bl^c}$.
\end{lemma}
We will study the relationship between $\mathbf{P^-}$ and $\mathbf{P^\natural}$ via the function taking $(s,t)$ to $(s^c, \hat{t})$.
If $[a_n, b_n]$ is the location of the principal $n$-block of $s^c$, we define $w_n^c$ to be 
the word $(u_n^c, \hat{t}_n)$ (in the language $\Sigma\times \Lambda$)
where $\hat{t}_n=\hat{t}\rest[A_n+a_n,A_n+b_n]$.  Rephrasing this, if $(u_n, \rev{v_n})$ are 
the principal $n$-subwords of $(s,t)$ then $w_n^c=(u_n,\rev{v_n})^c$.

\begin{prop} \label{naturally ergodic}
The sequence $\la (u_n,\rev{v_n}):n\in\nn\ra$ is a generic sequence (resp. an ergodic sequence) if and only if $\la w_n^c:n\in\nn\ra$ is a generic sequence (resp. an ergodic sequence).
\end{prop} 
\pf 
This follows immediately from equation \ref{just another miracle}.
\qed
It is worth remarking that Proposition \ref{naturally ergodic} can be restated in the language of Definition \ref{shifted genericity} as saying that $\la (u_n, \rev{v_n},0):n\in\nn\ra$ is a generic sequence if and only if $\la (u_n^c, \rev{v_n^c}, A_n):n\in\nn\ra$ is a generic sequence.

\bigskip

  The next theorem is the analogue of Theorem \ref{k and kc} adapted to  lifting joinings of $\bk$ with 
  $\bl\inv$ to joining of $\bk^c$ with $(\bl^c)\inv$. 
  In the theorem  the notation   $(\nu, \nu^c)$ and $(\mu, \mu^c)$ refer to pairs of \hyperlink{corresponding measures}{corresponding measures}. 
 We assume that $\bk$ is built in the language $\Sigma$ 
 and $\bl$ is built in the 
 language $\Lambda$.

\begin{theorem}\label{k and l inv}
Suppose that $\la \mcu_n:n\in\nn\ra$ and $\la \mcv_n:n\in\nn\ra$ are construction sequences for two ergodic odometer based systems $(\bk,\nu)$ and 
$(\bl,\mu)$ with the same sequence parameters 
$\la k_n:n\in\nn\ra$.  Let $(\bk^c, \nu^c)$ and  $(\bl^c,\mu^c)$ be the 
associated ergodic circular systems built with a circular coefficient sequence $\la k_n, l_n:n\in\nn\ra$. Then there is a canonical affine homeomorphism $\rho\mapsto \rho^c$ between the simplex of anti-synchronous joinings $\rho$ of $(\bk,\nu)$ and $(\bl^{-1},\mu)$ and the simplex of anti-synchronous joinings of $(\bk^c,\nu^c)$ and $((\bl^c)^{-1},\mu^c)$
such that equation
\ref{rho vs rhoc} holds between $\rho$ and $\rho^c$. 
\end{theorem}

\pf Suppose that we are given an anti-synchronous ergodic joining $\rho$ 
between $\bk$ and $\bl\inv$.
Let $(s,t)$ be generic for $\rho$. By lemma \ref{generic seqs are ergodic}, the sequence of principal $n$-blocks, $\la (u_n, \rev{v_n}):n\in\nn\ra$ is ergodic. By Proposition \ref{naturally ergodic} the sequence $\la w_n^c:n\in\nn\ra$ define an ergodic measure $\rho^c$. Since the $\la(u_n,\rev{v_n}): n\in\nn\ra$ satisfy equation \ref{deja fourth mess}, the Ergodic Theorem
implies that $\rho^c$ and $\rho$ satisfy equation \ref{rho vs rhoc}. It is easy to check that the definition of $\rho^c$ is independent of the choice of the generic pair $(s,t)$.

For the other direction we can assume that we are given a generic pair $(s^c,\hat{t})$ for an ergodic measure $\rho^c$ on 
$\bk^c\times \rev{\bl^c}$ that concentrates on pairs 
$(s^c, \rev{t^c})\in \bk^c\times \rev{\bl^c}$ such that 
$\pi(\rev{t^c})=\natural(\pi(s^c))$.   Taking principal subwords  gives us a generic sequence 
$\la (u_n^c, \hat{t}_n):n\in\nn\ra$. Each $\hat{t}_n$ is  a well-defined word $\rev{v_n^c}$ in $\rev{\mcv_n^c}$.

 As in the 
\hyperlink{def of UT}{definition of $UT$} the pair $(s^c, \rev{\hat{t}})$ gives a pair of sequences of genetic markers $(\la j_n:n\ge N\ra,\la j_n':n\ge N\ra$ for some $N$. Letting $u_n=c_n\inv(u_n^c)$ and 
$v_n=c_n\inv(\rev{\hat{t}_n})$ the sequences $\la u_n, j_n\ra$ and $\la v_n, j'_n\ra$ determine a pair in $\bk\times \bl$ up to finite translations. These sequences are defined independently of the exactly location of the zero of $\hat{t}$; the small shifts used in the definition of $\natural$ do not change the two sequences.

 If we let $(s,t)=(UT(s^c), UT(\rev{\hat{t}}))$, making small adjustments if necessary to make $(s,t)$ anti-synchronous,  we get an element of $\bk\times\bl\inv$. Applying Proposition \ref{naturally ergodic} again we see the theorem.

 We can extend this correspondence to non-ergodic joinings $\rho$ on $\bk\times \bl^{-1}$ and $\rho^c$ on $\bk^c\times \rev{\bl^c}$, exactly as in Theorem \ref{k and kc}; to go up we take an ergodic decomposition of $\rho$:
 \[\rho=\int\rho_id\mu(i)\]
  and define 
  \[\rho^c=\int\rho_i^cd\mu(i).\]
  To go down we use the ergodic decomposition theorem and  the measure $\mu(i)$ to reverse this process.
 
 Clearly the map $\rho\mapsto \rho^c$ is an affine bijection. It remains to show that it is continous. However, just as in Theorem \ref{k and kc}, we see from equation \ref{rho vs rhoc}, that for each $n$ there is a constant $C_n$, independent of $\rho$ such that for all $u\in \mcu_n, v\in \mcv_n$, 
 \begin{equation*}
 \rho^c(\la (u,rev(v))^c\ra)=C_n\rho(\la (u, v)\ra).
\end{equation*}
This clearly implies that the map $\rho\mapsto\rho^c$ is a weak* homeomorphism.
\qed

The proof of Theorem \ref{k and l inv} shows that $(s,t)$ is generic for $\rho$ if and only if the pair $(s^c, \hat{t})$ is generic for $\rho^c$. Moreover, the proofs of  Theorems \ref{k and kc} and \ref{k and l inv} are quite robust. In particular the constructions of the corresponding measures are independent of the various choices of generic points $s$ or $s^c$, $(s,t)$ or $(s^c,\hat{t}\ )$.

\section{The Main Result}\label{the categories}

We now turn to the main results of this paper. Fix an arbitrary circular coefficient sequence  
$\la k_n, l_n:n\in \nn\ra$ for the rest of the section.  Let ${\mathcal OB}$ be the category whose 
objects are ergodic 
odometer based systems with coefficients $\la k_n:n\in \nn\ra$. The morphisms between objects 
$(\bk,\mu)$ and $(\bl,\nu)$ will be
 synchronous graph joinings of $(\bk,\mu)$ and $(\bl,\nu)$ or  anti-synchronous graph joinings of 
 $(\bk,\mu)$ and $(\bl^{-1},\nu)$. We call this the 
\emph{category of odometer based systems.}

Let $\mcc B$ be the category whose objects consists of all ergodic circular systems with coefficients
$\la k_n,l_n:n\in\nn\ra$.  The morphisms between objects $(\bk^c,\mu^c)$ and $(\bl^c,\nu^c)$ will be synchronous graph joinings of $(\bk^c,\mu^c)$ and $(\bl^c,\nu^c)$ or anti-synchronous graph joinings of $(\bk^c,\mu^c)$ and $((\bl^c)^{-1},\nu^c)$.
We call this the \emph{category of circular systems.}

\begin{remark}
Were we to be completely precise we would take  objects in $\mco B$ to be \emph{presentations} of odometer based systems by construction  sequences $\la \mcw_n:n\in\nn\ra$ without spacers together with suitable generic sequences and the objects in $\mcc B$ to be \emph{presentations} by circular construction sequences and their generic sequences. This subtlety does not cause problems in the applications so we ignore it.
\end{remark}

The main theorem of this paper is the following:
\begin{theorem}\label{grand finale}
For a fixed circular coefficient sequence $\la k_n, l_n: n\in\nn\ra$ the categories $\mco B$ and $\mcc B$ are isomorphic by a function $\mcf$ that takes synchronous joinings to synchronous joinings, anti-synchronous joinings to anti-synchronous joinings,  
isomorphisms to isomorphisms and {weakly mixing extensions to weakly mixing extensions.}
\end{theorem}

\medskip

Elaborating on Example \ref{for struct}:
\begin{corollary}\label{dis sys}
The map $\mcf$ preserves systems of  factor maps (or alternatively extensions). Explicitly: let $\la I, \le_I\ra$ be a partial ordering, $\la X_i:i\in I\ra$ be a family of odometer based systems  and $\la \pi_{i,j}:j\le i\ra$ is a commuting family of factor maps with $\pi_{i,j}:X_i\to X_j$. Then $\la \mcf(\pi_{i,j}):j\le i\ra$ is a commuting family of factor maps among $\la \mcf(X_i):i\in I\ra$. Moreover the analogous statement holds for circular systems $\la X_i^c:i\in I\ra$, factor maps $\la \pi_{i,j}:j\le i\ra$ and $\mcf^{-1}$.
\end{corollary}
Theorem \ref{grand finale} can be interpreted as saying that the \emph{whole} isomorphism and factor structure of 
systems based on the odometer $\la k_n:n\in\nn\ra$ is canonically  isomorphic  to the isomorphism and factor structure of 
circular systems based on $\la k_n, l_n:n\in\nn\ra$. We call this a \emph{Global Structure Theorem}.

\subsection{The proof of the main theorem}
Before we prove theorem \ref{grand finale} we owe the following lemma:

\begin{lemma}\label{paying debt}
Both $\mco B$ and $\mcc B$ are categories, and the composition of synchronous joinings is synchronous, the composition of two anti-synchronous joinings is synchronous and the composition of a synchronous and an anti-synchronous joining (in either order) is anti-synchronous.
\end{lemma}
\pf To see that $\mco B$ and $\mcc B$ are categories we must see that the morphisms are closed under composition. This is equivalent to the statement that the composition of two synchronous or anti-synchronous joinings are synchronous or anti-synchronous. This, in turn follows from Proposition \ref{laziness} (item 2) applied to  joinings of odometers or rotations. \qed

We now prove Theorem \ref{grand finale}.
\pf By Proposition \ref{bijection from readability} the map $\mathcal F$ gives a bijection between the objects of $\mco B$ and $\mcc B$ and hence it remains to define the functor on the morphisms (i.e. joinings between systems $(\bk, \mu)$ and $(\bl^{\pm 1}, \nu)$) and show that it preserves composition. 

\subsubsection{Defining $\mcf$ on morphisms} \label{f on morphs}We split the definition of $\mcf(\rho)$ into two cases according to whether $\rho$ is synchronous or anti-synchronous. In both cases we define $\mcf$ for arbitrary joinings even though the only joinings we use as morphisms in the categories are graph joinings; in particular the morphisms in each category are ergodic. 
\medskip

\bfni{Case 1: $\rho$ is synchronous:}
\smallskip

Suppose  that $\rho$ a  synchronous joining of odometer based systems $\bk$ and $\bl$ with coefficient sequence $\la k_n:n\in\nn\ra$ that are constructed with symbols in $\Sigma$ and $\Lambda$ from construction sequences $\la \mcu_n:n\in\nn\ra$ and $\la\mcv_n:n\in\nn\ra$. We define a new construction sequence $\la \mcw_n:n\in\nn\ra$ with the symbol set $\Sigma\times \Lambda$.

Given $n$, we put a sequence 
\[\la (\sigma_0,\lambda_0), (\sigma_1,\lambda_1)\dots (\sigma_{K_{n}-1},\lambda_{K_{n}-1})\ra\]
into $\mcw_n$ if and only there are words $u=(\sigma_0,\dots \sigma_{K_{n}-1})\in \mcu_n$ and $v=(\lambda_0,\dots \lambda_{K_{n}-1})\in \mcv_n$.

It is easy to check that $\la \mcw_n:n\in\nn\ra$ is an odometer based construction sequence 
with coefficients $\la k_n:n\in\nn\ra$. Let $(\bk,\bl)^\times$ be the associated odometer based system. Since  $\rho$ is synchronous, it concentrates on members of 
$\bk\times \bl$ that correspond to elements of $(\bk,\bl)^\times$. We can canonically identify $\rho$ with a shift invariant measure $\nu$ on $(\bk,\bl)^\times$.

Let $((\bk,\bl)^\times)^c$ be the circular system associated with $(\bk,\bl)^\times$.  We can  apply Theorem \ref{k and kc} to 
find  shift invariant measure $\nu^c$ on $((\bk,\bl)^\times)^c$ associated with $\nu$ that is ergodic just in case $\nu$ is ergodic. Shift invariant 
measures on $((\bk,\bl)^\times)^c$ can be canonically identified with synchronous joinings on $\bk^c\times\bl^c$. Let $\rho^c$ be 
the joining of $\bk^c\times\bl^c$ corresponding to $\nu^c$.  We let $\mcf(\rho)=\rho^c$.%

Explicitly: A generic sequence $\la (u_n, v_n,0):n\in\nn\ra$  for the joining $\rho$,  can be viewed as a generic sequence $\la (u_n,v_n):n\in\nn\ra$ for $(\bk,\bl)^\times$ and transformed  into a generic sequence $\la (u_n^c, v_n^c):n\in\nn\ra$ for $((\bk,\bl)^\times)^c$. The latter corresponds to a generic sequence of the form $\la (u^c_n,v^c_n,0):n\in\nn\ra$ for the {joining} $\rho^c$. 
This process is clearly reversible so $\mcf$ is a bijection between the synchronous joinings of $\mco B$ and the synchronous joinings of $\mcc B$.

We must show that if $\rho$ is a graph joining then so is $\rho^c$. Once this is established it follows by symmetry that if $\rho$ is an isomorphism then $\rho^c$ is an isomorphism.  Namely if  $\rho^*$ is the adjoint joining of $\bl$ with  $\bk$ defined as $\rho^*(A)=\rho(\{(s,t):(t,s)\in A\})$, then $(\rho^*)^c=(\rho^c)^*$.  Hence 
$\rho^*$ is a graph joining iff $(\rho^c)^*$ is a graph joining.

Suppose that $\rho$ is a graph joining. We apply Proposition \ref{no proof}, part \ref{working graph}. It suffices to show that for all basic open sets in $\bk^c$ of the form $\la u^c\ra_0$ where $u^c\in \mcu_n^c$ and all $\epsilon>0$, there are words $ v_1^c,  v_2^c \dots  v^c_{k^*}$ that belong to $\bigcup_n\mcv_n^c$ and locations $l^c_1, \dots l^c_{k^*}$ such that: 
	\begin{equation}\label{upchuck}\rho^c((\la u^c\ra_0\times {\bl^c})\Delta (\bk^c\times \bigcup \la {v^c_j}\ra_{l^c_j}))<	
	\epsilon.
	\end{equation}

{Consider $u$ such that $c_n(u)=u^c$. Because $\rho$ is a graph joining, for all $\delta>0$ we can find words $v_1,\dots v_{k'}$ and locations $l_1, \dots l_{k'}$ such that 
	\begin{equation}\label{upchuck2}\rho((\la u\ra_0\times {\bl})\Delta (\bk\times 		\bigcup_{i\le k'} \la {v_i}\ra_{l_i}))<
	\delta.
	\end{equation}
Without loss of generality we can assume that  for some $m\ge n$ each $v_i$ is an $m$-word and that each $l_i\le 0$.}
	
Let $(s,t)$ be generic for $\rho$ and considering the pair $s^c=TU(s)$, $t^c=TU(t)$. Then by  Remark \ref{benjy's 
 remark}
$(s^c,t^c)$ is generic for $\rho^c$. We will choose words $v_j^c$ and  locations $l_j^c$ and compute the measure in inequality \ref{upchuck} by computing the density of locations representing points in the symmetric difference.

\noindent Let
	 \begin{eqnarray}
	 B_0&=&\{k:\mbox{$u$ occurs at $k$ in $s$, but for no $i$ does $v_i$ occur in $t$}\notag \\
	 && \mbox{at $l_i+k$}
	 \}\label{downstairs s}\notag\\
	 B_1&=&\{k:\mbox{for some $i$, $v_i$ occurs in $t$ at $k+l_i$  but $u$ does}\label{downstairs t}\notag\\
	 &&\mbox{not occur in $s$ at $k$}\}\notag
	  \end{eqnarray}
 By inequality \ref{upchuck2}, $B_0\cup B_1$ can be taken to have density less than 
 $\delta$.  
\medskip

Given words and locations $\{v_j^c, l^c_j:j\in J\}$ we can define two sets $B^c_0, B^c_1\subseteq \poZ$, as follows:
	\begin{eqnarray}
	B^c_0&=&\{k: u^c\mbox{ occurs in $s^c$ at $k$ but for no $j$ does $v_j^c$ occurs in $t^c$}\notag\\
	&&\mbox{at $l^c_j+k$}\}\notag\\
	B^c_1&=&\{k:\mbox{for some $j$\ \ $v^c_j$ occurs in $t^c$ at $l^c_j+k$ but $u^c$ does}\\
	&& \mbox{\ \ \ not occur in $s^c$ at $k$.}\} \notag
	\end{eqnarray}
We need to find the words and locations $v_j^c, l_j^c$ so that the density of $B^c_0 \cup B^c_1$ is less than $\epsilon$. 

For each $i$, if $-l_i$ is not the location of the beginning of an $n$-word in $v_i$ then dropping $\la v_i\ra_{l_i}$ reduces the measure of the  symmetric difference in inequality \ref{upchuck2}. Thus, without loss of generality we can assume that for all $i$, there is an $(n,m)$-genetic marker $\vec{j}(i)$ coding the location of the $n$-word in $v_i$ that starts at $-l_i$. Since $B_0\cup B_1$ has density less than $\delta$,  the density of $k$ such that either:
    \begin{enumerate}
    \item $u$ occurs at $k$ but for each $i$, $k$ is not the position of the beginning of an $n$-word with genetic marker $\vec{j}(i)$ in an occurrence of $v_i$ or 
    \item for some $i$,  $k$ \emph{is} the position of the beginning of an $n$-word with genetic marker $\vec{j}(i)$ in an occurrence of $v_i$, but $u$ does not occur at $k$,
    \end{enumerate}
 has density less than $\delta$.

We are in a position to define the $v_j^c$ and the $l_j^c$. For each $i$ we define index sets $J_i$ and a collection $\{l_j^c:j\in J_i\}$.  We arrange the $J_i$'s so that they are pairwise disjoint and for some $k^*$, $\bigcup_i J_i=\{j:1\le j\le k^*\}$. For 
$j\in J_i$, all of the $v_j^c$ are the same and equal to  $c_m(v_i)$.  For a fixed $i$, let $\{-l^c_j:j\in J_i\}$ be the collection of locations of  the beginnings of 
$n$-subwords of $c_m(v_i)$ that have genetic marker $\vec{j}(i)$. 

To compute the density of $B_0^c\cup B_1^c$,  it suffices to consider an extremely large $M$ and compute the density of $B_0^c\cup B_1^c$ inside the principal 
$M$-subword $(w^c_0,w^c_1)$ of $(s^c, t^c)$.  Let $(w_0,w_1)$ be the principal $M$-subword of $(s,t)$ and $c_M(w_0)=w_0^c$ and $c_M(w_1)=w_1^c$.

We now argue as in Lemma \ref{dracula}.  Let $d_0$ be the density of $B_0\cup B_1$ in 
$(w_0, w_1)$
 and $d^c_0$ be the density of $B^c_0\cup B^c_1$ in $(w^c_0, w^c_1)$. Among all $n$-words the proportion $d_p$  that begin with an element of $B_0\cup B_1$ is $d_0*K_n$. 
The density of $k\in \poZ$ that start $n$-words in $(s,t)$ is $(1-\mu(\bigcup_n^\infty\boundary_i))/q_n$. Letting $d^*$ be the density of $k\notin \bigcup_n^M\boundary_i$, we see that $d^*$ is bounded away from $0$ and $1$ independently of $M$.
The proportion $d_p^c$ of circular 
$n$-subwords of $(w^c_0, w^c_1)$ that begin with a $k\in B^c_0\cup B^c_1$ is 
\[{d_0^c*q_n\over (1-d^*)}.\]
 Since $\rho$ concentrates on $\{(s,t):\pi(s)=\pi(t)\}$ and $\rho^c$ concentrates on $\{(s^c,t^c):\pi(s^c)=\pi(t^c)\}$, the $n$-words with a particular genetic marker in $w_0$ occupy the position of the same genetic marker in $w_1$ and similarly for  $w_0^c$ and $w_1^c$. The $(n,M)$-genetic markers set up a one-to-one correspondence between $n$ subwords $u^*$ of $w_0$ and regions of $w_0^c$ that consist of occurrences of $(u^*)^c$ that have the same genetic marker. Each of the regions of $w_0^c$ with the same genetic marker have the same number of $n$-words in them.

Temporarily call an $n$-subword of $(w_0^c, w_1^c)$ \emph{bad} if it begins with a $k$ in $B_0^c\cup B_1^c$ and similarly for $n$-subwords of $(w_0,w_1)$ and $B_0\cup B_1$.  Then the property of being bad is determined by the $(n,M)$-genetic marker of the $n$-word: if $k$ is the beginning of $n$-subword of $w_0$ with genetic marker $\vec{j}$, and $k'$ is the beginning of an $n$-subword of $w_0^c$ with the same genetic marker in $w_0^c$, then $k\in B_0\cup B_1$ if and only iff $k'\in B_0^c\cup B_1^c$.

It follows the proportion of bad $n$-subwords of $(w_0,w_1)$ is the same as the proportion of bad subwords of  $(w_0^c, w_1^c)$. In otherwords:
\[d_p=d_p^c.\]
It follows that
\[{d_0* K_n}={d_0^c*q_n\over (1-d^*)}.\]
Thus by taking $\delta$ small enough and $M$ large enough we can make $d_0$ as small as we want, and thus arrange that $d_0^c\ll\epsilon$ as desired.
\smallskip

To finish showing  that $\mcf$ is a bijection between graph joinings in each category and isomorphisms in each category we must also show that if $\rho^c$ is a graph joining then so is $\rho$.  But this is very similar. Given a $u^c\in \mcu^c_n$, and an $\epsilon>0$ we can find $v_1^c, \dots v_{k^*}^c$ and locations $l^c_1,\dots l^c_{k^*}$ so that inequality \ref{upchuck} holds. Again we can assume that for some $m$, for all $j$, $v_j^c\in \mcw_m^c$. The numbers $|l_j^c|$ determine 
locations in $v_j^c$ of beginnings of $n$-words. We can augment our collection of locations by adding more $l_j^c$'s so that if
$l$ is the start of a location in $v_j^c$ that has the same $(n,m)$-genetic marker as $l_j^c$, then for some $j'$ we have $l_{j'}^c=-l$ and $v_{j'}^c=v_j^c$. In doing this we do not increase the density of $B_0^c\cup B_1^c$. Reversing the procedure above this gives words $v_j\in \bigcup_n\mcv_n$ and locations $l_j$ such that the density of $B_0\cup B_1$ is less than $\epsilon$. (Note the lack of boundary in $\bk\times \bl$ makes the computation easier by reducing the density of $B_0\cup B_1$.)

\medskip

\bfni{Case 2: $\rho$ is anti-synchronous}

On the anti-synchronous joinings we take $\mcf$ to be the bijection between anti-synchronous joinings
 of $(\bk,\mu)$ with $(\bl^{-1},\nu)$ and of the circular systems $(\bk^c,\mu^c)$ with $((\bl^c)^{-1},\nu^c)$ defined in Theorem \ref{k and l inv}.
We show that $\mcf$ takes  anti-synchronous graph joinings to anti-synchronous graph joinings and vice versa. Having done this it will follow by a symmetry argument that $\mcf$ sends anti-synchronous isomorphisms to anti-synchronous isomorphisms.

Suppose that $\rho$ is an anti-synchronous graph joining; i.e.  $\rho$ is a graph joining of $\bk$ with $\bl\inv$ that concentrates on $\{(s,t):\pi(t)=-\pi(s)\}$. The map $x\mapsto \rev{x}$ projects to the odometer map $\pi(x)\mapsto -\pi(x)$; in particular $\rev{\bl}$ is based on the same odometer that 
$\bl$ is. By Lemma \ref{joining correspondence} we can view $\rho$ as a graph joining of $\bk$ with 
$\rev{\bl}$ that concentrates on $\{(s,t):\pi(s)=\pi(t)\}$.  Similarly we view $\rho^c$ as concentrating on $\bk^c\times \rev{\bl^c}$.

We must show that for all basic open sets in $\bk^c$ of the form $\la u^c\ra_0$ where $u^c\in \mcu_n^c$ and all $\epsilon>0$, there are words $ v_1^c,  v_2^c \dots  v^c_{k^*}$ that belong to $\bigcup_n\mcv_n^c$ and locations $l^c_1, \dots l^c_{k^*}$ such that: 
	\begin{equation}
	\rho^c((\la u^c\ra_0\times \rev{\bl^c})\Delta (\bk^c\times \bigcup \la \rev{v^c_j}\ra_{l^c_j}))<	
	\epsilon.\notag
	\end{equation}

{Consider $u$ such that $c_n(u)=u^c$. Because $\rho$ is a graph joining for all $\delta>0$  and all large enough $m$ we can find words $v_1,\dots v_{k'}\in\mcv_m$ and locations $l_1, \dots l_{k'}$ such that 
	\begin{equation}\label{upchuck7}\rho((\la u\ra_0\times \rev{\bl})\Delta (\bk\times \bigcup \la\rev{v_i}\ra_{l_i}))<
	\delta.
	\end{equation}
Without loss of generality we can assume that each $l_i\le 0$.} We will take $m$ sufficiently large according to a restriction we define later.

Let $(s,t)$ be generic for $\rho$ and let $\hat{t}$ be as in Definition \ref{t-hat}.  Then $(s^c,\hat{t})$ is generic for $\rho^c$. We argue as before considering sets:

	\begin{eqnarray}
	 B_0&=&\{k:\mbox{$u$ occurs at $k$ in $s$, but for no $i$ does $\rev{v_i}$ occur in $\rev{t}$}\notag\\
	 && \mbox{at $l_i+k$.}\label{downstairs s natural}
	 \}\\
	 B_1&=&\{k:\mbox{for some $i$, $\rev{v_i}$ occurs in $\rev{t}$ at $k+l_i$  but $u$ does}\label{downstairs t natural}\\
	 &&\mbox{not occur in $s$ at $k$.}\}\notag
	  \end{eqnarray}
Then inequality \ref{upchuck7}, shows that $B_0\cup B_1$ can be taken to have density less than any positive $\delta$.

Given words and locations $\{v_j^c, l^c_j:j\in J\}$ we consider $B^c_0, B^c_1\subseteq \poZ$, as follows:
	\begin{eqnarray}
	B^c_0&=&\{k: u^c\mbox{ occurs in $s^c$ at $k$ but for no $j$ does $v_j^c$ occurs in $\hat{t}$}\notag\\
	&&\mbox{at $l^c_j+k$}\}\notag\\
	B^c_1&=&\{k:\mbox{for some $j$, $v^c_j$ occurs in $\hat{t}$ at $l^c_j+k$ but $u^c$ does}\notag\\
	&& \mbox{\ \ \ not occur in $s^c$ at $k$.}\} \notag
	\end{eqnarray}
Given $\{(v_i,l_i):1\le i\le k'\}$, we need to find the words and locations $v_j^c, l_j^c$ so that the density 
of $B^c_0 \cup B^c_1$ is less than $\epsilon$. As in the synchronous case, for each $i$ we build index 
sets $J_i$  so that the $J_i$'s to be disjoint and have union the interval $\{j:1\le j\le k^*\}$ for some 
$k^*$. For all $j\in J_i$ we take $v^c_j=c_m(v_i)$. We need to find a collection of locations 
$\{l_j:j\in J_i\}$.

 Fix an $i\le k'$. Without loss of generality we can assume that $l_i$ is the beginning of a reversed $n$-block $\rev{v'}$ in $\rev{v_i}$, since otherwise, discarding $\la \rev{v_i}\ra_{l_i}$ makes inequality \ref{upchuck7} sharper. 
If $(s_0,\rev{t_0})\in \bk\times \rev{\bl}$ is an arbitrary member of 
\[(\la u\ra_0\times {\bl})\cap(\bk\times \la \rev{v_i}\ra_{l_i})\]
with $\pi(s_0)=-\pi(t_0)$, then
there is an $m$-word $u^*$ such that $s_0\in \la u^*\ra_{l_i}$. Let $\vec{j}(i)$ be the genetic marker of $u$ in $u^*$.
We note that $\vec{j}(i)$ does not depend on $s_0$, since it is  determined entirely by the location of $u$ in $u^*$ and $u^*$ must be aligned with $\rev{v_i}$.

The genetic marker $\vec{j}(i)$ defines a region of $n$-words in $\mcu_n^c$ inside an $m$-word in $\mcu_m^c$. Let $L_i$ be the collection of $l$ that are at the beginning of an $n$-word in $\mcu_n^c$ with genetic marker $\vec{j}(i)$ in an $m$-word  in $\mcu_m^c$ and set
\begin{equation}
\{l^c_j:j\in J_i\}=\{A_m-l: l \in L_i\}.
\end{equation}
This determines the collection $\{v^c_j,l_j^c:1\le j\le k^*\}$.

We now compute the density of $B_0^c\cup B_1^c$ in terms of the density of $B_0\cup B_1$.
To do this  it suffices to consider a large enough $M$ that $s^c$ has a principal $M$-block $[a_M,b_M)$ and compute densities inside this principal $M$-block. If this is sufficiently small we can deduce that the density of $B_0^c\cup B_1^c$ is small in $\poZ$. By Remark 
\ref{natural coding}, we can also assume that $M$ is so large that $\natural$ restricted to this principal $M$-block is equal to 
$\bar{\Lambda}_M$ along this $M$-block; equivalently the principal $M$-block of $\hat{t}$ is $[a_M+A_m, b_M+A_M)$.

From Proposition \ref{fair and equal}, we know that if $I$ is an $m$-sub-block of $s^c\rest [a_M, b_M)$ then either:
\begin{enumerate}
\item the corresponding sub-block of $\hat{t}$ is at $sh^{A_m}(I)$ or
\item $I$ is part of the $(m,M)$-slippage.
\end{enumerate}
By item 2 of Proposition \ref{fair and equal}, the number of $m$-sublocks in each case that correspond to a given $(n, M)$-genetic marker does not depend on the genetic marker.  Further in the second case $sh^{A_m}(I)$ is entirely part of 
$\bigcup_{m+1}^M\boundary_i(\hat{t})$.

We compute the density $d_0^c$ of elements of $B^c_0\cup B^c_1$ by separating them into these two  sources. Explicity, we divide into:
\begin{description}
\item[Slippage:] Those $k\in B^c_0\cup B^c_1$
that begin an $n$-subword of a  location of an $m$-subword of $s^c$ that is in the $(m,M)$-slippage. 
\item[Mistakes:] those $k\in B^c_0\cup B^c_1$ such that $k$ is the location of the beginning of a circular $n$-subword inside $s^c\rest [a_M, b_M)$ and $[k+A_m,k+q_m+A_m)$ is the location of an $m$-word in $\hat{t}$.

\end{description}
We compute the density of the Mistakes and the Slippage separately.  Again we will call $n$-subwords that begin with elements of $B_0\cup B_1$ or $B_0^c\cup B_1^c$ \emph{bad}.

Both the Mistakes and the Slippage occur at the beginning of $n$-subwords of $s^c\rest[a_M, b_M)$. Define $d_b$ to be density of 
$\bigcup_{n+1}^M\boundary_i$ in $[a_M,b_M)$.
Then proportion of $k\in [a_M, b_M)$ that begin $n$-subwords is:
\[{1-d_b\over q_n}.\]
Of these a proportion $\varpi_m^M$ of the $n$-subwords are in the Slippage. Thus the collection of $k$ that belong to the Slippage has density

\[\varpi_m^M\left({1-d_b\over q_n}\right).\]
Since $\varpi_m^M$ goes to zero as $m$ goes to infinity we can make this term as small as desired by taking $m$ large enough.

Let $[a'_M, b'_M)$ be the location of the principal $M$-block of $s$ (and thus of $\rev{t}$).  Let $d_0$ be the density of $B_0\cup B_1$ in $[a'_M, b'_M)$.

Suppose now that $k$ belongs to the Mistakes. Let  $\vec{j}$ be the $(n,M)$-genetic marker of the word beginning with  $k$ in $s^c\rest[a_M, b_m)$. Then there is a unique $k'$ in $[a'_M,b'_M)$ that is at the beginning of an $n$-subword of $s\rest[a'_M, b'_M)$ and has genetic marker $\vec{j}$. By construction, for $k$ that are not in the Slippage:
\begin{equation}k\in B_0^c\cup B_1^c \mbox{ iff } k'\in B_0\cup B_1.\label{relative proportions equiv}
\end{equation}
Let $d_p$ be the proportion of $m$-subwords of $s\rest[a'_M, b'_M)$ that begin with a $k\in B_0\cup B_1$. Since every genetic marker is represented exactly the same number of times in the complement of the slippage (Proposition \ref{fair and equal}), the proportion of 
words that begin with $k$ in the Mistakes  is 
\begin{equation}d^c_p=d_p*(1-\varpi_m^M).\label{relative proportions}
\end{equation}
If $d_0$ is the density of $B_0\cup B_1$ in $[a'_M, b'_M)$ and $d^c_0$ is the density of the Mistakes, then  
\begin{eqnarray}
d_0&=&d_p/K_n \label{downstairs density} \\
d_0^c&=& d^c_p\left({1-d_b\over q_n}\right)\label{upstairs density}
\end{eqnarray}
Putting together equations \ref{relative proportions}, \ref{downstairs density} and \ref{upstairs density}, we see that if we make $d_0$ sufficiently small we can make $d_0^c$ as small as desired.
\bigskip

\bfni{Summarizing:} By taking $M$ large enough, the density of $B_0^c\cup B_1^c$ is well approximated by the density of $B_0^c\cup B_1^c$ inside $[a_m, b_m)$. This is the sum of the density of the $(m,M)$ slippage and the density of the Mistakes. We can make the density of the Slippage arbitrarily small by taking $m$ large enough and the density of the Mistakes arbitrarily small by taking $\delta_0$ sufficiently small. This establishes the claim that if $\rho$ is a graph joining then so is $\rho^c$.

\bigskip

We must show that if $\rho^c$ is a graph joining then so is $\rho$. We suppose that we are given a $u\in \mcu_n$, we must find $\{v_i,l_i:i\le k'\}$ so that equation \ref{upchuck7} holds. Let $u^c=c_n(u)$ and approximate $\la u^c\ra_0\times \rev{\bl^c}$ using $\{v_j^c, l_j^c:i\le k^*\}$. Again, we can assume that the collection of locations is saturated in the sense that if
$l$ is the start of a location in $v_j^c$ that has the same $(n,m)$-genetic marker as $l_j^c$, then for some $j'$ we have $l_{j'}^c=-l$ and $v_{j'}^c=v_j^c$. In doing this we do not increase the density of $B_0^c\cup B_1^c$. We can now use equations \ref{relative proportions}, \ref{downstairs density} and \ref{upstairs density} again to see that if $d_0^c$ is made sufficiently small then so is $d_0$.

\bigskip

Our next claim is that $\rho$ is an isomorphism if and only if $\rho^c$ is an isomorphism. Recall from Proposition \ref{invertible symmetry} that $\rho$ is an isomorphism iff both $\rho$ and $\rho^*$ are graph joinings. Thus if $\rho$ is an isomorphism, both $\rho^c$ and $(\rho^*)^c$ are graph joinings. Since $\natural$ is an involution:
\[(\rho^*)^c=(\rho^c)^*.\]
Thus if $\rho$ is an isomorphism, so is $\rho^c$. 

Reversing this line of reasoning shows that if $\rho^c$ is a graph joining then $\rho$ is.

\subsubsection{$\mcf$ preserves composition}\label{composition}
To finish the proof that $\mcf$ is a functor we must show that $\mcf$ preserves composition. The 
argument splits into four natural cases: composing synchronous joinings, composing a synchronous 
joining with an anti-synchronous joining on either side and composing two anti-synchronous joinings.  
We will carefully work out the case for compositions of synchronous embeddings, and discuss the 
appropriate modification in the cases involving at least one anti-synchronous embedding 
after Lemma \ref{compositions}.

The cases differ only that the shifts involved in the generic sequences have different forms. For ergodic synchronous joinings  generic sequences can be taken to be of the form $\la (u_n, v_n, 0):n\in\nn\ra$, whereas for anti-synchronous joinings of $\bk^c$ and  $\rev{\bl^c}$ a natural generic sequence is of the form $\la (u^c_n, \rev{v^c_n}, A_n):n\in\nn\ra$.\footnote{i.e. $\la (u_n,\rev{v_n})^c:n\in\nn\ra$.}

\bfni{Preparatory Remarks}

In  the characterization of the  relatively independent joining $\rho$ of $\rho_1$ and $\rho_2$ given in Lemma \ref{allow us to disintegrate} and Proposition \ref{rel ind join}, the partitions $\mca_k, \mca'_k$ and $\tilde{\mca}_k$ are given by $\la u_k\ra_{s_1}, \la v_k\ra_{s_2}$ and $\la w_k\ra_{s_3}$ for $s_1, s_2, s_3\in \poZ$. Formally the partitions 
$\mca_k\times\mca'_k, \mca_k\times\tilde{\mca}_k$ and $\mca'_k\times\tilde{\mca}_k$ and $\mca_k\times\mca'_k\times \tilde{\mca}_k$
consist of all possible products of these basic open sets. However, in the situation we are considering we have synchronous and anti-synchronous joinings.  For synchronous joinings we can build a generating family for the relatively independent joining $\rho$ of $\rho_1$ and $\rho_2$ by considering products of pairs of basic open intervals in the same locations; e.g. pairs of the form $\la u_k\ra_s\times \la w_k\ra_s$.  As a consequence, for verifying the hypotheses of Proposition \ref{rel ind join} we can restrict our attention to the case where $s^*=0$.

In the case of  anti-synchronous joinings we need to distinguish the odometer based from the circular systems. For anti-synchronous joinings of odometer based systems $\bk$ with $\bm^{-1}$ we can consider only intervals  of the form $\la u_k\ra_s\times \la \rev{w_k}\ra_{s+s^*}$ where $s^*=0$. For anti-synchronous joinings of the circular systems $\bk^c$ with $\bm^c$, asymptotically the Empirical Distances concentrate on words of the form $\la u_k^c\ra\times \la\rev{w_k^c}\ra_{A_k}$ (where $A_k$ is the amount of shift for $\natural$ at scale $k$). Moreover, translations of sets of this form generate the measure algebra of the anti-synchronous joining. 

Thus in the proof of the next lemma, to verify the hypothesis 3 of Proposition \ref{rel ind join} we can take $s^*=0$ or $s^*=A_k$ depending on whether $\rho_1\circ\rho_2$ is synchronous or anti-synchronous.

\bigskip

Fix odometer based systems $\bk$, $\bl$ and $\bm$ with construction sequences $\la \mcu_n:n\in\nn\ra$, $\la \mcv_n:n\in\nn\ra$ and $\la \mcw_n:n\in\nn\ra$ respectively. Let  $\rho_1$ and $\rho_2$ be synchronous graph joinings of $\bk$ and $\bl$, and $\bl$ and $\mathbb M$ 
respectively and $\rho$ their relatively independent joining over $\bl$.

Since $\rho_1$ and $\rho_2$ are graph joinings so is their composition. Thus the relatively independent joining is ergodic. Hence by Lemma \ref{ill existe}
 we can find generic sequences for $\rho_1, \rho_2$ and $\rho$ that satisfy the hypothesis of Proposition \ref{rel ind join}.
\begin{lemma}\label{compositions}Let $\la (u_n, v_n, w_n,0, 0):n\in\nn\ra$ be generic for $\rho$.
Then the sequence $\la (u_n^c, v_n^c, w_n^c,0,0): n\in\nn\ra$ is generic for the relatively independent joining $\rho^c$ of $\rho_1^c$ with $\rho_2^c$.
\end{lemma}

Assuming the lemma, we show that $\mcf$ preserves compositions.  Corollary \ref{en fin} shows that 
$\la (u_n, w_n,0):n\in\nn\ra$ is generic for $\rho_1\circ \rho_2$. From the way that $\mcf$ is 
constructed, if $\nu^c=\mcf(\rho_1\circ\rho_2)$, then 
$\la (u_n^c, w_n^c):n\in\nn\ra$ is generic for $\nu^c$ (viewed as a measure on a circular system). From Lemma \ref{compositions} and Corollary \ref{en fin}, we know that 
$\la (u_n^c, w_n^c,0):n\in\nn\ra$ is generic for $\rho_1^c\circ\rho_2^c$. Hence $\mcf(\rho_1\circ\rho_2)=\mcf(\rho_1)\circ\mcf(\rho_2)$ as desired. 
\bigskip

\noindent It  remains to prove Lemma \ref{compositions}.
\smallskip

\pf We claim that $\la (u_n^c, v_n^c, w_n^c,0,0):n\in\nn\ra$ satisfies the hypotheses of Proposition \ref{rel ind join} for the 
joinings $\rho^c_1$ and $\rho^c_2$.

The first two hypotheses follow immediately: $\rho_1^c$ and $\rho_2^c$ are constructed by taking the generic sequences  $\la (u^c_n, v^c_n,0):n\in\nn\ra$ and $\la (v_n^c, w_n^c,0):n\in\nn\ra$ determined by $\la (u_n, v_n,0):n\in\nn\ra$ and $\la (v_n, w_n,0):n\in\nn\ra$
respectively, and the measures did not depend on the precise generic sequence taken.  
Hypothesis \ref{hyp 3} remains to be shown. 

We are given $\epsilon>0$, $k$ and $s^*$ and need to  find $(k')^c, G^c_{(k')^c}$ and the $I_{v^c}$'s so that inequalitites \ref{i} and \ref{ii} hold.  
Since $\rho_1^c$ and $\rho_2^c$ are synchronous, so is the relatively independent joining.  By the preparatory remarks can take $s^*$, the relative location of words in $\bk$ and $\bm$ to be 0. Since the sequence of $(u_n, v_n,w_n,0,0)$'s is generic for the relatively independent product of $\rho_1$ and $\rho_2$, we can find $k', N, G_{k'}\subseteq \mcv_{k'}$ and for each $v\in G_{k'}$ a set $I_v\subset [0, K_{k'})$ such that the conditions in hypothesis \ref{hyp 3} hold in the odometer context.\footnote{For odometer systems, the length of the words in  $\mcu_{k'}, \mcv_{k'}$ and $\mcw_{k'}$ is $K_{k'}$, for circular systems the words at stage $k'$ have length $q_{k'}$.}

Choose $k'$ so large that the density $d_b$ of the boundary portions of circular $k'$-words is less than $\epsilon*10^{-6}$ and so that for each $v\in G_{k'}$, there is an $I_v$ with
\[|I_v|>\left({1-(\epsilon*10^{-6})\over 1-d_b}\right)*K_{k'}.\]

Let $(k')^c=k'$, and $G_{k'}^c=\{v^c:v\in G_{k'}\}$.
For each $v^c\in G_{k'}^c$ we define the set $I_{v^c}\subseteq [0, q_{k'})$. Each $I_v\subseteq[0,K_{k'})$ and  each $s\in I_v$ has a genetic marker $\vec{j}_s$ in $v$. We let $I_{v^c}=\{s^c:s^c$ has the same genetic marker in $v$ as  some $s\in I_v$ does in $v\}$. Equation \ref{up or down} implies that
\[{|I_v|\over K_{k'}}={|I_{v^c}|\over q_{k'}}(1-d_b)\]
and thus $|I_{v^c}|>(1-\epsilon)q_{k'}$.

Equation \ref{all is equal} implies that for $v\in G_{k'}$ and all large $n$, 
\[EmpDist(v_n)(v)=EmpDist(v_n^c)(v^c),\]
from which hypothesis \ref{i} follows immediately. 

Fix a $v_0^c\in G^c_{k'}$ and an 
$s^c\in I_{v_0^c}$. Let $v_0\in G_{k'}$ correspond to $v_0^c$, and $s\in I_v$ correspond to $s^c$. Let $(u^c,w^c)\in \mcu^c_k\times \mcw^c_k$. To see hypothesis \ref{ii}, we need to compute the  empirical distributions of $(u^c,w^c), u^c$ and $w^c$ conditioned on $v_0^c$.

Let $A^c$ be the collection of $((u')^c, v_0^c,(w')^c)\in \mcu_{k'}^c\times \mcv_{k'}^c\times \mcw_{k'}^c$ such that  $u^c$ occurs at $s^c$ in $(u')^c$ and $w^c$ occurs at $s^c$ in $(w')^c$. Let $B^c$ be the collection of all $((u')^c, v_0^c,(w')^c)\in \mcu_{k'}^c\times \mcv_{k'}^c\times \mcw_{k'}^c$. Then:
\begin{equation}\label{upper conditioning}
EmpDist_{k,k,s^c,s^c}(u_n^c, v_n^c, w_n^c|v_0^c)(u^c, w^c)={EmpDist_{k'}(u_n^c, v_n^c, w_n^c)(A^c)\over EmpDist_{k'}(u_n^c, v_n^c, w_n^c)(B^c)}.
\end{equation}
As in the definition of $\mcf$ in Section \ref{f on morphs}, we can view the relatively independent joining $\rho$ on $\bk\times_\bl\bm$ as concentrating on a single odometer system 
$(\bk,\bl,\bm)^\times$ and $\rho^c$, the relatively independent joining of $\rho_1^c, \rho_2^c$ as concentrating on 
$((\bk,\bl,\bm)^\times)^c$, which is canonically isomorphic to $\bk^c\times_{\bl^c}\bm^c$.

In the odometer system $(\bk,\bl,\bm)^\times$, consider the set $A$ consisting of those $k'$-words $(u',v_0,w')$ such that $u'$ and $w'$ have $u$ and $v$ in position $s$.  Then $A^c=\{((u')^c, v_0^c,(w')^c):(u', v_0,w')\in A\}$. Similarly $B^c=\{((u')^c, v_0^c,(w')^c):(u', v_0,w')\in B\}$.
Equation 
\ref{all is equal} implies that
\begin{equation}\label{frontal}
EmpDist(u_n,v_n,w_n)(A)=EmpDist(u_n^c, v_n^c, w_n^c)(A^c).
\end{equation}
and 
\begin{equation}\label{frontal2}
EmpDist(u_n,v_n,w_n)(B)=EmpDist(u_n^c, v_n^c, w_n^c)(B^c).
\end{equation}
Finally noting that
\begin{equation}
\label{lower conditioning}
EmpDist_{k,k,s,s}(u_n,v_n,w_n|v_0)(u,w)={EmpDist_{k'}(u_n,v_n,w_n)(A)\over EmpDist_{k'}(u_n,v_n,w_n)(B)},
\end{equation}
and using equations \ref{upper conditioning} and \ref{frontal} we see that 
\begin{eqnarray}
\label{nothing matters}
EmpDist_{k,k,s^c,s^c}(u_n^c, v_n^c, w_n^c|v_0^c)(u^c, w^c)=\\
EmpDist_{k,k,s,s}(u_n,v_n,w_n|v_0)(u,w).\notag
\end{eqnarray}
Arguing in the same manner we see:
\begin{eqnarray}\label{again again}
EmpDist_{k,s^c}(u_n^c,v_n^c|v_0^c)(u^c)&=&EmpDist_{k,s}(u_n,v_n|v_0)(u)\\
EmpDist_{k,s^c}(v^c_n,w^c_n|v_0^c)(v^c)&=&EmpDist_{k,s}(v_n,w_n|v_0)(v)
\label{again and again}
\end{eqnarray}
Since for large $n$, 
 	\begin{eqnarray}
          &\|EmpDist_{k,k,s,}(u_n, v_n,w_n|v_0)
         -EmpDist_{k,s}(u_n,v_n|v)*EmpDist_{k,s}(v_n,w_n)|v)\|\notag\\
         &<\epsilon,\notag
         \end{eqnarray}
from equations \ref{nothing matters}, \ref{again again} and \ref{again and again} we get the desired conclusion that
 \begin{eqnarray}
          \|EmpDist_{k,k,s^c,s^c}(u^c_n, v^c_n,w^c_n|v^v_0)\ \ \ \ &\notag\\
         -EmpDist_{k,s}(u^c_n,v^c_n|v_0^c)*EmpDist_{k,s}&\!\!\!\!\! (v^c_n,w^c_n|v_0^c)\|\notag
         \end{eqnarray}
 is less than $\epsilon$.
\qed

Lemma \ref{compositions} holds where one or both of the joinings $\rho_1$ and $\rho_2$ are anti-synchronous as well, however the shift coefficients for the circular systems are no longer all $0$ but belong to  $\{0,\pm A_n\}$ depending on which joinings are anti-synchronous. Similarly $s^*\in \{0, \pm A_{k}\}$. The argument follows the same path until it reaches equation \ref{frontal}. This equation relies, in turn on equation \ref{all is equal}. The analogue of equation \ref{all is equal} for anti-synchronous joinings is equation \ref{just another miracle}, which in turn carries over to the relatively independent product. The upshot is that equations \ref{nothing matters}, \ref{again again} and \ref{again and again} hold after applying the appropriate shifts of $u_n^c$ and $v_n^c$ relative to $u_n^c$.

\bigskip
\noindent This finishes the proof of Theorem \ref{grand finale}.\qed

\subsection{Weakly-Mixing  and Compact Extensions}\label{wm and comp}
We now show that $\mcf$ preserves weakly-mixing  and compact extensions.
 The fact that compact extensions are preserved is due to E. Glasner and we reproduce the proof here with his kind permission.
\begin{prop}\label{wm ext} Let $(\bk,\mu)$ and $(\bl,\nu)$ be ergodic and suppose that $\rho$ and $\rho^c$ are corresponding synchronous joinings determining factor maps
\begin{eqnarray*}
\pi&:&\bk\to \bl\\
\pi^c&:&\bk^c\to \bl^c.
\end{eqnarray*}
Then $\bk$ is a weakly mixing extension of $\bl$ (via $\pi$) if and only if $\bk^c$ is a weakly mixing extension of $\bl^c$ (via $\pi^c$).
\end{prop}
\pf Recall that if $\pi:X\to Y$ is a factor map from $(X,\mcb,\mu,T)$ to $(Y,\mcc,\nu,S)$, then the extension is weakly-mixing if the relatively independent joining $X\times_Y X$ of $X$ with itself over $Y$ is ergodic relative to $Y$. In case $Y$ is ergodic, this simply means that the relatively independent joining is ergodic.

Suppose that $\bk$ and $\bl$ are odometer based systems with construction sequences 
$\la \mcw_n:n\in\nn\ra$ and $\la \mcv_n:n\in\nn\ra$ respectively. If $\rho$ is a synchronous factor joining of $\bk$ over $\bl$, and the extension is weakly-mixing then we can find an ergodic sequence of words $\la (u_n, v_n, w_n)\in \mcw_n\times \mcv_n\times \mcw_n:n\in\nn\ra$ that is generic for the relatively independent joining of $\rho$ with itself over $\bl$, i.e. $\rho\times_\bl\rho$. This sequence will satisfy the hypotheses of 
Proposition \ref{rel ind join}. It follows that the sequence of $(u_n^c,v_n^c,w_n^c)$'s is also generic for an ergodic measure $\nu$.  As we argued in Lemma \ref{compositions}, the $(u_n^c, v_n^c, w_n^c)$'s also satisfy the hypothesis of Proposition \ref{rel ind join}. It follows that $\nu$ is the relatively independent joining  $\rho^c\times_\bl\rho^c$.  Since $\nu$ is ergodic $\rho^c$ is weakly mixing. 

If, on the other hand the sequence of $(u_n,v_n,w_n)$ is \emph{not} ergodic, then the sequence $(u_n^c, v^c_n,w^c_n)$ is also not ergodic. Hence if $\rho^c$ is weakly-mixing, then $\rho$ is weakly mixing.\qed

It is immediate from the Furstenberg-Zimmer structure theorem (\cite{glasbook}, Chapter 10, Proposition 10.14)  that $X$ is a 
relatively distal extension of $Y$ if and only if there is no intermediate extension $Z$ of $Y$,  with 
$X$ being  a non-trivial weakly-mixing extension of  $Z$. Thus $\mcf$ takes measure-distal extensions to 
measure-distal extensions.  

What requires more effort to establish is the following:
\begin{prop}(E. Glasner)\label{comp ext}
The functor $\mcf$ takes compact extensions to compact extensions.
\end{prop}
\pf Glasner's proof uses a result proved in the forthcoming \cite{part4}: If $(\bk, \mu)$ is an ergodic odometer based system the $X$ is a compact group extension of $(\bk, \mu)$ then there is a representation of $X$ as an odometer based system with the same coefficients.

Since $X$ is a compact  extension of $Y$ if and only if $X$ is a factor of a compact group extension of $Y$,\footnote{See \cite{Furstenberg-Weiss} for an explicit statement and  proof.}  it suffices to show that $\mcf$ takes compact group extensions to compact group extensions. 

To prove that $\mcf$ takes compact group extensions to compact group extensions we use a remarkable theorem of Veech that characterizes group extensions $\pi:X\to Y$ of ergodic systems. The criteria is that every ergodic joining of $X$ with itself that is the identity on $Y$ (i.e. $\rho$, as a measure, concentrates on those pairs $(x_1, x_2)$ such that $\pi(x_1)=\pi(x_2)$) comes from a graph joining which is an isomorphism of $\xbmt$ that projects to the identity map on $Y$.\footnote{This first appears in \cite{Veech}.} 

Explicity, Theorem 6.18, on page 136 of \cite{glasbook} shows that if, in the ergodic decomposition of the relatively independent product $X\times_Y X$, only graph joinings appear, then $X$ is a compact group extension. The converse follows from Proposition 6.15, part 2 in \cite{glasbook}, that if $X$ is a compact group extension of $Y$ then every ergodic self-joining of $X$ over $Y$ which is the identity on $Y$ is a graph joining.

The map $\mcf$ takes  ergodic joinings to ergodic joinings, and all graph joinings to graph joinings, and the identity joining to the identity joining. Thus we see it preserves group extensions.\qed

Furstenberg \cite{FuBook} and Zimmer \cite{zi} independently showed that for every ergodic system $\mathbb X$ there is an ordinal $\alpha$ and a tower of extensions $\la X_\beta:\beta\le \alpha\ra$ such that $X_0$ is the trivial system, $X_\alpha=X$ and for all $\beta<\alpha$, $X_{\beta+1}$ is a compact extension of $X_\beta$, unless $\alpha=\beta+1$ where $X_\alpha$ is either a compact or a weakly mixing extension of $X_\beta$. If there is no compact extension at the end of the tower, then $X$ is \emph{measure-distal} and $\la X_\beta:\beta<\alpha\ra$ is a \emph{distal tower} approximating $X$.
The least ordinal such that $\mathbb X$ can be represented this way is the \emph{distal height} or \emph{distal order} of $\mathbb X$.

Let $(\bk,\mu)$ be an odometer based system and consider the odometer factor $\mco$. Let $(\bk',\mu')$ be the Kronecker factor of $(\bk,\mu)$. Then we have 
\begin{center}
\begin{equation*}
\begin{diagram}
\node{(\bk,\mu)}\arrow{s,r}{\pi_1}\\
\node{(\bk',\mu')}\arrow{s,r}{\pi_2}\\
\node{\mco}
\end{diagram}
\end{equation*}
\end{center}
where $\pi_2$ may or may not be a trivial factor map. This tower is carried by $\mcf$ to
\begin{center}
\begin{equation*}
\begin{diagram}
\node{(\bk^c,\mu^c)}\arrow{s,r}{\pi_1}\\
\node{((\bk')^c,(\mu')^c)}\arrow{s,r}{\pi_2}\\
\node{\mcr_\alpha}
\end{diagram}
\end{equation*}
\end{center}
If $\bk'$ is a non-trivial extension of $\mco$, then Glasner's result tells us that $(\bk')^c$ is a compact extension of $\mcr_\alpha$, but is silent on the issue of whether $(\bk')^c$ is discrete spectrum; i.e. we do not know  whether $\mcf$ takes the Kronecker factor of $\bk$ to the Kronecker factor of $\bk^c$.

Suppose now that $\bk$ is given by a finite tower of factors:
\begin{equation*}
\begin{diagram}
\node{\mco}\node{\bk_0}\arrow{w}\node{\bk_1}\arrow{w}\node{\dots}\arrow{w}\node{\bk_{N-1}=\bk}\arrow{w}
\end{diagram}
\end{equation*}
where $\bk_0$ is the Kronecker factor of $\bk$ and for all $i, \bk_{i+1}$ is the maximal compact extension of $\bk_i$ in $\bk$. Then $\bk$ is distal of height $N$. The map $\mcf$ carries this to a tower of compact extensions
\begin{equation*}
\begin{diagram}
\node{\mcr_\alpha}\node{\bk^c_0}\arrow{w}\node{\bk^c_1}\arrow{w}\node{\dots}\arrow{w}\node{\bk^c_{N-1}=\bk^c}\arrow{w}
\end{diagram}
\end{equation*}
From this we see that the distal height of $\bk^c$ is either $N$ or $1+N$. 

We do not know an example whether the height of $\bk^c$ can be $1+N$. 
However the ordinary skew product construction applied to odometers gives examples of distal height $n$ where $\mco$ is the Kronecker factor. Hence from our analysis we see that \emph{there are} ergodic circular systems with distal height $N$ for all finite $N$.  

In \cite{BF},  Beleznay and Foreman proved that for all countable ordinals $\alpha$ there is an ergodic measure preserving transformation $T$ of distal height $\alpha$.  In that construction there are no eigenvalues of the operator $U_T$ of finite order. Hence if we let $\mco$ be an  odometer with coefficient sequence $\la k_n:n\in\nn\ra$ going to infinity, $T\times \mco$ is an ergodic transformation with distal height $\alpha$ and zero entropy. In the forthcoming \cite{part4} we see that this implies that $T\times \mco$ can be presented as an odometer based transformation.  By the analysis we just gave we see that $(T\times \mco)^c$ is a circular system with height $1+\alpha$. In \cite{part4} we see that $(T\times \mco)^c$ can be realized as a smooth transformation.  For infinite $\alpha$, $1+\alpha=\alpha$, hence we have:

\begin{theorem}
Let $N$ be a finite or countable ordinal. Then there is an ergodic measure distal diffeomorphism of 
$\bt^2$ of distal height $N$.
\end{theorem}

\subsection{Continuity}
Fix a measure space $(X,\mu)$. As noted in Section \ref{symbolic shifts}, we can identify symbolic shifts built from construction sequences with cut-and-stack constructions (whose levels generate $X$). 
By fixing a countable generating set in advance, we can make this association canonical. 
The levels in the cut-and-stack construction give the relationship with arbitrary partitions of $X$.  In this way the usual weak topology on measure preserving transformation of $X$ described in Section \ref{abstract measure spaces} determines a topology on the presentations of symbolic shifts as limits of construction sequences.

The finitary nature of the maps $\la c_n:n\in\nn\ra$ that give bijections between words in $\mcw_n$ and words in 
$\mcw_n^c$ easily shows that the map $\mcf$ is a \emph{continuous} map from the presentations of odometer based 
systems to presentations of circular systems. Thus we have:

\begin{corollary}\label{continuity}
The functor $\mcf$ is a homeomorphism from the objects in $\mco B$ to $\mcc B$.
\end{corollary}
For the purposes of the complexity of the isomorphism relation we note:
\begin{corollary}
The map $\mcf$ is a continuous reduction of conjugacy between odometer based systems and circular systems.
\end{corollary}
 
\subsection{Extending the main result}
In the main result we restricted the morphisms to graph joinings, largely because compositions of graph joinings are ergodic joinings. Unfortunately a composition of ergodic joinings is not necessarily ergodic, and non-ergodic joinings also arise naturally as relatively independent joinings of ergodic joinings. In this section we indicate how to extend our results to the broader categories that include non-ergodic joinings as morphisms. For convenience, we will continue to require that our objects are ergodic measure preserving systems. 

Let $\mco B^+$ and $\mcc B^+$ be the categories that have the same objects as $\mco B$ and 
$\mcc B$, but where the collections of morphisms are expanded to include \emph{all} synchronous and anti-synchronous joinings (rather than just graph joinings). 

In Section \ref{f on morphs}, the definition  of  $\mcf$ included all such joinings ($\mcf(\rho)$ for a non-ergodic $\rho$ was defined via an ergodic decomposition). Thus without modification we can view $\mcf$ as a map:
\[\mcf:\mco B^+\to \mcc B^+.\]
To show that $\mcf$ is a morphism between these categories, i.e. to show  preserves composition for arbitrary morphisms,  we develop a more combinatorial approach to lifting morphisms that coincides with the original definition.
\medskip

We start by generalizing the notion of a generic sequence of words to include non-ergodic measures. Suppose $\bk$ is a symbolic system with a construction sequence 
$\la \mcw_n:n\in\nn\ra$.  Let $\mu$ be a shift invariant measure which we assume is supported on the set  $S\subseteq\bk$ (where $S$ is given in definition \ref{def of S}). 
The ergodic decomposition theorem gives a representation of $\mu$ as $\int \mu_pd\lambda(p)$, where each $\mu_p$ is a shift invariant ergodic measure and $\lambda$ is a probability measure on a set $P$ parameterizing the ergodic components.  For each 
$p$, there is a generic sequence of words $\la w^p_n:n\in\nn\ra$ for the measure $\mu_p$. 
The main observation is that the set of probability measures on words of a fixed length is compact.
Thus for any fixed $k$ and $\epsilon>0$, we can find a finite set $P_k\subseteq P$ of parameters so that for all $p$, there is some $p'\in P_k$ with\footnote{The notions of $EmpDist$ and $\hat{\mu}_k$ are given  in the beginning of Section \ref{sequences and points}.}
\begin{equation}\label{ep} \|\hat{\mu}^p_k-\hat{\mu}^{p'}_k\|<\epsilon.\end{equation}
This gives a partition of the parameter space into sets $\{E_p:p\in P_k\}$ such that inequality
\ref{ep} holds for all $p'\in E_p$. 

Now let $n$ be sufficiently large such that for each $p\in P_k$, we can find an element $w_n^p\in \mcw_n$ with 
\begin{equation}\label{*}
\|EmpDist_k(w_n^p)-\hat{\mu}_k^p\|<\epsilon.
\end{equation}
If we denote $\lambda(E_p)$ by $\alpha(p)$, then $\alpha(p)\ge 0$ and 
$\sum_{p\in P_k}\alpha_p=1$. It is clear that one can obtain $\hat{\mu}_k$ up to a small error from 
the finite data 
$\{(w_n^p,\alpha(p)):p\in P_k\}$, which is a weighted finite collection of words.

For the symbolic sequences that we are interested in, such as the circular systems, the measure of 
the spacers is independent of the invariant measure $\mu$ (see Section \ref{GMs and coding}). 
This means that for all $n,p$, the sum $\sum_{w'\in\mcw_n}\mu^p_{q_n}(\la w'\ra)$ is the same.
In this context using  inequality \ref{*} we can arrange  the inequality: 
\begin{equation}
\notag
\|(\sum_{p\in P_k}\alpha(p)EmpDist_k(w_n^p))-\hat{\mu}_k\|<\epsilon.
\end{equation}
The measure $\lambda$ is defined on the extreme points of the simplex of shift invariant 
probability measures and if we choose the finite sets $P_k$ to consist of points that lie in the 
closed support of $\lambda$ then we an easily ensure that when we go from $(k,\epsilon)$ to a 
$(k', \epsilon')$ with $k'>k, \epsilon'<\epsilon$ that $P_{k'}\supseteq P_{k}$. Taking a sequence $k\to \infty$ and $\epsilon_k\to 0$ with $\sum\epsilon_k<\infty$, we get a set
 $\{\nu_1, \nu_2, \dots \}$ of ergodic measures  and finite sets $I_k\subseteq I_{k+1}$ of integers with probability measures $\alpha_k$ on $I_k$ such that $(\sum_{i\in I_k}\alpha_k(i)\nu_i)$ converges to  
 $\mu$ in the weak* topology.

\begin{definition}\label{generic non-ergodic} Let $n_k$ go monotonically to infinity and $\{(w^i_{n_k},\alpha_k(i))_k\}$ be a weighted sequence of words as above. 
Suppose that for each $k$ and $i\in I_k$, 
$\|EmpDist_k(w_{n_k}^i)-\hat{\nu}_{i,k}\|<\epsilon_k$, then we call
$\{(w^i_{n_k},\alpha_k(i))\}$ a \emph{generic sequence} for $\mu$. 
\end{definition}
We note that for a fixed $i$, as $k$ varies $\{w^i_{n_k}\}$ is a generic sequence for $\nu_i$--which is one of the ergodic measures in the support of $\lambda$.

In a manner exactly analogous to the analysis in Section \ref{sequences and points}, Definition 
\ref{generic non-ergodic} can be extended to products of symbolic systems, allowing for shifting of words in construction sequences.

Restricting our objects to ergodic systems  $\xbmt, \ycns$ and $(Z,\mcd,\tilde{\mu},\tilde{T})$ allows us to deal with the non-ergodic analogue of the material discussed  between Definition \ref{empdist for joinings}  and Lemma \ref{ill existe} in a relatively straightforward way which we now discuss.

For the analogue of Proposition \ref{rel ind join} in the non-ergodic case let us make the following 
observation.  Fix a non-ergodic joining $\rho$ of $X$ and $Y$ that has ergodic decomposition 
$\rho=\int \rho^pd\lambda(p)$, where, by the ergodicity of $X$ and $Y$, each $\rho^p$ is also a joining of $X$ with $Y$. Fix a $k$ and an $\epsilon>0$ and a cylinder set determined by a word $u\in \mcw^X_k$, at location $s^*$ and let $\phi$ represent its indicator function. For $k'$ large, by the Martingale convergence theorem, there is a subset $G$ of $Y$ of measure close to one such that when we look at the conditional expectation of $\phi$ with respect to the partition induced by the principal $k'$-words of $y\in G$,  for $\mca_{k'}$ and compare it to $\mathbb E(\phi|\mcd)$, the error is small.

The element of that partition that contains $y$ is given by a word $v_y\in \mcw_{k'}^Y$ and a location parameter $s_y$, and the conditional expectation is:
\begin{equation}\label{**}
{\rho(sh^{s^*}(\la u\ra)\cap sh^{s_y}(\la v_y\ra))\over \nu(\la v_y\ra)} 
\end{equation}
This easily gives a set $G_{k'}\subseteq \mcw_{k'}$ with $\hat{\nu}_k(G_{k'})>1-\epsilon$ and a $J_v\subseteq [0,q_{k'})$ such that for $v\in G_{k'}, j\in J_v$, formula \ref{**} gives a good approximation to $\rho_y(sh^{s^*}(\la u\ra))$ for most of the $y\in sh^{s_y}(\la v_y\ra)$.

If we have a generic sequence of weighted words for $\rho$, then we can use it to calculate the expression in \ref{**}. 
This observation makes it possible for us to formulate Proposition \ref{rel ind join} for non-ergodic joinings. 

We are given ergodic systems $X,Y,Z$ and are given construction sequences $\la \mcu_n,\mcv_n,\mcw_n:n\in\nn\ra$ such that for each $n$, the words in each 
$\mcu_n, \mcv_n, \mcw_n$ have the same length. Two joinings $\rho_1$ of $X$ and $Y$ and $\rho_2$ of $Y$ and $Z$ are given. The analogue of Proposition \ref{rel ind join} is now:

\begin{prop} \label{non-erg rel prod}
Let 
\begin{equation}\label{word soup}\la \{(u^i_{n_k},v^i_{n_k},w^i_{n_k},s^i_{n_k},t^i_{n_k}):i\in I_k\},\alpha_k\in Prob(I_k):
k\in\nn\ra
\end{equation}
be a sequence of weighted words and $\sum \epsilon_k<\infty$. Suppose that the following hypothesis are satisfied:
    \begin{enumerate}
    \item  $\la \{(u^i_{n_k},v^i_{n_k},s^i_{n_k}):i\in I_k\}, \alpha_k)\ra_{k}$ is generic for $\rho_1$,
    \item   $\la \{(v^i_{n_k},w^i_{n_k},t^i_{n_k}):i\in I_k\}, \alpha_k\ra_{k}$ is generic for $\rho_2$,
    \item For all $\epsilon, k, s^*$ there are $k', N$ and a set $G_{k'}\subset \mcw_{k'}^Y$ and for each $v\in G_{k'}$ there is a set $J_v\subseteq [0, q_{k'})$ such that 
    	\begin{enumerate}
	\item $\sum_{v\in G_{k'}}EmpDist(v_{n_k})(v)>1-\epsilon$
	\item $|J_v|>(1-\epsilon)q_{k'}$
	\item For all $v\in G_{k'}$ and $s\in J_v$, if $n_k>N$,
	\begin{align*}
\|\sum_{i\in I_k}EmpDist_{k_0,k_0, s,s+s^*}(u^i_{n_k}, sh^{s^i_{k_n}}(v^i_{n_k}),sh^{t^i_{k_n}}(w_{n_k})|v)\alpha_k(i)\ \ \ \ \ \ &-\\\sum_{i\in I_k}EmpDist_{k_0,s}(u^i_{n_k}, sh^{s^i_{k_n}}(v^i_{n_k})|v)\alpha_k(i)* \hskip 1.5in &\\
\sum_{i\in I_k}EmpDist_{k_0,s+s^*}(v^i_{n_k}, sh^{t^i_{n_k}-s^i_{k_n}}(w^i_{n_k})|v)\alpha_k(i)\|&
<\epsilon
	\end{align*}
	\end{enumerate}
    \end{enumerate}
 Then the weighted sequence given in \ref{word soup} is generic for the relatively independent joining $X\times_Y Z$.
\end{prop}

The analogues of Corollary \ref{en fin} and Lemma \ref{ill existe} are easily verified, giving us a characterization of compositions of non-ergodic joinings and the existence of generic sequences satisfying the hypothesis of Proposition \ref{non-erg rel prod}.

Verifying that $\mcf$ preserves composition is now straightforward in the manner of Section \ref{composition}: the $G^c_{k'}$ and $J_{v^c}$ are constructed in exactly the same way.  Checking the conditional distributions of short words relative to longer words ($k$ vs. $k'$) involves counting $k'$-words, and these are counted using Equation \ref{all is equal} for each component $(u_k,v_k, w_k)$ separately. The weighted average is then preserved.

\bigskip

\section{Lagend}
In this section we explore the interplay of the geometric, arithmetic and combinatorial aspects of the manner in which $\mcf$ wraps the odometer based words around the circle. The map $\mcf$  does not preserve the dynamics of the odometer when transforming it into a rotation, indeed it can't. The shift $sh^k$ of the odometer corresponds to a shift $sh^{k^c}$ of the rotation. The relationship between  $k$ and $k^c$ is characterized combinatorially as an optimal  wrapping property. The latter is defined in terms of the notion of a \emph{perfect match}. The results in this section can be used to give an alternate proof of the fact that if $(\bk,\mu)$ is ergodic then so is $(\bk^c,\mu^c)$ that does not use the notion of a generic sequence of words.

Central to our understanding circular systems is the manner in which an $s^c$ had its $n$-words aligned with $n$-words in $sh^k(t^c)$. A word $u$ occurs in $s^c$ lined up with a word $w$ in $sh^k(t^c)$ if and only if $u$ occurs at some location $l$ in $s$ and $w$ occurs at $k+l$ in $t^c$.

\begin{definition}\label{perfect match def}
Let $\vec{x},\vec{y}$ be  strings in the language $\Sigma\cup \{b, e\}$ and $u, v$ be words of the same length. A \emph{$k$-match} of $u$ and $v$ in $\vec{x}$ and $\vec{y}$ is a location $l$ in the domain of $\vec{x}$ such that $u$ occurs at $l$ in $\vec{x}$ and $v$ occurs at $l+k$ in $\vec{y}$. 

If $w^c_0, w^c_1$ are circular $m$-words then a \emph{perfect match} of $u^c,v^c$ in $w^c_0, w^c_1$ is a $k$ such that there are $(n,m)$-genetic markers $\vec{j}_u, \vec{j}_v$ such that $u^c$ occurs in $w^c_0$ and $v^c$ occurs in $w^c_1$ with genetic markers $\vec{j}_u$ and $\vec{j}_v$ respectively and $k$ is a match between all occurrences of 
$u^c$ and $v^c$ with these genetic markers.

\end{definition}
Thus $k$ is a perfect match of $u$ and $v$ if and only if the occurrences of $\vec{j}_u$ in $w_0^c$ are exactly aligned with the occurrences of $\vec{j}_v$ is $w_1^c$.

We will say that \emph{$k$ is a match between $u$ and $v$} if there is a location $l$ such that such that $k$ is a match between $u$ and $v$ at $l$,   
and that \emph{every $k$-match is perfect} when $k$ has the property that  for every  occurrence of a pair of words $u^c,v^c$ in $w^c_0,w^c_1$, if $k$ is a match between $u^c,v^c$ then $k$ is a perfect match between $u^c,v^c$. The astute reader will have already recognized that being a match or a perfect match only refers to the genetic markers and the underlying circular factor--thus the actual identities of $u^c,v^c, w^c_0$ and $w^c_1$ are not material--only the locations of the genetic markers.

The notion of a \emph{perfect match} is vacuous for odometer words; for if  $u,v$ are odometer $n$-words and $w_0, w_1$ are  odometer $m$-words then $u,v$ are the unique pair with a  genetic markers $\vec{j}_u$ and $\vec{j}_v$. Moreover,  if $k$ matches any pair of $n$-subwords,  $k$ matches every pair of corresponding $n$-subwords  in the overlap of $w_0$ and $sh^k(w_1)$. 

If $k>0$, then the $n+1$-subwords of $w_1$ in the overlap of $w_0$ and $sh^k(w_1)$ are split into two pieces by the  $n+1$-subwords of $w_0$; the left portion of each of the $n+1$-subwords of $w_0$ in the overlap coincides with the right portion of the corresponding $n+1$-subword of $w_1$. Call the matches in the left portion of $w_0$  \emph{left-matches}.
\bigskip

\bfni{Discussion.}
Let $u^c$ have genetic marker $\vec{j}_{u^c}=(j_n, j_{n+1},\dots j_{m-1})$ in $w_0^c$ and suppose that $u^c$ sits inside the $n+1$ word $(u')^c$ with genetic marker $( j_{n+1},\dots j_{m-1})$.  Then words with genetic marker $\vec{j}_u$ sit inside every 2-subsection of $u'$.  It follows that if $k^c>0$ and $k^c$ is a perfect match of $u^c$ with $v^c$ having genetic marker $\vec{j}_{v^c}=(j'_n, j'_{n+1},\dots j'_{m-1})$ win $w_1^c$, then $j_n\le j'_n$. Thus the relative position of $v^c$ in the $n+1$-subword of   $w_1^c$ with genetic marker 
$( j'_{n+1},\dots j'_{m-1})$ is to the right of the position of $u^c$ in $(u')^c$; i.e. the relative shift is to the  left to match $u^v$ with $v^c$.  For this reason, when $k^c>0$ we need only consider left shifts.

It is also easy to see that perfect matches between $n$-words  with genetic markers $\vec{j}$ and $\vec{j}'$  inside an $m$-words $w^c_0,w^c_1$ are those $k^c$ that match the first occurrence of an $n$-word with genetic marker $\vec{j}$ in $w^c_0$ with the first occurrence of an $n$-word with genetic marker $\vec{j}'$ in $w^c_1$.

\medskip
The next lemma says that perfect matches can be viewed as the locations of shifts of  odometer based words wrapped around the circle. 
\begin{lemma}\label{lost in translation}
Suppose that $w_0,w_1\in \mcw_m$ and $w_i^c=c_m(w_i)$. Let $n<m$ and $0\le k^c<q_m$ and suppose that $k^c$ is a perfect match between some pair of $n$-subwords of $w_0^c$ and $w_1^c$. Then there is a unique $k$ such that for all  genetic markers $\vec{j}$, $\vec{j}'$, 
\begin{itemize} 
\item $k^c$ is a perfect match between the $n$-subwords of the
$w_i^c$ with genetic markers $\vec{j}$ and $\vec{j}'$ iff 
\item $k$ is a left match between the $n$-subwords of $w$ with genetic markers $\vec{j}$ and $\vec{j}'$. 
\end{itemize}
\end{lemma}
The Lemma has an obvious analogue for negative $k^c$ and right matches.
\medskip

\pf
 Suppose that $k^c$ is a perfect match between $\vec{j}$ and $\vec{j'}$. Call the subwords of $w_0, w_1$ with genetic markers $\vec{j}$ and $\vec{j'}$ $u$ and $v$. Then $u^c, v^c$ are perfectly matched by $k^c$.  Let $k$ be the distance between the locations of $u$ and $v$. Since $k^c\ge 0$ we have $k\ge 0$. 
From our discussion we seen that $k$ is a left match of $u,v$. We claim that this $k$ satisfies the lemma. 
 
Let $u', v'$ be the $n+1$-subwords of $w_0, w_1$ inside which $u,v$ occur. Suppose that $u=u_0u_1\dots u_{k_n-1}$ and $v=v_0v_1\dots v_{k_n-1}$, so  $(u')^c=\mcc((u_0)^c,\dots, (u_{k_n-1})^c)$, $(v')^c=\mcc((v_0)^c, \dots, (v_{k_n-1})^c)$. If $u_i, v_j$ are left matched by $k$ in $u',v'$, then the first occurrences of $(u_i)^c$ and $(v_j)^c$ are matched by $k^c$, hence \emph{inside} $(u')^c, (v')^c$, $k^c$ is a perfect match of $(u_i)^c$ and $(v_j)^c$. 

The relative position of $(u')^c$ and $(v')^c$ is duplicated over all  $n+1$-words  with genetic markers $j_{u'}, j_{v'}$ in $w_0$ and  $sh^k(w_1)$. It follows that $k^c$ is a perfect match of $u^c$ and $v^c$ inside $w_0^c, w_1^c$. 

From the uniformity of the relative positions of $n+1$-words it also follows that any two $n$-subwords of  $n+1$-subwords $(u^*)^c, (v^*)^c$ in positions $i, j$ that are $k^c$ matched in $(u')^c, (v')^c$ are $k^c$-matched. Since these exactly coincide with the $n$-subwords of $w_0, w_1$ that are left-matched by $k$, we have proved the lemma.\qed

\begin{lemma}\label{perfection is possible} Let $w_0=w_1=w$. 
Let $n\in \nn$ and $k\in \mathbb Z$.Then:
\begin{enumerate}
\item {Let $M(k)$ be the least $M$ such that $k<q_M$. Let $u, v\in \mcw_n^c$. Then if $k$ matches $u,v$ inside $w\in \mcw_m^c$ with $m\ge M(k)$ then $u,v$ occur inside the same $M(k)$-subword of $w$.  }               \label{scale of shift}

\item {Let $m\ge M(k)$. Then $k$ is a perfect match of occurrences of $u$, $v$ inside an $m$-word iff $k$ is a perfect match inside the $M(k)$ word in which they appear. \label{new 2}}
\end{enumerate}
\end{lemma}
\pf Use Lemmas \ref{gap calculation} and      \ref{numerology lemma}.
\qed
Item \ref{new 2} means that we usually don't have to refer to a long words when we are discussing perfect matches of $u$ and $v$ and fixes the scale of the potential perfect matches.

We can identify perfect matches numerically:

\begin{lemma}\label{small nudge} Let $A<q_N$, $w^c\in \mcw_N^c$. Then there is a $(m,n)$-genetic marker $\vec{j}$ such that $A$ is the location of the first occurrence of some word 
genetic marker $\vec{j}$ if and only if

	\begin{equation}
	A=c_{N-1}l_{N-1}q_{N-1}+c_{N-2}l_{N-1}q_{N-2} + \dots +c_{m}l_mq_m \label{very pretty}
	\end{equation}

where $0\le c_i<k_i$. 
\end{lemma}

From Lemma \ref{lost in translation}, we see the correspondence between odometer translations and circular translation. We now address the question: given an arbitrary circular translation, how does one adjust it to get an odometer translation that gives the best fit among a given collection of $n$-words?

\begin{theorem}\label{perfection is sorta possible}
Let $n\in \nn$ and $k\in \mathbb Z$. Then 
if $\{(u_i, v_i):i\in I\}\subseteq \mcw_n^c \times \mcw_n^c$, $w\in \mcw_m^c$ ($m\ge M(k)$) then there is a $k'$ such that $|k'-k|<q_{m}$ and: 
	\begin{enumerate}
	\item all $k'$-matches of a $(u_i,v_i)$ in $w$ are perfect matches,\label{perfection}
	\item and 
	\[\sum_i |k'\mbox{-matches of a $u_i$ with a $v_i | \ge \sum_i |k$-matches of a $u_i$ with a $v_i|$}\]
	\end{enumerate}
\end{theorem}
\pf Without loss of generality $k\ge 0$ (otherwise we reverse the role of $u$ and $v$). 
{Words in $\mcw_m$ with $m>M(k)$ start with a block of $b$'s of length at least $q_{M(k)}$. Hence if $k$ matches $n$-words $u, v$ inside $w\in \mcw^c_m$, they both must occur in some $M(k)$-subword of $w$.}

To see item \ref{perfection}, we need to show how to improve $k$ to a $k'$ that is a perfect match. Changing $k$ will involve sacrificing some of the matches of pairs in $I$, but this will be compensated by the additional multiplicity of the remaining matches.

{We prove by induction on $d\ge 1$, that for all $m, n$ with $m-n=d$ and all collections of pairs of $m$-words 
$\{(w_0^j,w_1^j):j\in J\}$ and all $k$, all $m\ge M(k)$, all natural number weightings $\{\alpha_j:j\in J\}$ and all 
$\{(u_i,v_i):i\in I\}\subseteq \mcw_n^c\times \mcw_n^c$} we can find a $k'$ such that $|k'-k|<q_m$ such that (a) holds and 
\begin{eqnarray*}
\sum_j\sum_i \alpha_j|\{k'\mbox{-matches of a $u_i$ and a $v_i$ in  $(w_0^j,w_1^j)$}\}|&\ge\\
\sum_j\sum_i \alpha_j|\{k\mbox{-matches of a $u_i$ and a $v_i$ in  $(w_0^j,w_1^j)$}\}|
\end{eqnarray*} 

Suppose first that $d=1$. Then successive $2$-subsections of $m$-words 
 are separated by boundary sections of size 
\[j_i+(q-j_{i+1})\equiv p_n^{-1} \mbox{(mod $q$)}.\]
Because $q_n$ does not divide $p_n$, given a $2$-subsection $\vec{s}$ of $w_0^j$ there is a unique $2$-subsection $\vec{t}$ of $w^j_1$ within which $k$ can match $n$-words. Moreover this does not depend on $j$, but rather the underlying locations of the words.

We start by lining up blocks of the form $u_i^{l_n-1}$ with blocks of the form $v_i^{l_n-1}$. To do this we classify the $k$-matches of a pair $(u,v)=(u_i,v_i)$ into \emph{left block matches} if $u$ and $sh^k(v)$ align as\footnote{In both of these graphics the second row is a portion of $sk^k(w^j_1)$ and $B$ represents a boundary section. These pictures are independent of $j$.}:
\begin{center}
\includegraphics[width=.9\textwidth]{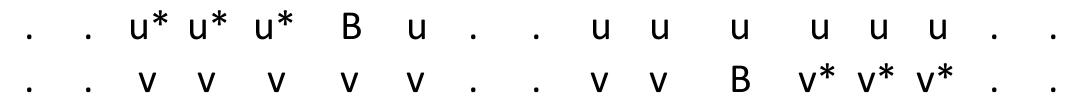}
\end{center}
and  \emph{right block matches} if $u$ and $sh^k(v)$ align as:
\begin{center}
\includegraphics[width=.9\textwidth]{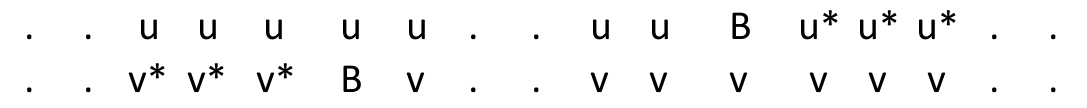}
\end{center}
Note by taking $k'$ to be $k+lq_n$ for some $l<l_n-1$ we can turn all left block matches of all {of the }$(u_i,v_i)$ into matches of entire $u_i^{l_n-1}$ with $sh^k(v_i^{l_n-1})$, but doing so destroys completely some of the right block matches. Similarly if we can shift to make all right block matches into matches of $u_i^{l_n-1}$ with $sh^k(v_i^{l_n-1})$ by destroying left block matches.

If we examine a particular left block match of a  pair $(u_i,v_i)$ in some $w_0^j$ and a right block match of another pair $(u'_i,v'_i)$ in $w_1^j$ and we change $k$ to $k'$ to make $u_i^{l_n-1}$ match with $sh^{k'}(v_i^{l_n-1})$ then the sum of $k'$-matches between $(u_i,v_i)$ and $(u_i',v_i')$ goes up by one: we lose the right block matches but we gain left block matches and we gain one more match from the boundary section.

Suppose that 
\begin{eqnarray*}
\sum_j\sum_j \alpha_j |\{\mbox{left block matches in }(w_0^j,w_1^j)|&\ge\\
\sum_j\sum_j \alpha_j |\{\mbox{right block matches in }(w_0^j,w_1^j)|
\end{eqnarray*}
 Then from the previous paragraph that if we take $k'=k+lq_n$ for some $l<l_n-1$ then we can make all left block matches have multiplicity $l_{n}-1$ (while removing right block matches) and have:
\begin{eqnarray*}\sum_j\sum_i \alpha_j |k'\mbox{-matches of a $u_i$ with a $v_i $ in some $(w_0^j,w_1^j)|$} &\ge\\ 
\sum_j\sum_i \alpha_j|k\mbox{-matches of a $u_i$ with a $v_i$ in some $(w_0^j,w_1^j)|$}.
\end{eqnarray*}
If, on the other hand, the weighted sum of the right block matches is greater than weighted sum of the left block matches, we shift the other direction to fix all right block matches and destroy all left block matches.

Thus we can assume that we have a $k$ such that for all $(u_i,v_i)$, $sh^k$ matches $(l_n-1)$-powers  of 
$u_i$ with $(l_n-1)$-powers of $v_i$. This $k$ would be a perfect match except that it matches $n$-words 
across 2-subsections. Writing each  $w_s^j=\mcc_n(w_1, \dots w_{k_n})$ then $sh^{k}$ matches blocks of 
the form $w_s^{l_n-1}$ in one $1$-subsection of $w_0^j$ with a block of the form $w_{s'}^{l_n-1}$ in a 
(potentially) different $1$-subsection of $w_1^j$. Moreover $s-s'$ is constant on all of these matches, since 
the differences between starts of $w_j^{l_n-1}$-blocks are of length $l_nq_n$. Fix such a pair $s, s'$. By 
changing $k$ so that it lines up $w_s^{l_n-1}$ with $w_{s'}^{l_n-1}$ in the first 1-subsection we create a 
perfect match of $n=m-1$-words and increase the total number of matches of the form $(u_i,v_i)$.
This establishes the case where $d=1$.

We now do the induction step. Let $d=m-n$ and assume the result holds for $d-1$. Suppose that we are given $\{\alpha_j:j\in J\}$. 

We can decompose a $k$-match between $n$-subwords of $w^j_0$ and $w^j_1$ as  $k_1+k^*$ where $k^*\in [-q_{m-1}+1, q_{m-1}-1]$ and $k_1$ is a match of $m-1$ subwords of $w^j_0$ and $w^j_1$.

Here is a picture of a pair $(u',v')\in \mcw^c_{m-1}\times \mcw^c_{m-1}$ comparing $w_0$ in the upper row with the $k_1$-shift of $w_1$ in the lower row.

\begin{center}
\includegraphics[width=.9\textwidth]{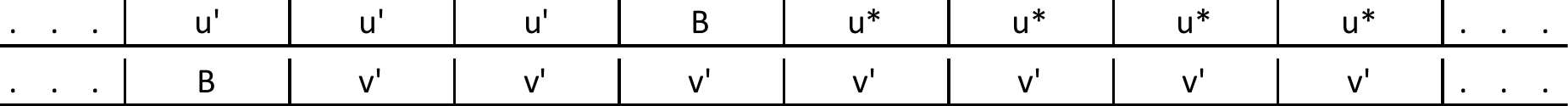}
\end{center}

Here is a picture after the $k=k_1+k^*$ shift of $w_1$:

\begin{center}
\includegraphics[width=.9\textwidth]{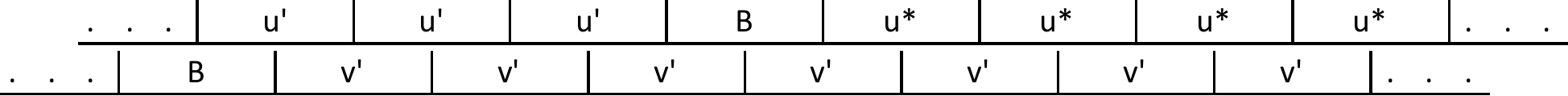}
\end{center}

Let $\{(u',v')_{i'}:i'\in I'\}$ be the collection of pairs $(u',v')$ from $\mcw^c_{m-1}$ sitting inside a pair $(w_0^j, w_1^j)$ that contain $k$-matches of words $(u_i,v_i)$. Arguing as in the case $d=1$ we can adjust $k_1$ to a $k_1'$ so that it is a perfect match of $m-1$-words in $I'$ and, summing over $I$ and $J$, the weighted sum of $k_1'+k^*$-matches 
of pairs in  $I$ does not decrease.\footnote{We note that it is not enough to increase the weighted sum of the number of matches of pairs in $I'$, because various  $I'$ matches may contain different number of $I$-matches. Nonetheless, arguing as in the case $d-1$, one of the two possibilities for lining up the $m-1$ subwords does not decrease the weighted sum of the number of $k_1'+k^*$-matches of $I$-words.}

This is how the $m-1$-words look after shifting by $k_1'+k^*$:
\begin{center}
\includegraphics[width=.9\textwidth]{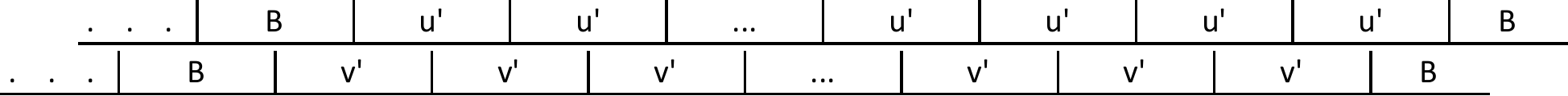}
\end{center}
The offset of the copies of $u'$ and $v'$ is $k^*$. Note that the boundary sections line up.

We now are in the position of having shifted by $k_1'$ so that the powers of pairs $\{(u',v')_{i'}:{i'\in I'}\}$ are lined up. The additional shift $k^*$ has absolute value less than $q_{m-1}$. Moreover all of the words $\{(u',v')_{i'}:{i'\in I'}\}$ are lined up the same way when shifted by $k^*$. 

We call an occurrence of a $(u',v')_{i'}$ that is lined up in $(w_0^j, sh^{k'_1}(w_1^j))$ \emph{good}. Let $\beta_{j,i'}$ be the number of good occurrences of $(u',v')_{i'}$ and
\[\alpha'_{i'}=\sum_j \alpha_j\beta_{j,i'}.\]
Note that 
\begin{eqnarray*}\sum_{i'}\sum_i \alpha_{i'}|(k_1'+k^*)\mbox{-matches of a $u_i$ with a $v_i $ in some good occurrence of $(u',v')_{i'}|$}\\ 
=\sum_j\sum_i \alpha_j |(k_1'+k^*)\mbox{-matches of a $u_i$ with a $v_i$ in some $(w_0^j,w_1^j)|$}.
\end{eqnarray*}

We now view the pairs $\{(u',v')_{i'}:{i'\in I'}\}$ as sitting on the intervals $[0,q_{m-1}-1]$ and then shifting $v'$  by $k^*$:
\begin{center}
\includegraphics[width=.2\textwidth]{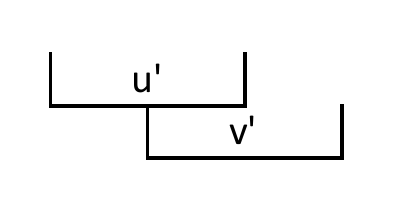}
\end{center}
We are in a position to apply our induction hypothesis with $I'$ playing the role of $J$, the $\alpha'_{i'}$'s being the $\alpha_j$'s, $d-1=(m-1)-n$ and the shift being $k^*$. 

The result is a $k^{**}$ such that every $k^{**}$-match of a $(u_i,v_i)$ in a $(u',v')_{i'}$ is perfect and 
\begin{eqnarray*}\sum_{i'}\sum_i\alpha_{i'} |(k^{**})\mbox{-matches of a $u_i$ with a $v_i $ in some $(u',v')_{i'}|$} &\ge\\ 
\sum_{i'}\sum_i \alpha_{i'}|(k^*)\mbox{-matches of a $u_i$ with a $v_i$ in some $(u',v')_{i'}|$}.
\end{eqnarray*}
We note that every $k_1'+k^{**}$-match of a $(u_i,v_i)$ in a $(w_0^j,w_1^j)$ is perfect. Since
\begin{eqnarray*}\sum_{i'}\sum_i \alpha_{i'}|(k_1'+k^{**})\mbox{-matches of a $u_i$ with a $v_i $ in some good occurrence of $(u',v')_{i'}|$}\\ 
=\sum_j\sum_i \alpha_j |(k_1'+k^{**})\mbox{-matches of a $u_i$ with a $v_i$ in some $(w_0^j,w_1^j)|$}.
\end{eqnarray*}
we see that
\begin{eqnarray*}\sum_j\sum_i \alpha_j|(k_1'+k^{**})\mbox{-matches of a $u_i$ with a $v_i $ in some $(w_0^j,w_1^j)|$} &\ge\\ 
\sum_j\sum_i \alpha_j |k\mbox{-matches of a $u_i$ with a $v_i$ in some $(w_0^j,w_1^j)|$}.
\end{eqnarray*}
This completes the proof of Lemma \ref{perfection is sorta possible}.\qed

\section{Open Problems}
We finish with two open problems that we find interesting and believe to be feasible. The first is to characterize the class of transformations isomorphic to circular systems in  Ergodic-theoretic terms. All circular systems have common properties such that can be described in terms of rigidity sequences or zero entropy. The suggestions is to find a complete characterization in using this type of notion. 

The second problem can be stated as follows. For the realization problem, the underlying rotation $\alpha$ of a circular system must be Liouvillian; however realization is   not necessary for the results in this paper.  Can an arbitrary irrational $\alpha$ be the underlying rotation of a circular system?

\bibliography{citations}

\begin{thebibliography}{10}

\bibitem{AK_original}
D.~V. Anosov and A.~B. Katok.
\newblock New examples in smooth ergodic theory. {E}rgodic diffeomorphisms.
\newblock {\em Trudy Moskov. Mat. Ob\v s\v c.}, 23:3--36, 1970.

\bibitem{BF}
Ferenc Beleznay and Matthew Foreman.
\newblock The complexity of the collection of measure-distal transformations.
\newblock {\em Ergodic Theory Dynam. Systems}, 16(5):929--962, 1996.

\bibitem{downar}
Tomasz Downarowicz.
\newblock The {C}hoquet simplex of invariant measures for minimal flows.
\newblock {\em Israel J. Math.}, 74(2-3):241--256, 1991.

\bibitem{feldman}
Jacob Feldman.
\newblock Borel structures and invariants for measurable transformations.
\newblock {\em Proc. Amer. Math. Soc.}, 46:383--394, 1974.

\bibitem{prequel}
M.~Foreman and B.~Weiss.
\newblock A symbolic representation of {A}nosov-{K}atok systems.
\newblock {\em To appear in Journal d'Analyse Mathematique}, pages 1--63, 2015.

\bibitem{FRW}
Matthew Foreman, Daniel~J. Rudolph, and Benjamin Weiss.
\newblock The conjugacy problem in ergodic theory.
\newblock {\em Ann. of Math. (2)}, 173(3):1529--1586, 2011.

\bibitem{part3}
Matthew Foreman and Benjamin Weiss.
\newblock Measure preserving diffeomorphisms of the torus are unclassifiable.
\newblock {\em TO APPEAR}, pages 1--102, 2018.

\bibitem{part4}
Matthew Foreman and Benjamin Weiss.
\newblock Odometer based systems.
\newblock {\em TO APPEAR}, 2019.

\bibitem{FuBook}
H.~Furstenberg.
\newblock {\em Recurrence in ergodic theory and combinatorial number theory}.
\newblock Princeton University Press, Princeton, N.J., 1981.
\newblock M. B. Porter Lectures.

\bibitem{Furstenberg-Weiss}
Hillel Furstenberg and Benjamin Weiss.
\newblock A mean ergodic theorem for {$(1/N)\sum^N_{n=1}f(T^nx)g(T^{n^2}x)$}.
\newblock In {\em Convergence in ergodic theory and probability ({C}olumbus,
  {OH}, 1993)}, volume~5 of {\em Ohio State Univ. Math. Res. Inst. Publ.},
  pages 193--227. de Gruyter, Berlin, 1996.

\bibitem{glasbook}
Eli Glasner.
\newblock {\em Ergodic theory via joinings}, volume 101 of {\em Mathematical
  Surveys and Monographs}.
\newblock American Mathematical Society, Providence, RI, 2003.

\bibitem{halmos}
Paul~R. Halmos.
\newblock {\em Lectures on ergodic theory}.
\newblock Chelsea Publishing Co., New York, 1960.

\bibitem{HvN}
Paul~R. Halmos and John von Neumann.
\newblock Operator methods in classical mechanics. {II}.
\newblock {\em Ann. of Math. (2)}, 43:332--350, 1942.

\bibitem{katoksbook}
Anatole Katok.
\newblock {\em Combinatorial constructions in ergodic theory and dynamics},
  volume~30 of {\em University Lecture Series}.
\newblock American Mathematical Society, Providence, RI, 2003.

\bibitem{Peterson}
Karl Petersen.
\newblock {\em Ergodic theory}, volume~2 of {\em Cambridge Studies in Advanced
  Mathematics}.
\newblock Cambridge University Press, Cambridge, 1989.
\newblock Corrected reprint of the 1983 original.

\bibitem{DansBook}
Daniel~J. Rudolph.
\newblock {\em Fundamentals of measurable dynamics}.
\newblock Oxford Science Publications. The Clarendon Press, Oxford University
  Press, New York, 1990.
\newblock Ergodic theory on Lebesgue spaces.

\bibitem{Shields}
Paul~C. Shields.
\newblock {\em The ergodic theory of discrete sample paths}, volume~13 of {\em
  Graduate Studies in Mathematics}.
\newblock American Mathematical Society, Providence, RI, 1996.

\bibitem{Veech}
William~A. Veech.
\newblock A criterion for a process to be prime.
\newblock {\em Monatsh. Math.}, 94(4):335--341, 1982.

\bibitem{vN}
J.~von Neumann.
\newblock Zur {O}peratorenmethode in der klassischen {M}echanik.
\newblock {\em Ann. of Math. (2)}, 33(3):587--642, 1932.

\bibitem{walters}
Peter Walters.
\newblock {\em An introduction to ergodic theory}, volume~79 of {\em Graduate
  Texts in Mathematics}.
\newblock Springer-Verlag, New York-Berlin, 1982.

\bibitem{Benjy}
Benjamin Weiss.
\newblock {\em Single orbit dynamics}, volume~95 of {\em CBMS Regional
  Conference Series in Mathematics}.
\newblock American Mathematical Society, Providence, RI, 2000.

\bibitem{zi}
Robert~J. Zimmer.
\newblock Ergodic actions with generalized discrete spectrum.
\newblock {\em Illinois J. Math.}, 20(4):555--588, 1976.

\end{thebibliography}
\bibliographystyle{plain}
\end{document}